%% file: 2003-19.tex
\input gtmacros

\input amsnames

\input amstex
\let\cal\Cal          
\catcode`\@=12        
%
%
\input gtoutput

\lognumber{291}
\volumenumber{7}\papernumber{19}\volumeyear{2003}
\pagenumbers{645}{711}
\received{8 November 2002}
\published{28 October 2003}
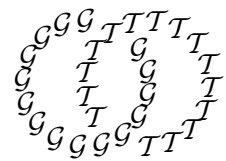
\accepted{20 October 2003}
\proposed{Haynes Miller}
\seconded{Ralph Cohen, Gunnar Carlsson}

\let\\\par
\def\topmatter{\relax}
\def\endtopmatter{\maketitlepage}
\let\gttitle\title
\def\title#1\endtitle{\gttitle{#1}}
\let\gtauthor\author
\def\author#1\endauthor{\gtauthor{#1}}
\let\gtaddress\address
\def\address#1\endaddress{\gtaddress{#1}}
\def\affil#1\endaffil{\gtaddress{#1}}
\let\gtemail\email
\def\email#1\endemail{\gtemail{#1}}
\def\subjclass#1\endsubjclass{\primaryclass{#1}}
\let\gtkeywords\keywords
\def\keywords#1\endkeywords{\gtkeywords{#1}}
\def\heading#1\endheading{{\def\S##1##2{\relax}\def\\{\relax\ignorespaces}
    \section{#1}}}
\def\head#1\endhead{\heading#1\endheading}
\def\headnonumber#1\endheadnonumber{
\vskip-\lastskip\penalty-800\vskip 20pt plus10pt minus5pt 
{\large\bf#1}         
\vskip 8pt plus4pt minus4pt
\nobreak\resultnumber=1}

\def\subhead#1\endsubhead{\sh{#1}}
\def\subsubhead#1\endsubsubhead{\sh{#1}}
\def\specialhead#1\endspecialhead{\sh{#1}}
\def\demo#1{\rk{#1}\ignorespaces}
\def\enddemo{\ppar}
\let\remark\demo
\def\endremark{}
\let\definition\demo
\def\enddefinition{\ppar}
\let\example\demo
\def\endexample{\ppar}
\def\qed{\ifmmode\quad\sq\else\hbox{}\hfill$\sq$\par\goodbreak\rm\fi}  
\def\proclaim#1{\rk{#1}\sl\ignorespaces}
\def\endproclaim{\rm\ppar}
\def\cite#1{[#1]}
\newcount\itemnumber
\def\roster{\items\itemnumber=1}
\def\endroster{\enditems}
\let\itemold\item
\def\item{\itemold{{\rm(\number\itemnumber)}}%
\global\advance\itemnumber by 1\ignorespaces}
\def\S{Section~\ignorespaces}  
\def\date#1\enddate{\relax}
\def\thanks#1\endthanks{\relax}   
\def\dedicatory#1\enddedicatory{\relax}  
\let\footnote\plainfootnote

\def\Refs{\ppar{\large\bf References}\ppar\bgroup\leftskip=25pt
\frenchspacing\parskip=3pt plus2pt\small}       
\def\endRefs{\egroup}
\def\widestnumber#1#2{\relax}
\def\endrefitem{}
\def\refdef#1#2#3{\def#1{\leavevmode\unskip\endrefitem#2\def\endrefitem{#3}}}
\def\ref{\par}
\def\endref{\endrefitem\par\def\endrefitem{}}
\refdef\key{\noindent\llap\bgroup[}{]\ \ \egroup}
\refdef\no{\noindent\llap\bgroup[}{]\ \ \egroup}
\refdef\by{\bf}{\rm, }
\refdef\manyby{\bf}{\rm, }
\refdef\paper{\it}{\rm, }
\refdef\book{\it}{\rm, }
\refdef\jour{}{ }
\refdef\vol{}{ }
\refdef\yr{(}{) }
\refdef\ed{(}{, Editor) }
\refdef\publ{}{ }
\refdef\inbook{from: ``}{'', }
\refdef\pages{}{ }
\refdef\page{}{ }
\refdef\paperinfo{}{ }
\refdef\bookinfo{}{ }
\refdef\publaddr{}{ }
\refdef\moreref{}{ }
\refdef\eds{(}{, Editors)}
\refdef\bysame{\hbox to 3 em{\hrulefill}\thinspace,}{ }
\refdef\toappear{(to appear)}{ }
\refdef\issue{no.\ }{ }
\newcount\refnumber\refnumber=1
\def\refkey#1{\expandafter\xdef\csname cite#1\endcsname{\number\refnumber}%
\global\advance\refnumber by 1}
\def\cite#1{[\csname cite#1\endcsname]}
\def\Cite#1{\csname cite#1\endcsname}  
\def\no#1{\noindent\llap{[\csname cite#1\endcsname]\ \ }}

\refkey4 
\refkey{13} 
\refkey{14} 
\refkey6 
\refkey3 
\refkey{12} 
\refkey{11} 
\refkey{15} 
\refkey1 
\refkey2 
\refkey5 
\refkey{17} 
\refkey{16} 
\refkey7 
\refkey9 
\refkey{10} 
\refkey8 

\define\a{\rightarrow}
\define\lra{\longrightarrow}

\define\st{^{\text{st}}}
\define\holim{\operatornamewithlimits{holim}}
\define\colim{\operatornamewithlimits{colim}}
\define\hocolim{\operatornamewithlimits{hocolim}}
\define\hofiber{\operatornamewithlimits{hofiber}}
\define\hocofiber{\operatornamewithlimits{hocofiber}}
\define\tfiber{\operatornamewithlimits{tfiber}}
\define\fiber{\operatornamewithlimits{fiber}}

\topmatter
\title Calculus III: Taylor Series\endtitle
\author Thomas G Goodwillie\endauthor
\affil Department of Mathematics, Brown 
University\\Box 1917, Providence RI 02912--0001,
USA\endaffil
\email tomg@math.brown.edu\endemail
\keywords Homotopy functor, excision, Taylor tower\endkeywords
\subjclass 55P99\endsubjclass\secondaryclass{55U99}
\abstract We study functors
from spaces to spaces or spectra that preserve weak homotopy
equivalences. For each such functor we construct a universal $n$-excisive
approximation, which may be thought of as its $n$-excisive part.
Homogeneous functors, meaning $n$-excisive functors with trivial
$(n-1)$-excisive part, can be classified: they correspond to symmetric
functors of $n$ variables that are reduced and $1$-excisive in each variable.
We discuss some important examples, including the identity functor and
Waldhausen's algebraic
$K$-theory.
\endabstract

\asciiabstract{%
We study functors from spaces to spaces or spectra that preserve weak
homotopy equivalences. For each such functor we construct a universal
n-excisive approximation, which may be thought of as its n-excisive
part.  Homogeneous functors, meaning n-excisive functors with trivial
(n-1)-excisive part, can be classified: they correspond to symmetric
functors of n variables that are reduced and 1-excisive in each
variable.  We discuss some important examples, including the identity
functor and Waldhausen's algebraic K-theory.}

\endtopmatter

\catcode`\@=\active

\headnonumber Introduction\endheadnonumber

This paper should have been finished many years ago. It is a continuation of
\cite{1} and \cite{2} (which were themselves a little late in coming). The
author has no excuse for this procrastination and wishes to apologize to
anyone who has been inconvenienced, especially his students and former
students.

As in \cite{1} and \cite{2} we are concerned with functors from $\Cal C$, which
may be the category of spaces or more generally spaces over a fixed space, to
$\Cal D$, which may be either the category of spaces or the category 
of spectra.
It is clear that one could extend the ideas to many other settings. The
functors will always be
\it homotopy functors\rm\ in the sense that they preserve weak homotopy
equivalences. We are also concerned with (natural) maps between such 
functors. We
call such a map of functors
$F@>>>G$ a weak equivalence (or for emphasis an objectwise weak equivalence, or
for brevity an equivalence) if for every object
$X$ the map
$F(X)@>>>G(X)$ is a weak equivalence. The goal is to shed light on homotopy
functors in general (and on particular ones) by systematically comparing them
with homotopy functors of some very special kinds, namely those which satisfy
\lq\lq$n^{th}$ order excision\rq\rq.

In taking this point of view we are led naturally to speak of the 
category $\Cal F
=\Cal F(\Cal C,\Cal D)$ of homotopy functors. Some of our conclusions 
refer to the homotopy category $h\Cal F$ of homotopy functors,
meaning the universal example of a category equipped with a functor 
from $\Cal F$
that takes every weak equivalence to an isomorphism. For example, some of our
results are most neatly expressed as statements to the effect that 
such and such
a functor from one category of functors to another induces an equivalence of
homotopy categories. Homotopy categories in this sense can be
made by the \lq\lq generators and relations\rq\rq\ construction:
a morphism in the homotopy category is an equivalence class of words, 
where each
letter in a word is either a forward arrow (a morphism in
$\Cal F$) or a backward arrow  (a formal inverse to a weak 
equivalence in $\Cal F$)
and two words are declared to represent the same morphism if they are 
related by
a sequence of basic moves (composing two forward or two backward arrows,
deleting an identity arrow, cancelling a forward arrow and its formal inverse).

The author has recently been reminded that there are set-theoretic 
objections to
bandying this kind of language about too freely, but he does not want 
to stop doing
so. Some ways of resolving this
difficulty will be discussed at the end of this introduction.

Of course, if one wants to prove theorems about a homotopy category then one
should be prepared to do most of the work in the original category. 
For example,
commutative diagrams in the category of functors will be a ubiquitous tool 
here,
but commutative diagrams in the homotopy category are relatively useless.

We now recall the definition of $n^{th}$ order excision, from section 3
of \cite{2}. Whereas ordinary, or first order, excision has to do 
with the behavior
of a functor on certain square diagrams,
$n^{th}$ order excision has to do with the behavior of a functor on certain
cubical diagrams of dimension
$n+1$.  We briefly recall the definitions. The reader can consult 
the opening sections of \cite{2} for details, including the specific 
models we are using for homotopy limits and colimits and the basic 
definitions and facts about cubical diagrams.

The homotopy functor $F$ is \it excisive\rm\  if it takes
homotopy pushout squares (also called homotopy cocartesian squares) to
homotopy pullback squares (also called homotopy cartesian squares). It is
\it reduced\rm\ if it takes the final object to an object weakly 
equivalent to the
final object. It is
\it linear\rm\ if it is both excisive and reduced.

A cubical diagram is called strongly homotopy cocartesian if all of its
two-dimensional faces are homotopy pushouts. It is called homotopy cartesian
if it is a homotopy pullback in the sense that the \lq\lq 
first\rq\rq\ object is
equivalent (by the obvious map) to the homotopy limit of all the 
others.  In this
paper we have generally omitted the word \lq\lq
homotopy\rq\rq\ in the expressions \lq\lq (strongly) homotopy
(co)cartesian\rq\rq.
$F$ is
said to be
$n$-\it excisive\rm , or to satisfy $n^{th}$ order excision, if it 
takes strongly
homotopy cocartesian
$(n+1)$-cubical diagrams to homotopy cartesian diagrams. An $n$-excisive
functor is always
$(n+1)$-excisive.

In \cite{1} we studied the approximation of homotopy functors by excisive
functors and codified this in the notion of (first) derivative of a
functor. We calculated the derivative, first in some basic examples and
then in the example which had given rise to the whole project: stable
pseudoisotopy theory (and with it Waldhausen's algebraic $K$-theory functor
$A$).

In \cite{2} we introduced the class of
\it analytic\rm\ functors. These are homotopy functors whose deviation from
being $n$-excisive is bounded in a certain way for all
$n$. Some functors are more analytic
than others; an analytic functor is $\rho$-analytic for some integer 
$\rho$, where
a smaller $\rho$ means stronger bounds. We showed that
$\rho$-analytic functors behave rather rigidly on the category of
$\rho$-connected spaces (and more generally on the category of spaces
equipped with
$(\rho +1)$-connected
maps to a fixed space), in the sense, roughly, that if a
map between two such functors induces an equivalence of first derivatives
everywhere then its homotopy fiber is constant up to homotopy within 
the category of
$\rho$-connected spaces.

In \cite{3} these results were used to relate $A$ to some functors built
out of the free loopspace.

Here we move from excisive approximation to
$n$-excisive approximation, obtaining functors $P_nF$ which can be thought of
as Taylor polynomials of $F$. We show that if $F$ is $\rho$-analytic 
and  $X$ is
$\rho$-connected then as $n$ tends to infinity the $n^{th}$ approximation
converges to $F(X)$ in the strong sense that the connectivity of a certain
map $F(X)@>>>(P_nF)(X)$ tends to
infinity. In other words, the number $\rho$ gives a sort of radius of
convergence.

What we were doing in \cite{2} was like showing that a function
$f(x)$ is determined, in some region, by
$f(0)$ and $f'(x)$. Continuing with the same analogy, what we are 
doing here is like
showing that $f(x)$ is determined by $f(0)$, $f'(0)$, $f''(0)$ and so on.

Here is a sketch of the main results of this paper, presented in a rather
different order from the one in which they will be proved.

First consider homotopy functors $\Cal T@>F>>\Cal Sp$ from based
spaces to spectra.
For any such functor, and for any
$n\geq 0$, we can make an
$n$-excisive functor $\Cal T@>{P_nF}>>\Cal Sp$ with a map $F@>>>P_nF$ that
is universal (in a homotopy category) among maps from
$F$ to $n$-excisive functors. The functors $\lbrace P_nF\rbrace_{n\geq 0}$ fit
together to form a tower, and $F$ maps into the limit of the tower.

If the role of
$n^{th}$ Taylor polynomial is being played by
$P_nF$, then the role of
$n^{th}$ term in the series is played by the (homotopy) fiber of the map
$P_nF@>>>P_{n-1}F$, which will be denoted by $D_nF$. We refer to $P_nF$ as the
$n^{th}$ \it stage\rm\ of the tower and to $D_nF$ as the $n^{th}$
\it layer\rm\. These \it homogeneous
polynomial functors\rm\ are the focus of much of the paper.

The constant term $(P_0F)(X)$ is the same, up to natural
weak equivalence, as the constant functor
$F(*)$.

Up to the same kind of equivalence, the linear (or homogeneous of degree one)
functor
$(D_1F)(X)$ necessarily has the form
$C_1\wedge X$ for some fixed spectrum $C_1$, at least
when restricted to finite complexes $X$. The coefficient spectrum $C_1$,
which of course is
$(D_1F)(S^0)$, is called the derivative of $F$ at the one-point space.

Likewise the homogeneous quadratic functor $(D_2F)(X)$ is necessarily given by
$$
(C_2\wedge (X\wedge X))_{h\Sigma _2}
$$
where $\Sigma _2$ is the symmetric group, $C_2$ is some spectrum
with $\Sigma _2$-action, $X\wedge X$ has the obvious action, and the
subscript $h\Sigma _2$ denotes homotopy orbit spectrum. The coefficient
spectrum $C_2$, with its  $\Sigma _2$-action, is
called the second derivative of $F$ at the one-point space.

The explanation of why homogeneous quadratic functors always have this
form involves \it symmetric bilinear\rm\ functors. If $H$ is any $2$-excisive
functor (still from based
spaces to spectra), then its second order cross-effect, defined as the total
fiber (= homotopy fiber of homotopy fibers) of the square diagram
$$\CD
                           H(X_1\vee X_2) @>>> H(X_1)\\
                            @VVV @VVV\\
                               H(X_2) @>>>  H(*),
\endCD$$
will be a functor $L(X_1,X_2)$ of two variables, linear in each variable, and
symmetric with respect to interchanging the variables. It turns out that the
homogeneous quadratic part of $H$ can be recovered as $L(X,X)_{h\Sigma_2}$. A
bilinear functor necessarily has the form
$$
L(X_1,X_2)=C\wedge (X_1\wedge X_2),
$$
and a symmetry on $L$ means a $\Sigma_2$-action on $C$.

This pattern persists. For any $n\geq 1$ an  $n^{th}$ degree homogeneous
functor $H$ must have the form
$$
H(X)=L(X,\dots,X)_{h\Sigma_n}
$$
where
$L$, the $n^{th}$ order cross-effect of $H$, is a symmetric multilinear functor
of
$n$ variables. Thus $(D_nF)(X)$ will have the
form
$$
  (C_n\wedge (X\wedge\dots\wedge X))_{h\Sigma _n}
$$
at least for finite $X$. The object $C_n$, a spectrum with an action of
$\Sigma_n$, will be called the
$n^{th}$ \it derivative\rm\  of $F$ at the one-point space. We also 
refer to it as
the \it coefficient\rm\  of the multilinear functor, or of the homogeneous
functor.

It is worth emphasizing what brand of stable equivariant homotopy
theory is appearing here. Let $G$ be a finite
group. To do serious homotopy theory in the category of
$G$-spaces one usually takes the weak equivalences to be those equivariant maps
which induce weak equivalences on spaces of $H$-fixed points for all subgroups
$H$. (For CW objects this is equivalent to saying that the map has an 
inverse up to
$G$-homotopy.) There is also a \it very\rm\ weak notion of equivalence:
equivariant maps which are nonequivariantly weak equivalences. Using these as
the weak equivalences leads to a theory that is easier and less 
interesting, but
which has its uses. For example, a very weak equivalence of $G$-spaces
always induces a weak equivalence of homotopy orbit spaces and also of homotopy
fixed point spaces. Every $G$-space is equivalent in the very weak 
sense to a free
$G$-space (product with
$EG$), and for free
$G$-spaces the two kinds of equivalence coincide.

In the course of
investigating homogeneous functors we encounter spectra with finite groups
acting on them, and this leads us to an obvious stable analogue of 
this easier brand
of equivariant homotopy theory: Make a category whose objects are spectra
equipped with a $G$-action and whose maps are maps of spectra respecting the
action, and call such an equivariant map a (very) weak equivalence if it is
nonequivariantly a weak equivalence of spectra. It is these
objects (with $G=\Sigma_n$) up to this kind of equivalence which correspond
to homogeneous functors of degree $n$ up to weak natural equivalence.

In order to extend the main ideas of \it serious\rm\
equivariant homotopy theory to the stable setting, May and his collaborators
created a beautiful and powerful theory of
$G$-spectra \cite{17}. We mention this only to say that we do not need it here.
 From that
sophisticated point of view our spectra with $G$-action are what are sometimes
called
\it naive\rm\
$G$-(pre)spectra, namely
$G$-spectra involving only trivial actions of $G$ on suspension 
coordinates. The
$G$-spectrum expert will know what else to say about the place of the naive
theory in the sophisticated one, but neither the expert nor the 
novice should have
to think about the sophisticated theory to read this paper (unless 
the expert just
cannot help thinking about it).

Returning to the towers, a simple and familiar example is the \lq\lq Snaith
splitting\rq\rq\ of the functor
$F(X)=\Sigma^{\infty}\Omega\Sigma X$. Its
$n^{th}$ homogeneous part is $(D_nF)(X)\sim \Sigma^{\infty}(X^{\wedge n})$, and
its tower splits: $P_nF\sim \prod_{1\leq k\leq n}D_kF$. The limit of the tower
is the product of all the layers. If
$X$ is connected then the tower converges to $F(X)$ in the sense that the
map from $F(X)$ to the homotopy limit of the tower is a weak equivalence,
and in fact in the stronger sense that the map from $F(X)$ to $(P_nF)(X)$ has a
connectivity tending to infinity with $n$.

A related example is $F(X)= \Sigma^{\infty}\Omega X$. Here we
have
$$
(D_nF)(X)\sim \Omega^n\Sigma^{\infty}(X^{\wedge n}).
$$
(We will give a quick proof of this using the previous example and the general
fact
$P_n(F\circ\Sigma)\sim (P_nF)\circ\Sigma$. See Example 1.20 below.) The
$n^{th}$ derivative is the wedge of $n$ factorial copies of the 
spectrum $S^{-n}$
permuted transitively by $\Sigma_n$. The tower does not split. It converges
(in the strong sense again) when the space $X$ is $1$-connected.

A more general example is $F(X)= \Sigma^{\infty}Map_*(K,X)$ where $K$ is a
based finite complex and $Map_*$ means the space of based maps. Here the
$n^{th}$ derivative is the S-dual of a based complex with
$\Sigma_n$-action, namely the quotient of the smash product $K^{\wedge n}$
by the fat diagonal. $F$ is $dim(K)$-analytic, and the tower 
converges to $F(X)$
if
$X$ is
$dim(K)$-connected. Arone \cite{4} thoroughly explored this class of
examples, giving a concrete description of the stages $P_nF$ and not
just the layers $D_nF$.

The $n^{th}$ coefficient of Waldhausen's $A(X)$ is again the S-dual of a based
finite complex with $\Sigma_n$-action, namely
$({\Sigma_n})_+\wedge_{C_n}S^{n-1}$, where $C_n\subset \Sigma_n$ is the
transitive cyclic subgroup of order $n$ and the sphere $S^{n-1}$ is the
one-point compactification of the reduced regular representation of $C_n$.

So far we have been concerned with functors from based spaces
to spectra. We now discuss three different variants of this setup: 
functors from unbased
spaces, functors from spaces over a fixed space, and functors to spaces.
In these new settings the Taylor tower construction goes through with 
no change,
but some additional work is needed to understand what a homogeneous functor
looks like.

The switch from functors $\Cal T@>>>\Cal Sp$ of based spaces to functors
$\Cal U@>>>\Cal Sp$ of unbased spaces is fairly innocuous. If $\Cal
T@>\phi >>\Cal U$ is the forgetful functor then
$(P_nF)\circ\phi=P_n(F\circ\phi)$ and 
$(D_nF)\circ\phi=D_n(F\circ\phi)$. 

A snag appears when one tries to 
relate homogeneous functors to
symmetric multilinear functors, since the definition of the
cross-effect requires basepoints.
The good news is that in the end this
does not matter: homogeneous functors $\Cal T@>>>\Cal Sp$ of any
degree extend
uniquely (in an appropriate up-to-natural-weak-equivalence sense) to
$\Cal U$, and the same is true for symmetric multilinear functors in any
number of variables. The proof involves the left adjoint $\psi$ of $\phi$, in
other words the functor that adds a disjoint basepoint to an unbased space. For
example, although the usual definition of the suspension spectrum of a space
$X$ requires $X$ to be based, there is a well-known extension to unbased
spaces. It associates to
$X$ the homotopy fiber of
$$
\Sigma^{\infty}\psi (X)@>>>\Sigma^{\infty}\psi (*).
$$
If $X$ is based then this is naturally (with respect to based maps) equivalent
to $\Sigma^{\infty}X$, and on the other hand if $\Cal U@>L>>\Cal Sp$ is \it
any\rm\ linear functor such that $L\circ\phi$ is (naturally
equivalent to) $\Sigma^{\infty}$ then
$L$ must be as defined above. This was explained in \cite{1} and is
generalized to higher degrees here.

Note that for inhomogeneous functors this goes very wrong. For example,
the
$1$-excisive functor
$$
\Sigma^{\infty} \psi (*)\vee\hofiber(\Sigma^{\infty}\psi (X)
@>>>\Sigma^{\infty}\psi
(*)).
$$
from $\Cal U$ to $\Cal Sp$ is genuinely different from $\Sigma^{\infty}\psi
(X)$, but the difference disappears after restriction to
$\Cal T$. (The empty set is the only space $X$ at which the functors
disagree, but even for nonempty spaces the equivalence cannot be
chosen to be natural with respect to unbased maps.) An example of a
$2$-excisive functor $\Cal T@>>>\Cal Sp$ which does not extend to $\Cal U$
at all is the suspension spectrum of
$$
F(X)=((X\wedge X)/\Delta_X)_{\Sigma_2},
$$
the orbit space for the
$\Sigma_2$-action on the quotient of $X\wedge X$ by the diagonal. If $X$ is the
disjoint union of a point and a circle then the rational homology of $F(X)$
depends on where the basepoint is placed in $X$.

The second switch, from functors of spaces to functors of spaces over a
fixed space, is something like the switch from MacLaurin series to general
Taylor series. Instead of building, for each space $X$, a tower that has
$F(*)$ at the bottom and attempts to converge to $F(X)$, one can build, for
each map of spaces $X@>>>Y$, a tower that has $F(Y)$ at the bottom
and attempts to converge to $F(X)$. If we fix $Y$ and think of everything in
sight as a functor of spaces over $Y$ then the $n^{th}$ stage of the tower is
$n$-excisive. As in the case $Y=*$, there are two options, each with its own
technical advantages: the category
$\Cal U_Y$ of plain spaces over $Y$ and the category $\Cal T_Y$ of spaces over
$Y$ equipped with a section. We  sometimes refer to the latter as
fiberwise based spaces over $Y$.

The correspondence between homogeneous and symmetric multilinear
functors works the same for spaces over $Y$ as it does for spaces, but the
business of describing a multilinear functor by coefficient spectra becomes
more complicated when $Y$ has more than one point. As a trivial example, to
describe a linear functor of spaces over the two-point space $\lbrace
y_1,y_2\rbrace$ one needs two spectra. In general a linear functor 
$L$ of spaces
over
$Y$ determines a spectrum for each point $y\in Y$, namely the
coefficient $L(Y\vee_y S^0)$ of the linear functor
$Z\mapsto L(Y\vee_y Z)$ from $\Cal T$ to $\Cal Sp$, where $Y\vee_y Z$
means the wedge sum of $(Y,y)$ with the based space $Z$, viewed as a space
over $Y$. This spectrum depends
\lq\lq continuously\rq\rq\ on the point $y$ in some sense. If $Y$
is path-connected then up to weak homotopy equivalence the spectrum is
independent of the point; but there is a twisting that must not be overlooked.

In \cite{1} we defined
the \it differential\rm\ $D_YF$ of a functor $\Cal U@>F>>\Cal Sp$ at a space
$Y$ to be a linear functor from spaces over $Y$ to spectra, the linear
approximation to the functor
$$
(X@>>>Y)\mapsto \hofiber(F(X)@>>>F(Y))
$$
We also defined the \it derivative\rm\ to be a spectrum $\partial_yF(Y)$
depending on a space
$Y$ and a point
$y\in Y$,
namely $(D_YF)(Y\vee_y S^0)$; it might be called the $y$ coefficient of
$D_YF$ or the partial derivative of $F$ at $Y$ in the $y$ direction.
The relationship between the differential $D_YF$ (a
linear functor of spaces over $Y$) and the derivative
$\partial_yF(Y)$ (a spectrum) is something like the relationship between the
differential of a function (a linear function on the tangent space) 
and a partial
or directional derivative (a number, which specifies the behavior of that
linear function on a certain one-dimensional tangent subspace). The spectrum
$\partial_yF(Y)$ records the derivative of $F$ at $Y$ in the \lq\lq
direction\rq\rq\ given by $y$. Here we
make the multilinear generalization, defining the
$n^{th}$ differential to be the symmetric multilinear functor corresponding to
the
$n^{th}$ layer of the Taylor tower and defining the $n^{th}$ derivative
$$
\partial^{(n)}_{y_1,\dots,y_n}F(Y)
$$
to be its value at $(Y\vee_{y_1}S^0,\dots,Y\vee_{y_n}S^0)$.

The third
switch, from spectrum-valued to space-valued functors, brings a real surprise.
The construction of the Taylor series goes through with no change. If the
functor
$F$ happens to be $\Omega^{\infty}G$ for some spectrum-valued functor $G$,
then we have $P_nF=\Omega^{\infty}P_nG$ and
$D_nF=\Omega^{\infty}D_nG$. The surprise is that although most functors are
not of the form $\Omega^{\infty}G$ this is not reflected at the homogeneous
level: every homogeneous functor $\Cal T@>>>\Cal T$ of degree $n\geq 1$ is
infinitely deloopable, in fact in a unique and functorial way. Thus $D_nF$
always has the form $\Omega^{\infty}((C_n\wedge X^{\wedge
n})_{h\Sigma_n})$ and even in the space-valued case we can introduce a
spectrum called the
$n^{th}$ derivative of the functor, or the coefficient spectrum of the
$n^{th}$ homogeneous layer.

A central example of a functor $\Cal T@>>>\Cal T$ that is interesting from
this point of view is the identity. Let us call it $I$. Its first 
derivative is the
sphere spectrum. It is easy to see, in any of a number of ways, that 
the $n^{th}$
derivative is equivalent to the wedge sum of $(n-1)!$ copies of the
$(1-n)$-sphere spectrum, with
$\Sigma_n$ acting in such a way that the subgroup $\Sigma_{n-1}$ freely
permutes the summands. To identify $D_nI$ one has to know the action of the
full group
$\Sigma_n$.
Johnson \cite{5} gave an explicit finite complex with
$\Sigma_n$-action whose S-dual is the answer.  Arone
and Mahowald \cite{6} gave a different answer of that kind, showed that it
was equivalent to Johnson's, and used it to make some very interesting
calculations.

A basic example of an inhomogeneous functor that does not deloop is
$P_2I$, the
$2$-excisive approximation of the identity functor.
The functor $P_1I$ is of course $Q=\Omega^{\infty}\Sigma^{\infty}$.
$(P_2I)(X)$ is another standard object, the homotopy fiber of the James-Hopf
map
$$
QX@>>>Q((X\wedge X)_{h\Sigma_2}).
$$
(which is not induced by a spectrum map from $\Sigma^{\infty}X$ to
$\Sigma^{\infty}(X\wedge X)_{h\Sigma_2})$.

One notable feature of this last example turns out to be quite general: the
fibration sequence
$$
D_nF@>>>P_nF@>>>P_{n-1}F
$$
can always be extended one step to the right, as long as the constant
functor $P_0F$ is contractible. This falls out of the proof of (and is
essentially equivalent to) the deloopability of homogeneous functors.

The construction of the Taylor tower also goes through for functors taking
values in unbased instead of based spaces, but in that case it is nonsense to
speak of the layers $D_nF$ as functors, since that would involve using \lq\lq
the\rq\rq\ fiber of a map of unbased spaces. This small fussy point is really
not so small. The identity functor $I$ above was from
based spaces to based spaces. Let $J$ be the identity functor from unbased
spaces to unbased spaces. Then $(P_1J)(X)$ is the homotopy fiber of
$$
Q(X_+)@>>>QS^0,
$$
but not with respect to the usual basepoint in $QS^0$. In particular
$(P_1J)(\emptyset)$ is empty, so that there is no natural basepoint in
$(P_1J)(X)$. This seriously interferes with defining $D_2J$.
$P_1J$ is excisive and reduced, but perhaps too badly twisted to be rightly
called linear. For such reasons we hesitate to even speak of homogeneous
functors to unbased spaces.

The paper is organized as follows:

\S $1$ defines the tower $\lbrace P_nF\rbrace$ in the general case,
proves that $P_nF$ is $n$-excisive, establishes the universal
mapping property of $P_nF$, and notes the convergence of the tower in the
case of an analytic functor $F$.

\S 2 shows that homogeneous space-valued functors
can be functorially delooped and concludes that they correspond precisely to
homogeneous spectrum-valued functors.

\S 3 establishes the correspondence between homogeneous functors and
symmetric multilinear functors in the case of functors from
fiberwise based spaces over $Y$ to spectra. By \S 2
this result extends to space-valued functors.

\S 4 shows that for homogeneous functors, and also for symmetric
multilinear functors, it does not matter whether the domain category is
(fiberwise) based or unbased spaces, so that some of the results of
\S 3 extend to functors of plain spaces over $Y$.

\S 5 establishes notation for the symmetric multilinear functors
that (according to sections 2--4) encode the homogeneous layers of a
Taylor tower, and develops the idea of
coefficient spectra for multilinear functors.

\S 6 establishes a useful tool for working out examples. The key point is
that the $n^{th}$ cross effect of the homogeneous functor $D_nF$ can be made
by \lq\lq multilinearizing\rq\rq\ the $n^{th}$ cross effect of $F$ itself.

\S 7 works out the
$n^{th}$ derivatives of functors like $\Sigma^{\infty}Map_*(K,-)$.

\S 8 recalls and
discusses known results on the Taylor tower of the identity.

\S 9 indicates how to get the
$n^{th}$ derivative of Waldhausen's $A$, taking the first derivative 
as starting
point.

A few words about set theory are in order. We all know that we must 
not speak of
the set of all sets or the set of all topological spaces; set theory, 
as formulated
to avoid Russell-type paradoxes, does not allow it. And since the 
category $\Cal
C$ of spaces is not small, there may be
objections to speaking of the category of all functors, or homotopy 
functors, from
$\Cal C$ to $\Cal D$ in that case. Even if $\Cal C$ and
$\Cal D$ are both small, the functor category will not be small in 
general, so that
it is illegal to make a new category by using generators and 
relations to invert
some morphisms.

We choose to dodge this as follows. As usual
when basing mathematics on set theory, we begin by fixing a universe
$\Cal U$ of sets. Now by topological spaces we mean those whose
point sets are sets in this strict sense. The category of such spaces is not a
small category, any more than the set of all ($\Cal U$-)sets is a 
($\Cal U$-)set.
Nevertheless, if we are willing to work in two universes [\Cite{7}; I.6],
we are not in such bad shape. Introduce a larger universe $\Cal U'$ 
in which the set
of all $\Cal U$-sets is a set. The category of all ($\Cal U$-)spaces is then
$\Cal U'$-small. In
$\Cal U'$ there can be no objection to speaking of the category of all functors
from spaces to spaces, or of the full subcategory of homotopy functors. To
invert the weak equivalences by generations and relations, one can 
always pass to
a third universe in which the category of homotopy functors is small.

This solution will not suit all tastes. It may be a bit wasteful and it
may be a bit crude. For some more refined purposes it will be inadequate. For
example, if one wants to introduce some sort of space of maps between two
homotopy functors such that the set of components will be the set of
morphisms in the homotopy category, then it will certainly be a 
drawback to find
that this \lq\lq function space\rq\rq\ is not a space in the original sense. In
general, if one wants to work very seriously with the homotopy category of
functors then one will probably want to introduce a closed model structure on
the functor category, with all the benefits that that brings. (In 
fact,  although
we have not attempted to do so, there are certainly many
reasons for reworking this whole theory in the context of closed 
model categories.
The objects of study should be functors from one (small?) model category to
another, subject to some mild axioms, and the category of homotopy functors
should turn out to be a model category, too.)

On the other hand, in this paper
the role of the homotopy categories is a modest one: they are used here
mainly as shorthand, to make some sentences briefer and more readily
comprehensible than they otherwise would be. The \lq\lq multiple
universes\rq\rq\ device allows us to use this shorthand without telling any
lies, but without making the catogory of spaces artificially small.

\eject

\head \S1. The Taylor tower\endhead

Let $\Cal C@>F>> \Cal D$ be a homotopy functor, where $\Cal C=\Cal C_Y$ is
either $\Cal U_Y$ or $\Cal T_Y$ and
$\Cal D$ is either $\Cal U$, $\Cal T$, or $\Cal Sp$. Let $n\geq 0$.

The $n$-excisive approximation $P_nF$ will be constructed by the infinite
iteration of another construction $T_n$ which is designed to
bring the functor $F$ a little closer to being $n$-excisive.  The
special case when $n=1$ was treated in [\Cite{1}; Def.\ 1.10], where
$T_1$ and $P_1$ were called $T$ and $P$. Thus if $Y=*$ and
$F(*)=*$ then $T_1F\sim\Omega F\Sigma$ and
$P_nF\sim\hocolim_{i\geq 0}\Omega^i F\Sigma^i$ .

The definition of $T_n$ uses the \it fiberwise join \rm over $Y$,
as introduced in the proof of [\Cite{2}; 5.1]. Let $X$ be a space over $Y$
and let $U$ be a space. (In most applications $U$ will be a discrete
finite set.)  The fiberwise join is the space
$$
X*_YU=\hocolim (X\leftarrow X\times U \a Y\times U)
$$
considered as a space over $Y$. The name
signifies that the functor
\lq\lq join with
$U$\rq\rq\ is being applied to all the fibers of $X@>>>Y$.

If $U$ has one element, then $X*_YU@>>>Y$ is the \it fiberwise cone\rm\ of
$X$ over $Y$ (the mapping cylinder of $X@>>>Y$ considered as a space over $Y$).
If $U$ has two elements then
$X*_YU@>>>Y$ is the \it fiberwise suspension\rm\ and will sometimes be
denoted by
$\Sigma_YX$.

It
should be noted that this construction has the best of both worlds, in
that on the one hand each fiber of
$X*_YU@>>>Y$ is homeomorphic to the join with $U$ of the corresponding fiber
of
$X@>>>Y$, and on the other hand each homotopy fiber is weakly homotopy
equivalent to the join with $U$ of the corresponding homotopy fiber of
$X@>>>Y$.

Recall that even if $\Cal C$
is
$\Cal T_Y$ rather than
$\Cal U_Y$ then $X*_YU$ is still in $\Cal C$;
it inherits a section from $X$. In other words, if $X$
is a fiberwise based space over $Y$ then
$X*_YU$ is also canonically fiberwise
based, without for example choosing a point in $U$.

The object $X*_YU$ depends
(bi-)functorially on
$X$ and $U$. Note also that there is a natural isomorphism of spaces over $Y$:
$$
                     (X*_YU)*_YV\cong X*_Y(U*V)
$$
where $U*V$ is the ordinary join of two spaces.

Let $\Cal P(\underline{n+1})$ be the poset of subsets of
$\underline{n+1}=\lbrace 1,\hdots ,n+1\rbrace$. Any object
$X\in \Cal C$ yields an
$(n+1)$-dimensional cubical
diagram
  in $\Cal C$.
$$
\CD
\Cal P(\underline{n+1}) \a \Cal C\\
U\mapsto X*_YU
\endCD
$$
Consider the composed
functor
$U\mapsto
F(X*_YU)$, a cubical diagram in $\Cal D$.  The homotopy limit of
its restriction to $\Cal P_0(\underline{n+1})$,
the poset of nonempty subsets of $\underline{n+1}$, will be called
$(T_nF)(X)$. Clearly this yields a homotopy functor $\Cal C @>{T_nF}>> \Cal
D$. There is a natural
map $F @>{t_nF}>> T_nF$, since any cubical diagram determines a map from the
\lq\lq initial\rq\rq\  object to the homotopy  limit of the others (see
[\Cite{2}; Def.\ 1.2]):
$$
F(X) = F(X*_Y\emptyset )  @>{t_nF}>>  \holim_{U\in\Cal
P_0(\underline{n+1})}F (X*_YU) = (T_nF)(X)
$$
Let $(P_nF) (X)$ be the sequential homotopy colimit of the diagram
$$
F (X) @>{(t_nF)(X)}>> (T_nF)(X) @> {(t_n T_nF)(X)}>> (T_n^2F )(X) @>{(t_n
T_n^2F)(X)}>>
\hdots.
$$
Clearly $P_nF$ is
a homotopy functor and we have a natural map $F@>{p_nF}>> P_nF$.

The cubical diagram $U\mapsto
X*_YU$ is strongly cocartesian for all
$X$. Therefore if the functor $F$ is $n$-excisive then the maps 
$t_nF$ and $p_nF$
will be weak equivalences for all $X$ . In this sense $n$-excisive functors are
unchanged by
$T_n$ and $P_n$.

We will see below that $P_nF$ is always $n$-excisive, and that (in a 
homotopy category) it is the best
$n$-excisive approximation to $F$ in a categorical sense.

\remark{1.1 Remark} If the map $X@>>>Y$ is $m$-connected, then for
nonempty $U$ the map $X*_YU@>>>Y$ is $(m+1)$-connected, and therefore
$(T_n^iF)(X)$ depends only on the behavior of $F$ on objects whose
structural maps to $Y$ are $(m+i)$-connected. In this sense the Taylor
approximations of a functor
$\Cal C_Y@>F>>\Cal D$ depend only on the restriction of the functor to
objects \lq\lq arbitrarily close to $Y$\rq\rq\ (just as the Taylor
expansion of $f(x)$ in powers of $x-y$ depends only on the behavior 
of $f(x)$ in
an arbitrarily small neighborhood of $y$). A related comment is that $P_n$
commutes (up to natural isomorphism) with fiberwise suspension:
$$
P_n(F\circ\Sigma_Y)\cong(P_nF)\circ\Sigma_Y.
$$
This follows from the natural isomorphism
$$
\Sigma_Y(X*_YU)\cong(\Sigma_YX)*_YU.
$$
\endremark

Before giving the general proof that $P_nF$ is $n$-excisive, we give a
different proof under the added assumption that $F$ is \it stably
$n$-excisive\rm\ [\Cite{2}; 4.1], and in this case we show that
$P_nF$ approximates
$F$ not only in a categorical sense but also from the point of view of
connectivity.

The first task is to show that if
$F$ is stably
$n$-excisive then
$T_nF$ is more nearly $n$-excisive than $F$ is, and that $T_nF$ agrees
with $F$ to $n^{th}$ order in the following sense:

\definition{1.2 Definition}A map $u:F\rightarrow G$ between
two functors from
$\Cal C_Y$ to $\Cal D$ is said to satisfy
$O_n(c,\kappa)$ if, for every $k \geq\kappa$,
for every object $X$ of $\Cal C_Y$ such that the map
$X\a Y$ to the final object is
$k$-connected, the map
$u_X:F(X)\a G(X)$ is $(-c+(n+1)k)$-connected. We say that
$F$ and $G$ \it agree to order \rm $n$ (via the map $u$) if this
holds for some constants $c$ and $\kappa$.
\enddefinition

\remark{1.3 Remark}The letter $O$
stands for \lq\lq osculation\rq\rq\. This condition on (a map 
between) two functors is analogous to a condition on two
\it functions\rm\
(say real functions of one or several variables). The functions $f$ 
and $g$ can be said to agree to
$n^{th}$ order at $y$ if there are constants $C$ and $K$ such that 
for every $x$ such that $|x-y|\leq K$ we have $|f(x)-g(x)|\leq 
C|x-y|^{n+1}$ .  \endremark

The condition $E_n(c,\kappa)$, stable $n^{th}$ order excision, was defined in
[\Cite{2}: Def.\ 4.1].

\proclaim{1.4 Proposition}If $F$ satisfies $E_n(c,\kappa)$, then
\roster
\item $T_nF$ satisfies
$E_n({c-1},\kappa -1)$ and
\item $t_nF:F\a T_nF$ satisfies
$O_n(c,\kappa)$.
\endroster
\endproclaim

\demo{Proof}
The second conclusion is immediate from the definitions. For
the first, note that the functor
$X\mapsto X*_YU$ from $\Cal C$ to $\Cal C$ preserves
strongly cocartesian cubes and (unless $U$ is empty) increases the
connectivity of maps. It follows immediately that for each nonempty $U$ the
functor
$X\mapsto F(X*_YU)$ from $\Cal C$ to $\Cal D$ satisfies
$E_n(c-(n+1),\kappa -1)$.  By [\Cite{2}; 1.20], $T_nF$
satisfies $E_n(c-1,\kappa -1)$.
\enddemo

\proclaim{1.5 Proposition}If $F$ is stably
$n$-excisive, then
\roster
\item $P_nF$ is $n$-excisive and
\item $F$ agrees with
  $P_nF$ to order $n$ (via $p_nF$).
\endroster
\endproclaim
\demo{Proof}Suppose that $F$ satisfies $E_n(c,\kappa)$. By
1.4(1) and induction, $T_n^iF$
satisfies
$E_n(c-i,\kappa -i)$. This easily implies that $P_nF$ is
$n$-excisive. We also find, by 1.4(2), that
$t_nT_n^iF$ satisfies
$O_n(c-i,\kappa -i)$, and in particular $O_n(c,\kappa)$,
for all $i$. It follows easily that
$O_n(c,\kappa)$ is also satisfied by the composed maps
$(t_nT_n^iF)\circ\hdots\circ
t_nF$ and in the limit by $p_nF$.
\qed\enddemo

\proclaim{1.6 Proposition}Let $F@>u>>G$ be a map between
homotopy functors. If $F$ and $G$ agree to $n^{th}$ order
via $u$, then the induced map
$P_nF@>{P_nu}>>P_nG$ is an equivalence. The converse holds if
$F$ and $G$ are stably
$n$-excisive.
\endproclaim
\demo{Proof}Suppose $u$ satisfies $O_n(c,\kappa)$. For each
nonempty finite set $U$ the resulting map of functors
$$
F(-*_YU)\a G(-*_YU)
$$
satisfies $O_n(c-(n+1),\kappa -1)$; again we have used the fact
that the functor
$X\mapsto X*_YU$ from $\Cal C$ to $\Cal C$ preserves
strongly cocartesian cubical
diagrams and increases the connectivity of maps.
Therefore, using [\Cite{2}; 1.20] as in the
proof of 1.4, $T_nF@>{T_nu}>> T_nG$ satisfies
$O_n(c-1,\kappa-1)$. By induction,
$T_n^iF@>{T_n^iu}>> T_n^iG$ satisfies $O_n(c-i,\kappa -i)$. It follows,
letting $i$ tend to infinity, that $P_nF@>>>P_nG$ is an equivalence.
For the converse, use 1.5(2) and the commutative diagram
$$
\CD
F @>u>> G\\
@V{p_nF}VV @VV{p_nG}V\\
P_nF @>{P_nu}>> P_nG
\endCD
$$
\vglue -12mm\hbox{}{\qed}
\medskip
Note that when $F$ is stably $n$-excisive $P_nF$ can be
characterized, up to natural equivalence, as the only $n$-excisive
functor that agrees to
$n^{th}$ order with $F$.

The remaining results of \S 1 (except for the last sentence of 1.13)
have nothing to do with connectivity. In particular, functors are not assumed
to satisfy any kind of stable excision hypothesis.

\proclaim{1.7 Proposition}Up to equivalence,
\roster
\item $T_n$ commutes with $\holim$
\item $P_n$ commutes with finite $\holim$
\item $T_n$ and $P_n$ commute with $\hofiber$.
\item $T_n$ and $P_n$ commute with filtered (in particuloar sequential)
$\hocolim$.
\item for spectrum-valued functors, both $T_n$
and $P_n$ commute with any $\hocolim$.
\endroster
\endproclaim

\demo{Proof and explanation}The main point is that $\holim$
commutes with $\holim$, $\hocolim$ commutes with $\hocolim$,
and, up to equivalence, finite $\holim$ commutes with filtered
$\hocolim$.

(1)\qua This one is true up to isomorphism: If
$\lbrace F_\alpha
\rbrace$ is any sort of diagram of homotopy functors
$\Cal C@>{F_\alpha}>> \Cal D$ and $F$ is given by $F(X) =
\holim_\alpha F_\alpha(X)$ then we have
$$\align
                    (T_nF)(X)
  &\cong \holim_U \holim_\alpha F_\alpha(X*_YU)\\
&\cong \holim_\alpha \holim_U F_\alpha(X*_YU)\\
&\cong \holim_\alpha  (T_nF_\alpha)(X)
\endalign
$$

(2)\qua In the same situation we have
a natural map
$$\align
                   (P_nF)(X) &\cong \hocolim_i (T_n^iF)(X) \\
                             &\cong \hocolim_i \holim_\alpha (T_n^iF_\alpha)(X)
\\
                             &\a \holim_\alpha \hocolim_i (T_n^iF_\alpha)(X)
\\
                               &\cong \holim_\alpha (P_nF_\alpha(X)
\endalign
$$
It is an equivalence in the case when $\lbrace F_\alpha
\rbrace$ is a finite diagram.

(3)\qua This follows from (1) and (2). We have
$$\align
P_n\hofiber (F@>>>G)&\cong P_n\holim (F@>>>G@<<<*)\\
&\sim\holim (P_nF@>>>P_nG@<<<P_n*)\\
&\sim\holim (P_nF@>>>P_nG@<<<*)\\
&\cong\hofiber (P_nF@>>>P_nG)
\endalign
$$
and likewise for $T_n$.

(4)\qua $T_n$ is a finite $\holim$ and $P_n$ is a $\hocolim$ of finite
$\holim$'s.

(5)\qua Since hocolims commute, $\hocolim$ preserves
cocartesian cubes. Thus
$\hocolim$ of spectra preserves cartesian cubes.
This implies (using [\Cite{2}; 1.19])
that $T_n$ commutes with $\hocolim$ up to equivalence.
Therefore the
same holds for $P_n=\hocolim_iT^i$.
\qed\enddemo

Let $\Cal F(\Cal C,\Cal D)$ be the category whose objects are the
homotopy functors from $\Cal C$ to $\Cal D$, and whose morphisms are the
natural maps. Let $h\Cal F(\Cal C,\Cal D)$ be its
\it homotopy category\rm ; it has the same objects and is obtained by
formally inverting the equivalences. Because the functors
$T_n$ and
$P_n$ from
$\Cal F(\Cal C,\Cal D)$ to itself take equivalences to equivalences,
they give rise to
functors from the homotopy category to itself. Morphisms in
$h\Cal F(\Cal C,\Cal D)$ are called \it weak maps\rm\. Any functor
weakly isomorphic to an
$n$-excisive functor is $n$-excisive.

\proclaim{1.8 Theorem}For any homotopy functor $\Cal C@>F>>\Cal
D$, the functor $P_n F$ is $n$-excisive. In the homotopy category
$h\Cal F(\Cal C,\Cal D)$,
$p_nF$ is the universal map from $F$ to an
$n$-excisive functor.
\endproclaim

The key to 1.8 is the
following:

\proclaim{1.9 Lemma}Let $\Cal X$ be any strongly cocartesian
$(\underline{n+1})$-cube in
$\Cal C$, and let $F$ be any homotopy functor. Then the map of cubes
$$
                  F(\Cal X) @>{(t_nF)(\Cal X)}>> (T_nF)(\Cal X)
$$
  factors through some cartesian cube.\endproclaim

\demo{Proof of 1.8, assuming 1.9}Let $\Cal X$ be any strongly
cocartesian
$(\underline{n+1})$-cube in
$\Cal C$, and consider the diagram
of cubes
$$
                  F(\Cal X) \a (T_nF)(\Cal X) \a (T_n^2F)(\Cal X) \a \dots
$$
which leads by $\hocolim$ to the cube $(P_nF)(\Cal X)$.
By 1.9 each of the maps of cubes
displayed above factors through some cartesian cube. Therefore the cube
$(P_nF)(\Cal X)$ is equivalent to a sequential $\hocolim$ of
cartesian cubes, and so it is itself cartesian. This shows that $P_nF$ is
$n$-excisive.

For the existence half of the universal mapping property, let
$F@>u>>P$ be any weak map to an $n$-excisive functor. The diagram
$$\CD
F@>u>>P\\
@Vp_nFVV @Vp_nPVV\\
P_nF@>P_nu>> P_nP
\endCD
$$
shows that $u$ factors through $p_nF$, since $p_nP$ as a weak map is
invertible.

For the uniqueness we must show that if $P$ is $n$-excisive then a weak
map
$P_nF@>v>>P$ is determined by the composition $v\circ p_nF$. It
suffices if in the diagram of weak maps
$$\CD
F@>p_nF>>P_nF@>v>>P\\
@Vp_nFVV @Vp_nP_nFV\sim V@Vp_nPV\sim V\\
P_nF@>P_n(p_nF)>\sim>P_nP_nF@>P_nv>> P_nP
\endCD
$$
the marked ($\sim$) arrows are invertible,
for then $v$ is determined by $P_nv$, which is determined by
$P_nv\circ P_n(p_nF)=P_n(v\circ p_nF)$, which is of course determined by
$v\circ p_nF$.

The marked vertical arrows are equivalences because
$P_nF$ and
$P$ are $n$-excisive. In order for  $P_n(p_nF)$
to be an equivalence, it will be enough (by 1.7(4)) if
$P_n(t_nF)$ is an equivalence. Define the functor $J_UF$ by
$(J_UF)(X)=F(X*_YU)$. In the diagram
$$
\align
P_nF@>P_n(t_nF)>> &P_nT_nF
\cong P_n(\holim_{U\in\Cal
P_0(\underline{n+1})}J_UF)\\ @>\sim >>&\holim_{U\in\Cal
P_0(\underline{n+1})}P_nJ_UF
\cong\holim_{U\in\Cal
P_0(\underline{n+1})}J_UP_nF
\endalign
$$
the second arrow is an equivalence by 1.7(2). The
composition is also an equivalence; this simply means that for
every $X$ the functor
$P_nF$ takes the cube $\lbrace X*_YU\rbrace$ to a
cartesian cube, and it is true because $P_nF$ is $n$-excisive.
It follows that $P_n(t_nF)$ is an equivalence.
\qed\enddemo

The fact that $P_n(p_nF)=p_nP_nF$ in the
homotopy category follows from 1.8; if we had had it before 1.8,
then we could have skipped the last part of the proof of 1.8.

\remark{1.10 Remark}Theorem 1.8 is one of several statements in this paper
for which we have two different proofs: an older one which requires some
kind of stable excision hypothesis, and a newer one which has nothing to do
with connectivity. Another such statement is Theorem 2.1, which also depends
(near the end of the proof of 2.2) on 1.9. Another is 3.1, and
another is 6.1. The older proofs have the advantage of a certain common-sense
quality, and if what we really care about here is convergent Taylor towers then
the older proofs are good enough; but the fact that the theorems are still true
without connectivity hypotheses is striking and the newer proofs are perhaps
worth looking at, too.
\endremark

\demo{Proof of 1.9}(Any reader who, like the author, finds this proof a little
opaque, may wish to look at it again after reading the proof of 3.2, 
which is related
but simpler.) Let us write $n$ instead of $n+1$. For subsets
$T,U_1,\dots,U_n$ of
$\underline n$, define $\Cal{\hat X}(T,U_1,\dots,U_n)$ to be
$$
\hocolim{(\Cal X(T)@<<<\coprod_{1\leq s\leq n}(\Cal X(T)\times
U_s)@>>>\coprod_{1\leq s\leq n}(\Cal X(T\cup{\lbrace
s\rbrace})\times U_s))}.
$$
This can also be described as the union, along
$\Cal X(T)$, of the spaces\nl
$\Cal X(T)*_{\Cal X(T\cup \lbrace s\rbrace)}U_s$.

Clearly $\Cal {\hat X}$ is a functor from
$\Cal P(\underline n)\times P(\underline n)^n$ to $\Cal C$. We have
$\Cal X(T)\cong\Cal {\hat X}(T,\emptyset,\dots,\emptyset)$, and there
is a natural map $\Cal {\hat X}(T,U,\dots,U)@>>>\Cal X(T)*_YU$ that
corresponds to the identity map when $U=\emptyset$.

Let $\Cal E$ be the set of all
$(U_1,\dots,U_n)\in {\Cal P_0(\underline{n})}^n$ such that, for at
least one $s\in
\underline{n}$,
$s\in U_s$. Since $\Cal E$ contains the image of the diagonal map
$\Cal P_0(\underline{n})@>>>{\Cal P_0(\underline{n})}^n$,
we can
factor the map
$t_{n-1}F(\Cal X (T))$ as follows:
$$
\align
F(\Cal X(T))  &\a \holim_{(U_1,\dots,U_n)\in \Cal E}  F(\Cal
{\hat X}(T,U_1,\dots,U_n))\\
                    &@>>> \holim_{U\in \Cal P_0(\underline n)}  F(\Cal
{\hat X}(T,U,\dots,U))\\
                    &@>>> \holim_{U\in \Cal P_0(\underline n)}  F((\Cal
X(T)*_YU)\\
&=(T_{n-1}F)(\Cal X(T)).
\endalign
$$
We claim that the cube
$$
T\mapsto \holim_{(U_1,\dots,U_n)\in \Cal E}  F(\Cal
{\hat X}(T,U_1,\dots,U_n))
$$
is cartesian (for any homotopy functor $F$).

Let $\Cal E^*$ be the set of all
$(U_1,\dots,U_n)\in {\Cal P_0(\underline{n})}^n$ such that, for
some
$s$, $U_s=\lbrace s\rbrace$. This poset is
left cofinal in $\Cal E$ (exercise; see [\Cite{2}; page 298] for the 
definition), and
it follows that the restriction map from the $\holim$ over $\Cal E$ to the
$\holim$ over $\Cal E^*$ is an equivalence. Therefore it is enough if,
whenever $(U_1,\dots,U_n)\in\Cal E^*$, the cube
$$
T\mapsto   F(\Cal
{\hat X}(T,U_1,\dots,U_n))
$$
is cartesian. In fact, it is cartesian for a very basic reason: if 
$U_{s_0}=\lbrace
s_0\rbrace$ then the map
$$
\Cal{\hat X}(T,U_1,\dots,U_n)@>>>
\Cal{\hat X}(T\cup\lbrace s_0\rbrace,U_1,\dots,U_n)
$$
is an equivalence. To see this, examine the diagram
$$\CD
\Cal X(T)@<<<\coprod_{s}(\Cal X(T)\times
U_s)@>>>\coprod_{s}(\Cal X(T\cup{\lbrace
s\rbrace})\times U_s))\\ @VVV @VVV @VVV\\
\Cal X(T\cup\lbrace s_0\rbrace)@<<<\coprod_{s}(\Cal
X(T\cup\lbrace s_0\rbrace)\times U_s)@>>>\coprod_{s}(\Cal
X(T\cup{\lbrace s, s_0\rbrace})\times U_s))
\endCD
$$
The fact that the induced map from $\hocolim$ of top row to
$\hocolim$ of bottom row is an equivalence follows from the trivial fact
that for $s\neq s_0$ the square
$$\CD
\Cal X(T)\times
U_s@>>>\Cal X(T\cup{\lbrace s\rbrace})\times U_s\\
@VVV @VVV\\
\Cal
X(T\cup\lbrace s_0\rbrace)\times U_s@>>>\Cal
X(T\cup{\lbrace s, s_0\rbrace})\times U_s
\endCD
$$
is cocartesian, plus the even more trivial fact that the square
$$\CD
\Cal X(T)@<<<\Cal X(T)\times
U_{s_0}\\
@VVV @VVV\\
\Cal X(T\cup\lbrace s_0\rbrace)@<<<\Cal
X(T\cup\lbrace s_0\rbrace)\times U_{s_0}
\endCD
$$
is cocartesian.
\qed\enddemo

\proclaim{1.11 Corollary}If $0\leq m\leq n$ then the map
$$
P_mF @>{P_m(p_nF)}>> P_mP_nF
$$
is an equivalence.
\endproclaim

\demo{Proof}This is formal; using the universal mapping properties of
$P_m$ and $P_n$ and the fact that $m$-excisive functors are $n$-excisive, one
sees that $P_mP_nF$ has the universal mapping property that characterizes
$P_mF$.
\qed\enddemo

We now collect the
\lq\lq Taylor polynomials\rq\rq\
$P_nF$ into a \lq\lq Taylor series\rq\rq\, by showing that there is
a natural map $P_nF@>{q_nF}>>P_{n-1}F$ satisfying
$q_nF\circ p_nF
= p_{n-1}F$.

An effortless way to produce $q_nF$ as a weak map would be to reason as in
the proof of 1.11:
$p_nF$ is the universal map from $F$ to an $n$-excisive functor
and
$(n-1)$-excisive functors are $n$-excisive, so there is a unique map $q_nF$ in
$h\Cal F(\Cal C,\Cal D)$ such that $q_nF\circ p_nF = p_{n-1}F$. On the other
hand, it is desirable to define
$q_nF$ in such a way that
$q_nF\circ p_nF = p_{n-1}F$ on the nose and not just weakly. Moreover, the
explicit construction for $q_nF$ will be useful in its own right in
proving Lemma 2.2.

We will make a (commutative) diagram
$$
\CD
F @>{t_nF}>> T_nF @>{t_nT_nF}>> T_n^2F
@>{t_nT_n^2F}>> \hdots\\ @VVV @V{q_{n,1}}VV @V{q_{n,2}}VV \\
F @>{t_{n-1}F}>> T_{n-1}F @>{t_{n-1}T_{n-1}F}>>
T_{n-1}^2F @>{t_{n-1}T_{n-1}^2}>> \hdots\\
\endCD\tag1.12
$$
and then define $q_nF$ as the induced map of horizontal homotopy colimits. We
must define $q_{n,i}$ and then verify that the squares commute.

Notice that $T_n^iF$ is naturally isomorphic to a homotopy limit indexed by a
product of $i$ copies of the partially ordered set $\Cal
P_0(\underline{n+1})$:
$$
  (T_n^iF)(X) \cong \holim_{(U_1,\hdots ,U_i)\in {\Cal
P_0(\underline{n+1})}^i}
F(X*_Y(U_1*...*U_i))
$$
From this point of view there is an obvious map
$T_n^iF@>{q_{n,i}}>> T_{n-1}^iF$, induced by
the inclusion of
${\Cal P_0(\underline {n}})^i$ in ${\Cal P_0(\underline
{n+1}})^i$.  The first square in (1.12) now obviously commutes.  The
${(i+1)}\st$ square will commute if both squares commute in
$$
\CD
T_n^iF @>{t_nT_n^iF}>> T_n^{i+1}F @= T_n^{i+1}F\\
@V{q_{n,i}F}VV @V{T_nq_{n,i}F}VV @V{q_{n,i+1}F}VV\\
T_{n-1}^iF @>>{t_nT_{n-1}^iF}> T_nT_{n-1}^{i}F @>>{q_{n,1}T_{n-1}^iF}>
T_{n-1}^{i+1}F
\endCD
$$
The left square commutes because $t_nF$ is natural in $F$. The other commutes
because it is induced by a commutative diagram of posets
$$
\CD
{\Cal P_0(\underline {n+1})}^{i+1} @= {\Cal P_0(\underline {n+1})}^{i+1}\\
@AAA @AAA\\
\Cal P_0(\underline {n+1})\times{\Cal P_0(\underline {n})}^i @<<< 
{\Cal P_0(\underline
{n})}^{i+1}
\endCD
$$
Summing up, we have:

\proclaim{1.13 Theorem}A homotopy functor $F$ from
spaces over $Y$ (with or without section) to either spaces, based 
spaces, or spectra, determines a tower
of such functors
$$
\CD
\vdots\\
@VVV\\
P_nF\\
@VV{q_nF}V \\
P_{n-1}F\\
@VV{q_{n-1}F}V\\
\vdots\\
  @VV{q_2F}V\\
P_1F @. {} \\
@VV{q_1F}V \\
  P_0F\\
\endCD
$$
and a map $F@>{\lbrace p_nF\rbrace}>>lim_n P_nF$. Each $P_n$ is a
homotopy functor (from homotopy functors to homotopy functors), and
$p_n$ and $q_n$ are natural. The functor $P_nF$ is
always $n$-excisive, and (in the homotopy category of homotopy functors)
$p_nF$ is universal among maps from $F$ to $n$-excisive functors.
If $F$ is
$\rho$-analytic and the structural map
$X\rightarrow Y$ is $(\rho +1)$-connected,
then the connectivity of the map
$F(X)@>{p_nF}>>(P_nF)(X)$ tends to $+\infty$ with $n$, so
that $F(X)$ is equivalent to the homotopy limit 
$(P_{\infty}F)(X\rightarrow Y)$ of the
tower.\endproclaim

\demo{Proof}
The last statement is the only new one. Recall from [\Cite{2}; 4.2] that
$\rho$-analyticity means that there is a number
$q$ such that, for all $n\geq 0$, $F$ satisfies $E_n(n\rho -q,\rho
+1)$. If $X\in \Cal C$ is such that the map $X\a Y$ is $k$-connected
with $k>\rho$ then by the proof of 1.5 the connectivity of the map
$F(X)@>{p_nF}>>(P_nF)(X)$ is at least $(q+k+n(k-\rho))$, which
tends to
$+\infty$ with
$n$.
\qed\enddemo

If $P_nF$ is analogous to an $n^{th}$ Taylor polynomial, then the homotopy
fiber
$$
D_nF = \hofiber(P_nF@>{q_nF}>>P_{n-1}F)
$$
is analogous to the $n^{th}$ term in a
Taylor series.  Notice that the definition of $D_nF$ is meaningful if the
category $\Cal D$ is either based spaces or spectra, but not if it is unbased
spaces. From $
P_n(F\circ\Sigma_Y)\cong(P_nF)\circ\Sigma_Y
$ (Remark 1.1) we have
$$
D_n(F\circ\Sigma_Y)\cong(D_nF)\circ\Sigma_Y.
\tag 1.14
$$
\definition{1.15 Definition}A homotopy
functor
$F:\Cal C\a \Cal D$ is called $n$-reduced if $P_{n-1}F\sim
*$. It is called $n$-homogeneous, or homogeneous of degree $n$, if it is both
$n$-excisive and $n$-reduced.
\enddefinition

Thus $1$-reduced means reduced, and $1$-homogeneous means
linear.

\remark{1.16 Remark}
If $n>1$ then it is not easy in general to test whether a functor $F$ is
$n$-reduced. Perhaps the main difficulty is that $P_{n-1}F\sim *$ does not
imply $T_{n-1}F\sim *$. On the other hand, by 1.6 a sufficient condition for
  $F$ to be
$n$-reduced is that
$F$ agrees to order
$n-1$ with the constant functor $*$, in other words that the connectivity of
$F(X)$ tends to infinity at least $n$ times faster than the connectivity of the
map
$X@>>>Y$. (If
$F$ is analytic, or even just stably
$n$-excisive, then this condition is also necessary for $F$ to be $n$-reduced.)
Thus, for example, for any spectrum $C$ that is bounded below ($k$-connected
for some
$k$) the
$n$-excisive functors $X\mapsto C\wedge X^{\wedge n}$ from
based spaces to spectra and $X\mapsto \Omega^{\infty}(C\wedge X^{\wedge n})$
from based spaces to based spaces are homogeneous. In fact this holds for all
spectra $C$, either by expressing $C$ as a homotopy colimit of bounded-below
spectra or by 3.1 below.
\endremark
\proclaim{1.17 Proposition}$D_nF$ is always $n$-homogeneous.
\endproclaim

\demo{Proof}The functor $D_nF=\holim
(P_nF@>>>P_{n-1}F@<<<*)$ is $n$-excisive because it is a
homotopy limit of $n$-excisive functors.
To see that it is also $n$-reduced, use 1.7(4) to identify
$P_{n-1}D_nF$ with the homotopy fiber of the map
$$
P_{n-1}(P_nF)@>{P_{n-1}(q_nF)}>>P_{n-1}(P_{n-1}F);
$$
by (1.11) this map is an equivalence.
\qed\enddemo

\proclaim{1.18 Proposition}Up to equivalence,
\roster
\item $D_n$ commutes with finite $\holim$.
\item $D_n$ commutes with $\hofiber$.
\item $D_n$ commutes with filtered $\hocolim$.
\item for spectrum-valued functors, $D_n$ commutes with
arbitrary
$\hocolim$.
\endroster
\endproclaim

\demo{Proof}This follows easily from 1.7.
\qed\enddemo

\example{1.19 Example}Let the functors $F_a$ and $F_b$ be
$a$-homogeneous and
$b$-homog\-en\-eous respectively, with $a<b$, and let $F_a@>f>> F_b$
be any natural map.  A
simple example is the diagonal inclusion $QX\a Q(X\wedge X)$.
Then the functor
$$
                        F(X)= \hofiber (F_a@>f>> F_b)
$$
is $b$-excisive.  Its tower has only two nontrivial layers
$D_aF\sim F_a$ and
$D_bF\sim \Omega F_b$; we have
$$\align
P_nF &\sim *\quad n< a \\
P_nF &\sim F_a  \quad a\leq n< b \\
P_nF &\sim F \quad n\leq b
\endalign$$
All of this follows from $P_nF\sim\hofiber
(P_nF_a@>>>P_nF_b)$. Note that when
$a>b$ there can be no interesting natural map
$F_a\a F_b$, since by (1.8) any such map factors (in $h\Cal F(\Cal C,\Cal D)$)
through
$P_{a-1}F_a\sim *$.
\endexample

\example{1.20 Example}There is a weak
equivalence (\lq\lq Snaith splitting\rq\rq)
$$
\Sigma^{\infty}\Omega\Sigma X \sim \prod_{n\geq
1}\Sigma^{\infty}X^{\wedge n}
\tag 1.21
$$
for based connected spaces $X$. The functor $\Sigma^{\infty}X^{\wedge
n}$ of $X$ is $n$-homogeneous. The $m^{th}$
Taylor polynomial of the right-hand side of 1.21 is
$\prod_{1\leq n\leq m}\Sigma^{\infty}(X^{\wedge n})$, by 1.6, since this finite
product agrees with the infinite product to order $m$. The same therefore holds
for the left-hand side, by 1.1.  In particular the
$n^{th}$ homogeneous layer of $
\Sigma^{\infty}\Omega\Sigma X$ is $\Sigma^{\infty}X^{\wedge
n}$. We will find later that this is enough to determine the $n^{th}$
homogeneous layer of $F(X)=
\Sigma^{\infty}\Omega X$. In fact, using 1.14 we have
$$
(D_nF)(\Sigma X)\sim (D_n(F\circ \Sigma))(X)
\sim \Sigma^{\infty}X^{\wedge n}\sim\Omega^n\Sigma^{\infty}(\Sigma
X)^{\wedge n},
$$
and by 3.8 this will imply
  $$
(D_nF)(X)\sim\Omega^n\Sigma^{\infty}X^{\wedge n}.
$$
Incidentally, the naturality of Taylor towers gives a quick way to get from the
James model of $\Omega\Sigma X$ to that splitting of
$\Sigma^{\infty}\Omega\Sigma X$. The space $JX$, free monoid on the based
space $X$, is naturally filtered by word length as an increasing union
$*=J_0X\subset J_1X\subset \dots$  in such a way that the subquotient
$J_n/J_{n-1}$ is
$X^{\wedge n}$. It follows that for each $n$ there is a natural 
fibration sequence
$$
\Sigma^{\infty}J_{n-1}X\a \Sigma^{\infty}J_n X \a \Sigma^{\infty}X^{\wedge
n}
$$
One sees by induction that $\Sigma^{\infty}J_n$
is
$n$-excisive, since it fibers over a (homogeneous)
$n$-excisive functor and the fiber is ($(n-1)$-excisive, hence) 
$n$-excisive. The
sequence above must split. Indeed, any natural fibration sequence
$$
F(X)\a G(X)\a H(X)
$$
of spectra in which
$H$ is $n$-homogeneous and $F$ is $(n-1)$-excisive must split: a
retraction from
$G$ to
$F$ (in the homotopy category of homotopy functors) is
given by the diagram
$$
\CD
F @>>> G\\
@V{p_{n-1}F} V\sim V
@VV{p_{n-1}G} V\\ P_{n-1}F @>>\sim
>  P_{n-1}G
\endCD
$$
The left arrow is an equivalence because $F$ is
$(n-1)$-excisive; the lower arrow because $P_{n-1}H$ is contractible.

\endexample

\remark{1.22 Remark}The Taylor tower construction extends easily to
functors
$$
\Cal C_{Y_1}\times\dots\times \Cal C_{Y_n}@>F>> \Cal D
$$
of several variables. Let us say that $F$ is $(d_1,\dots,d_n)$-excisive
[resp.\
$(d_1,\dots,d_n)$-reduced] if for
$1\leq j\leq n$ it is
$d_j$-excisive [resp.\ $d_j$-reduced] as a functor of the $j^{th}$ variable.
There is an
$n$-variable Taylor polynomial construction which will be denoted
$P_{d_1,\dots,d_n}F$, and which gives the universal
$(d_1,\dots,d_n)$-excisive functor under
$F$ (in the homotopy category of homotopy functors). It may be defined
either as
$P_{d_1}^{(1)}\dots P_{d_n}^{(n)}F$ where
$P_{d}^{(j)}$ is $P_d$ with respect to the $j^{th}$ variable, or 
directly as the
homotopy colimit of $(T_{d_1,\dots,d_n})^kF$ where
$$
(T_{d_1,\dots,d_n}F)(X_1,\dots,X_n)=\holim_{(U_1,\dots,U_n)\in \prod_{j}
\Cal P_0(\underline
{{d_j}+1})} F(X_1*_{Y_1}U_1,\dots,X_n*_{Y_n}U_n);
$$
these are naturally equivalent. A functor of $n$ variables will be called
multilinear if it is linear in each variable, in other words both
$(1,\dots,1)$-excisive and $(1,\dots,1)$-reduced. If
$F$ is $(1,\dots,1)$-reduced then the functor $P_{1,\dots,1}F$ is
multilinear and may be called the multilinearization of $F$. Loosely, this is
the homotopy colimit of
$$
\Omega^{{k_1}+\dots{k_n}}F(\Sigma_{Y_1}^{k_1}X_1,\dots,\Sigma_{Y_n}^{k_n}X_n)
$$
as $(k_1,\dots,k_n)@>>>(\infty,\dots,\infty)$.
\endremark

\head \S2. Delooping homogeneous functors\endhead

We will show that all homogeneous
space-valued functors arise from homogeneous spectrum-valued functors. The
main step is to show that they have natural deloopings.

Let $\Cal H_n(\Cal C,\Cal D)$ be the category of homogeneous functors
of degree $n$ from
$\Cal C$ to $\Cal D$, a full subcategory of $\Cal F(\Cal
C,\Cal D)$. The homotopy
category
$h\Cal H_n(\Cal C,\Cal D)$) is obtained from
$\Cal H_n(\Cal C,\Cal D)$ by formally inverting the
(objectwise) equivalences.

The functor
$\Cal Sp@>{\Omega^{\infty}}>>\Cal T$ preserves both weak equivalences
and cartesian cubes, and therefore composing with it yields a
functor
$$
\Cal H_n(\Cal C,\Cal Sp)@>{{\Omega^{\infty}}}>>\Cal H_n(\Cal C,\Cal T)
$$
which itself takes weak equivalences to weak equivalences.

\proclaim{2.1 Theorem}The functor
$
\Cal H_n(\Cal C,\Cal Sp)@>{{\Omega^{\infty}}}>>\Cal H_n(\Cal C,\Cal T)
$
has an inverse up to weak equivalence.
\endproclaim

\demo{Proof}The key is to get a functor
$F\mapsto BF$ from
$\Cal H_n(\Cal C,\Cal T)$ to itself such that $\Omega
BF$ is naturally equivalent to $F$. This will be given by the next lemma.

Assume for now that we have this. Then any object $F$ of $\Cal H_n(\Cal
C,\Cal T)$ yields a sequence
$\lbrace B^pF\rbrace$  of such objects
related by equivalences $B^pF\sim \Omega B^{p+1}F$. Call the resulting
spectrum-valued functor $B^{\infty}F$.
The fact that $B^{\infty}F$ is an
object of $\Cal H_n(\Cal C,\Cal Sp)$ follows easily
from the fact that each $B^pF$ is
an object of $\Cal H_n(\Cal C,\Cal T)$. The fact that the functor $B^{\infty}$
takes equivalences to equivalences follows from the fact that
each functor $B^p$ does so.
Clearly ${\Omega}^{\infty}B^{\infty}$ is naturally equivalent to the 
identity. To
check that this is also true for
$B^{\infty}{\Omega}^{\infty}$, let $F$ be an object of
$\Cal H_n(\Cal C,\Cal Sp)$ and write $F(X)=\lbrace F_q(X)\rbrace$.
The  bispectrum $\lbrace B^pF_q(X)\rbrace$ shows that the two spectra
$$
\lbrace B^pF_0(X)\rbrace = \lbrace B^p {\Omega}^{\infty}F(X)\rbrace =
B^\infty {\Omega}^{\infty}F(X)
$$
and
$$
\lbrace B^0F_q(X)\rbrace = \lbrace F_q(X)\rbrace = F(X)
$$
are naturally equivalent.
\qed\enddemo

The next result provides the desired delooping of $F$, namely $R_nF$. As usual
$\Cal C_Y$ is either
$\Cal U_Y$ or $\Cal T_Y$.

\proclaim{2.2 Lemma}Let $n>0$. If $F:\Cal C_Y\a\Cal T$ is
any reduced ($F(Y)\sim *$) homotopy functor, then up to natural equivalence
there is a fibration sequence
$$\CD
                            P_nF @>{q_nF}>> P_{n-1}F @>>> R_nF
\endCD
$$
in which the functor $R_nF$ is $n$-homogeneous.
\endproclaim

\demo{Proof}More precisely, we will obtain a
natural  diagram of homotopy functors
$$
\CD
P_nF @>{q_nF}>> P_{n-1}F\\
@A{\sim}AA @AA{\sim}A\\
\hat P_nF @>>> \tilde P_{n-1}F\\
@VVV @VVV\\
K_nF @>>> R_nF
\endCD\tag2.3
$$
in which the marked arrows are equivalences, the lower square
is cartesian,
$R_nF$ is
$n$-homogeneous, and $K_nF$ is contractible. (If $F$ were not reduced,
then in fact $K_nF$ would be equivalent to the constant functor
$F(Y)$.)

The proof is based on a close examination of the maps $q_{n,i}F$
which were used in
defining $q_nF$. The first step is to define, for each $i\geq 0$,
a diagram
$$
\CD
T_n^iF @>>> S_{n-1}^iF\\
@VVV @VVV\\
K_{n,i}F @>>> R_{n,i}F
\endCD\tag2.4(i)
$$
Roughly, this will become the lower square of (2.3) when $i$ goes to
infinity.

Define posets
$$\gather
\Cal B _n=\Cal P_0(\underline{n+1})-\lbrace\lbrace{n+1}\rbrace\rbrace\\
\Cal A _{n,i}
={\Cal P_0(\underline{n+1})}^i-{\Cal P_0(\underline{n})}^i.
\endgather
$$
Define the functor
$S_{n-1}F$ by
$$
               (S_{n-1}F)(X) = \holim_{U\in \Cal B_n} { F(X*_YU)}
$$
Note that the inclusions
$\Cal P_0(\underline{n+1})\supset \Cal B_n \supset \Cal
P_0(\underline{n})$ induce maps
$$\CD
T_nF @>>> S_{n-1}F @>{\sim}>> T_{n-1}F
\endCD\tag2.5
$$
whose composition is $q_{n,1}F$. The second map is an equivalence because
$\Cal P_0(\underline{n})$ is left cofinal in
$\Cal B_n$. Now let (2.4(i)) be obtained from
the posets and inclusions:
$$
\CD
{\Cal P_0(\underline{n+1})}^i @<<< {\Cal B_n}^i\\
@AAA @AAA\\
\Cal A_{n,i} @<<< {\Cal A_{n,i}}\cap {\Cal B_n}^i
\endCD\tag2.6(i)
$$
by forming the homotopy limit of
$(U_1,\dots,U_i)\mapsto F(X*_Y(U_1*\dots *U_i))$ over each
poset.

Diagram (2.4(i)) is cartesian by [\cite{2}; 1.9],
since $\Cal A_{n,i}$ and ${B_n}^i$ are
concave and their union is
${\Cal P_0(\underline{n+1})}^i$.
\proclaim{2.7 Claim}$K_{n,i}F\sim *$.
\endproclaim
\demo{Proof}Compare $(K_{n,i}F)(X)$, a
$\holim$ over $\Cal A_{n,i}$, with the
corresponding $\holim$ over the smaller poset
$$
\Cal A^*_{n,i}={\Cal P_0(\underline{n+1})}^i-{\Cal B_n}^i
$$
On the one hand, the comparison map is an equivalence, because
$\Cal A^*_{n,i}$ is left cofinal in $\Cal A_{n,i}$.
On the other hand, the $\holim$ of $F(X*_Y(U_1*\dots *U_i))$ over
$\Cal A^*_{n,i}$ is
contractible: for each $(U_1,\dots ,U_i) \in \Cal A^*_{n,i}$ ,
  we have $|U_j|=1$ for
some $j$, so that $U_1*\dots*U_i$ is contractible and
$$
             F(X*_Y(U_1*\dots*U_i)) \sim F(Y)    \sim *
\eqno{\qed}
$$

\proclaim{2.8 Claim}$R_{n,i}F$ is $n$-reduced.
\endproclaim
\demo{Proof}
There is a natural equivalence
$$
P_{n-1}R_{n,i}F\sim R_{n,i}P_{n-1}F
$$
by 1.7, so it will be enough if
$R_{n,i}F$ is contractible whenever $F$ is $(n-1)$-excisive (and reduced).

There is an isomorphism of posets
$$
{\Cal A_{n,i}}\cap {\Cal B_n}^i \cong
{\Cal P_0(\underline{n})}^i \times {\Cal P_0(\underline i})
$$
given by
$$
\align
                      (U_1,\dots ,U_i) &\mapsto (V_1,\dots ,V_i,W)\\
                                           V_i &= U_i - \lbrace{n+1}\rbrace\\
                                           W &= \lbrace{ i : {n+1}\in 
U_i }\rbrace
\endalign
$$
Therefore $(R_{n,i}F)(X)$ can be written as
$$
             \holim_{(V_1,\dots ,V_i,W)\in
{\Cal P_0(\underline{n})}^i \times
{\Cal P_0(\underline i})}
  { F(X*_Y(U_1*\dots*U_i)) }
$$
Analyze this as an iterated homotopy limit:
First fix $(V_2,\dots ,V_i,W)$ and take homotopy limit with
respect to $V_1$. Since $F$ is $(n-1)$-excisive, this yields,
up to equivalence,
$F(X*_Y(e_1*U_2*\dots *U_i))$ where $e_1$ is
$\lbrace{n+1}\rbrace$ or $\emptyset$
according as $1$ is or is not in $W$.
Next take the homotopy limit with respect to $V_2$, then $V_3$,
and so on through $V_i$, obtaining
$$
               (R_{n,i}F)(X) \sim
\holim_{W\in  {\Cal P_0(\underline i)}}
{ F(X*_Y(e_1*\dots *e_i)) }
$$
This is contractible because for each $W\in
{\Cal P_0(\underline i)}$ the space
${e_1*\dots *e_i}$ is contractible.
(It is a simplex of dimension $|W|-1\geq 0$.)
\qed\enddemo

The next step ought to be to take the homotopy colimit of (2.4(i))
as $i$ tends to $\infty$. This is not possible, since we do not have 
natural maps
$K_{n,i}F\a K_{n,i+1}F$
or $R_{n,i}F\a R_{n,i+1}F$. We do, however, have maps defined
up to homotopy, and with a
little care these will suffice.

Here are two variations on (2.6(i+1)):
$$\CD
{\Cal P_0(\underline {n+1})}^{i+1} @<<< {\Cal P_0(\underline
{n+1})}\times {\Cal B_n^i}\\
@AAA @AAA\\
{\Cal P_0(\underline {n+1})}\times {\Cal A_{n,i}} @<<<
{\Cal P_0(\underline {n+1})}\times ({\Cal A_{n,i}\cap \Cal B_n^i})
\endCD\tag2.9(i+1)
$$
$$\CD
{\Cal P_0(\underline {n+1})}^{i+1} @<<<
{\Cal P_0(\underline {n+1})}\times {\Cal B_n^i}\\
@AAA @AAA\\
{\Cal A_{n,i+1}} @<<<
\Cal A_{n,i+1}\cap
({\Cal P_0(\underline {n+1})}\times {\Cal B_n^i})
\endCD\tag2.10(i+1)
$$
We have maps of square diagrams of posets
$$
(2.9(i+1))\rightarrow (2.10(i+1))\leftarrow (2.6(i+1)).
$$
 From each of the three diagrams we get a square diagram of
functors by taking the $\holim$
of
$F(X*_Y(U_1*\dots *U_i))$. From (2.6(i+1)) we get (2.4(i+1)).
 From (2.9(i+1)) we get what might be called $T_n(2.4(i)))$.
From
(2.10(i+1)) we get something new; call it (2.11(i+1)). These are
related by maps (of square diagrams)
$$
\dots \a2.4(i)\a T_n(2.4(i))\leftarrow
2.11(i+1)\a2.4(i+1)\a\dots
$$
Now consider the (pointwise) homotopy colimit of this,
another square diagram. This will be the
lower half of (2.3). In view of the following, it is essentially
a limit of the cartesian
squares $2.4(i)$ and therefore it is itself cartesian:
\proclaim{2.12 Claim}The backwards arrow
$T_n(2.4(i))\leftarrow
2.11(i+1)$ is an equivalence (in all four corners of the square).
\endproclaim
\demo{Proof}In the upper corners it is an isomorphism.
In the lower corners it is induced
by the inclusions
$$\CD
{\Cal P_0(\underline {n+1})}\times
{\Cal A_{n,i}} @>>> \Cal A_{n,i+1}\\
@. @.\\
{\Cal P_0(\underline {n+1})}\times ({\Cal A_{n,i}\cap \Cal B_n^i}) @>>>
\Cal A_{n,i+1}\cap ({\Cal P_0(\underline {n+1})}\times {\Cal B_n^i})
\endCD$$
In each of the two cases the larger poset is the union of the smaller 
one, which
is concave (in the sense of [\Cite{2}; page 298]), with the concave set
$$
\Cal Q = \Cal A_{n,1}\times \Cal B_n^i
$$
and in each case the intersection is
$$
\Cal Q^0 =\Cal A_{n,1}\times(\Cal A_{n,i}\cap\Cal B_n^i)
$$
Thus by [\Cite{2}; 0.2] it will be enough if the inclusion $\Cal Q^0 \a \Cal
Q$ induces an equivalence of holims. But each
of these holims is contractible. (The argument is as in the proof of 2.7;
replace $\Cal Q$ by $\Cal A^*_{n,1}\times \Cal B_n^i$ and $\Cal Q^0$ by $\Cal
A^*_{n,1}\times(\Cal A_{n,i}\cap\Cal B_n^i)$.)
\qed\enddemo

It follows from 2.7 and 2.12 that $K_nF$ is contractible; it is the
$\hocolim$ of a diagram
$$\CD
\dots @>>> K_{n,i}F @>{t_n}>> T_nK_{n,i}F @<{\sim}<< ? @>>> K_{n,i+1}F
@>{t_n}>> \dots
\endCD$$
of contractible objects.

It follows from 2.8 and 2.12 that
$R_nF$ is $n$-reduced, being the
$\hocolim$ of the functors
$$\CD
\dots @>>> R_{n,i}F @>{t_n}>> T_nR_{n,i}F @<{\sim}<< ? @>>> R_{n,i+1}F
@>{t_n}>> \dots
\endCD$$
Moreover, $R_nF$ is $n$-excisive, by Lemma 1.9; if $\Cal X$ is any strongly
cocartesian $(n+1)$-cube in $\Cal C$ then $(R_nF)(\Cal X)$ is cartesian
because for each $i$ the map
$$
  (R_{n,i}F)(\Cal X) @>{t_n}>> (T_nR_{n,i}F)(\Cal X)
$$
factors through a cartesian cube.

Finally we construct the upper half of (2.3). Note
that
$\hat P_nF$ is the
$\hocolim$ of
$$
\CD
\dots @>>> T_n^iF @>{t_n}>> T_n^{i+1}F @<{=}<< T_n^{i+1}F @>{=}>> 
T_n^{i+1}F @>{t_n}>>
\dots
\endCD$$
Eliminating the identity maps we obtain an equivalence from $\hat 
P_nF$ to the $\hocolim$
of
$$\CD
\dots @>>> T_n^iF @>{t_n}>> T_n^{i+1}F @>{t_n}>> \dots
\endCD$$
which is $P_nF$. This is the upper left arrow in (2.3). Likewise $\tilde
P_{n-1}F$ is the
$\hocolim$ of
$$
\dots \lra S_{n-1}^iF \buildrel{t_n}\over\lra T_nS_{n-1}^iF 
\buildrel{=}\over\longleftarrow T_nS_{n-1}^iF 
\lra S_{n-1}^{i+1}F\buildrel{t_n}\over\lra 
 \dots
$$
  and so has an equivalence to the $\hocolim$ of the upper row in
$$\CD
\dots @>>> S_{n-1}^iF @>{t_n}>> T_nS_{n-1}^iF @>>> S_{n-1}^{i+1}F @>>> \dots\\
@. @V\sim VV @V\sim VV @V\sim VV  @.\\
\dots @>>> T_{n-1}^iF @>{t_n}>> T_nT_{n-1}^iF @= T_{n-1}^{i+1}F @>>> \dots
\endCD$$
This in turn has an equivalence (2.7) to the $\hocolim$ of the lower
row, and hence to $P_{n-1}F$. The composed map is the upper right arrow in
(2.3). The square commutes. This concludes the proof of (2.2).
\qed\enddemo

\head\S3.Symmetric multilinear functors\endhead

Let $\Cal C@>{\Delta}>>\Cal C ^n$ be the diagonal functor. It was shown
in [\Cite{2}; 3.4] that the composed functor $F\circ\Delta$ is
$(d_1+\dots +d_n)$-excisive if the functor $\Cal C ^n@>F>>\Cal
D$ is $d_j$-excisive in the $j^{th}$ variable, or in the language of 1.22
$(d_1,\dots,d_n)$-excisive. (The latter is easier to write but harder 
to pronounce.)
In particular
$F\circ\Delta$ is
$n$-excisive if $F$ is $(1,\dots,1)$-excisive.

\proclaim{3.1 Lemma}If
$\Cal C ^n@>F>>\Cal D$ is a $(1,\dots,1)$-reduced homotopy functor, then
$F\circ\Delta$ is $n$-reduced. Thus $L\circ\Delta$ is
$n$-homogeneous if $\Cal C ^n@>L>>\Cal D$ is
multilinear.
\endproclaim

The proof resembles that of 1.8. The key is:

\proclaim{3.2 Lemma}If $\Cal C ^n@>F>>\Cal D$ is a $(1,\dots,1)$-reduced
homotopy functor, then for any $X\in \Cal C$ the map
$$
(F\circ\Delta)(X) 
@>t_{n-1}(F\circ\Delta)>>T_{n-1}(F\circ\Delta)(X)
$$ 
factors through 
a (weakly) contractible
object.
\endproclaim

\demo{Proof of 3.1, assuming 3.2}
Let $X\in \Cal C$. The object $P_{n-1}(F\circ\Delta)(X)$ is
defined as the sequential $\hocolim$ of a diagram whose
$(i+1)\st$ map is
$$
\CD
T_{n-1}^i(F\circ\Delta)(X) @>t_{n-1}T_{n-1}^i(F\circ\Delta)>>
T_{n-1}^{i+1}(F\circ\Delta)(X)
\endCD
$$
It is enough if each of these maps factors through a weakly
contractible object. Lemma 3.2 takes care of this, not only for the first map
but also for the others, since the functor
$T_{n-1}(F\circ\Delta)$ is equal to
$G\circ\Delta$ for a functor $G$ that satisfies the same hypotheses as
$F$.
  \qed\enddemo
\demo{Proof of 3.2}
Let $\Cal E$ and $\Cal E^*$ be as in the proof of 1.9. Since $\Cal E$
contains the image of the diagonal map
$\Cal P_0(\underline{n})@>>>{\Cal P_0(\underline{n})}^n$, the
map
$t_{n-1}(F\circ\Delta)(X)$ can be factored
$$
\align
F(X,\dots,X) &\cong F(X*_Y\emptyset,\dots,X*_Y\emptyset)\\
                    &\a \holim_{(U_1,\dots ,U_n)\in \Cal E }{
F(X*_YU_1,\dots,X*_YU_n)} \\
                    &\a \holim_{U\in \Cal P_0(\underline{n})}
F(X*_YU,\dots ,X*_YU)
  \endalign
$$
But
$$
\holim_{(U_1,\dots ,U_n)\in \Cal E }{
F(X*_YU_1,\dots,X*_YU_n)}\sim \holim_{(U_1,\dots ,U_n)\in \Cal
E^* }{  F(X*_YU_1,\dots,X*_YU_n)}
$$
and this is a $\holim$ of weakly contractible objects:
if $(U_1,\dots,U_n)\in \Cal E^*$ then some $U_s$ is a one-element
set, giving
$X*_YU_s\sim Y$ and
$$F(X*_YU_1,\dots,X*_YU_n)\sim *.\eqno{\qed}
$$

The functor
$L:\Cal C^n\a\Cal D$ is
\it symmetric\rm\ if it has additional
structure  consisting of isomorphisms
$L(\pi):L(X_1,\dots ,X_n)\a L(X_{\pi(1)},\dots ,X_{\pi (n)})$
for all $\pi\in\Sigma_n$, with
$L(\sigma\circ\pi)=L(\pi)\circ L(\sigma)$. (In other words $L$ is extended
from $\Cal C^n$ to a wreath product category.)
If $L$ is symmetric and multilinear then the homogeneous functor
$L\circ
\Delta$ has a
$\Sigma_n$-action. In the case of spectrum-valued functors (that is, when
$\Cal D=\Cal Sp$), the object
$$
(\Delta_nL)(X) = L(X,\dots,X)_{h\Sigma_n}
$$
is then again an $n$-homogeneous functor of $X$, by 1.7(5), since homotopy
orbit space (or spectrum) is a special case of $\hocolim$. We are headed toward
proving that \it all\rm\
$n$-homogeneous functors arise in just this way.

The inverse of $\Delta_n$ will be provided by a construction called the
$n^{th}$ cross-effect, which takes a homotopy
functor $F$ and produces a symmetric homotopy functor $cr_nF$ of
$n$ variables.
To see how this inverse
construction should go, we recall an algebraic analogue. $\Delta_n$ is
analogous to the construction in algebra which uses a symmetric multilinear
function $l(x_1,\dots,x_n)$ to make a homogeneous function
$f(x)=l(x,\dots,x)/n!$. For example, if
$f$ is a degree two polynomial then from the bilinear form
$$
l(x_1,x_2)=f(x_1+x_2)-f(x_1)-f(x_2)+f(0)
$$
we recover the purely quadratic part of $f$ as $l(x,x)/2$.
The second cross-effect $cr_2F$ of a functor $F$
of based spaces will take the based spaces $X_1$ and $X_2$ to the
total homotopy fiber of the diagram
$$\CD
                           F(X_1\vee X_2) @>>> F(X_1)\\
                            @VVV @VVV\\
                               F(X_2) @>>>  F(*)
\endCD$$
Strictly speaking, in order to make $cr_2F$ preserve weak equivalences
we must first replace each $X_i$ by an equivalent object
having nondegenerate basepoint, perhaps by using \lq\lq whiskers\rq\rq.
In the general case (functors of fiberwise based spaces over $Y$) the
wedge $X_1\vee X_2$ is replaced by a categorical sum in
$\Cal T_Y$. Again, before forming the sum we should use fiberwise whiskers,
replacing each object
$X$ by the mapping cylinder of the  structural coretraction
$Y@>>>X$.

Thus the $n^{th}$ cross-effect is defined as follows: Let $\Cal T_Y@>F>>\Cal D$
be a homotopy functor, $\Cal D=\Cal T$ or $\Cal Sp$. For objects
$(X_1,\dots,X_n)$ of
$\Cal T_Y$, let
$\Cal S(X_1,\dots ,X_n)$ be the evident
$n$-cube taking $\underline{n}-T$ to the (whiskered) sum, in $\Cal T_Y$, of
the objects $X_s$ for ${s\in T}$. Define the cross-effect by first applying
$F$ to the cube and then taking the total homotopy fiber:
$$
(cr_n F)(X_1,\dots ,X_n)=\tfiber F(\Cal S(X_1,\dots ,X_n)).
$$
It is easy to see that $cr_n F$ is a homotopy functor (in each
variable), and symmetric and $(1,\dots,1)$-reduced.

The first cross-effect is the
\lq\lq reduced functor\rq\rq: 
$$
(cr_1 F)(X)=\fiber(F(X)@>>>F(Y)).
$$
The $0^{th}$ cross-effect,  should we ever need it, is a functor
of no variables: the object $F(Y)$.

\proclaim{3.3 Proposition}If $F$ is $n$-excisive
then for $0\leq m\leq n$ the functor $cr_{m+1} F$ is
$(n-m)$-excisive in each variable.
In particular, the $n^{th}$ cross-effect of an $n$-excisive functor
is symmetric multilinear and the $n^{th}$ cross-effect of an
$(n-1)$-excisive functor is trivial (equivalent to a point).
\endproclaim

\demo{Proof}Induction on $m$. The case $m=0$ is clear. To pass
from
$m-1$ to $m$, write
$$
(cr_{m+1} F)(X_1,\dots X_m,A)\cong (cr_{m}F_{+A})(X_1,\dots ,X_{m}),
$$
where
$F_{+A}(X)=\hofiber(F(X+A)@>>>F(X))$ and $+$ denotes whiskered sum in
$\Cal T_Y$. Use the rather obvious fact that
$F_{+A}$ is
$(n-1)$-excisive if $F$ is $n$-excisive.
\qed\enddemo

In general, then,
$$
\align
cr_nD_nF&=cr_n\hofiber(P_nF@>>>P_{n-1}F)\\
&\sim\hofiber(cr_nP_nF@>>>cr_nP_{n-1}F)\\
&\sim cr_nP_nF.
\endalign
$$
Thus if $F$ is $n$-excisive then $cr_nF$ is naturally equivalent to $cr_nD_nF$;
the
$n^{th}$ cross-effect of an $n$-excisive functor
\lq\lq sees\rq\rq\ only the $n$-homogeneous part.
The following simple result suggests that it sees it quite clearly,
at least in the spectrum-valued case.

\proclaim{3.4 Proposition}
If $F$ and $G$ are $n$-homogeneous functors $\Cal T_Y @>>> \Cal
Sp$, then any (natural) map $F@>>>G$ that induces an
equivalence $cr_nF@>>>cr_nG$ must be an equivalence
itself.
\endproclaim
\demo{Proof}Let $H$ be the homotopy fiber of $F@>>>G$. Thus $H$ is
$n$-homogeneous and
$cr_n H$ is the homotopy fiber of $cr_nF@>>>cr_nG$. Since a map of spectra
must be an equivalence if its homotopy fiber is contractible, we have only to
show that, for an $n$-homogeneous functor $\Cal T_Y @>H>>
\Cal Sp$,
$cr_n H\sim *$ implies $H\sim *$. In fact we will do a little better:
assuming only that $H$ is $n$-excisive and $cr_n H\sim *$, we will show that
$H$ is
$(n-1)$-excisive.

If $\Cal X$ is any strongly cocartesian
$n$-cube in $\Cal T_Y$, we must show that $H(\Cal X)$ is
cartesian. Because we are dealing with spectra, it will be enough if
we show that $\tfiber  H(\Cal X)\sim *$.

By assumption, this holds in the case where $\Cal X$ is the cube
$\Cal S(X_1,\dots,X_n)$ for objects  $X_1,\dots,X_n$ of $\Cal T_Y$.

It also holds
for the related cube, call it $\Cal S ^* (X_1,\dots ,X_n)$, which has the
same objects but with reversed arrows, sending $T$ rather than $\underline n-T$
to the sum of the
objects $X_s$ for ${s\in T}$. In fact we have
$$\align
\tfiber H(\Cal S ^* (X_1,\dots ,X_n))
&\sim
{\Omega}^n\tfiber H(\Cal S(X_1,\dots,X_n))\\&=
{\Omega}^n (cr_n H)(X_1,\dots ,X_n).\endalign
$$
This takes care of all cubes $\Cal X$ in which $\Cal X (\emptyset)\sim Y$,
because such a cube is naturally equivalent to $\Cal S ^* (\Cal X(\lbrace
1\rbrace),\dots ,\Cal X(\lbrace
n\rbrace))$ (see [\Cite{2}; 2.2]).

Given an arbitrary strongly cocartesian cube $\Cal X$, put
$$\Cal X'(T)=
\hocolim (Y@<<<\Cal X(\emptyset)\a \Cal X (T)).
$$
The obvious map of $n$-cubes
$\Cal X\a \Cal X'$ is a strongly cocartesian $(n+1)$-cube, so the 
resulting cube
$H(\Cal X)\a H(\Cal X')$ is cartesian. The cube
$H(\Cal X')$ is cartesian by the case already treated. Therefore by 
[\Cite{1}; 1.6]
$H(\Cal X)$ is cartesian.
\qed\enddemo

Let $\Cal L_n(\Cal C,\Cal D)$ be the category
of symmetric multilinear functors $L:{\Cal C}^n\a\Cal D$. The maps
are the natural maps that respect the symmetry. As usual the homotopy
category $h\Cal L_n(\Cal C,\Cal D)$ means the category obtained by
inverting the (objectwise) equivalences.

If $\Cal D$ is either $\Cal Sp$ or $\Cal T$ then there is the cross-effect
functor
$$
\Cal H_n(\Cal T_Y,\Cal D)@>{cr_n}>>\Cal L_n(\Cal T_Y,\Cal D)
$$
If $\Cal C$ is either $\Cal U_Y$ or $\Cal T_Y$ then there is the functor
$$
\Cal L_n(\Cal C,\Cal Sp)@>{\Delta_n}>>\Cal H_n(\Cal C,\Cal Sp)
$$
going the other way. Both $cr_n$ and $\Delta_n$ preserve weak equivalences
and so induce functors on homotopy categories.

\proclaim{3.5 Theorem}The functors
$$\align
&\Cal H_n(\Cal T_Y,\Cal Sp)@>{cr_n}>>\Cal L_n(\Cal T_Y,\Cal Sp)\\
&\Cal L_n(\Cal T_Y,\Cal Sp)@>{\Delta_n}>>\Cal H_n(\Cal T_Y,\Cal Sp)
\endalign$$
are mutual inverses up to natural (weak) equivalence.
\endproclaim

\demo{Proof of 3.5}
To prove that the composition
$$
\Cal L_n(\Cal T_Y,\Cal Sp)@>{\Delta_n}>>
\Cal H_n(\Cal T_Y,\Cal Sp)@>{cr_n}>>
\Cal L_n(\Cal T_Y,\Cal Sp)
$$
is equivalent to the identity, we look to the algebraic analogue. If $l$ is a
symmetric multilinear function of $n$ variables and $f$ is the homogeneous
polynomial
$$
f(x)=l(x,\dots,x)/n!,
$$
then $l$ can be recovered from $f$. It is given by
$$
l(z_1,\dots,z_n)=f(z_1+\dots+z_n)-f(z_1\dots+z_{n-1})-\dots+(-1)^nf(0)
$$
(an alternating sum of $2^n$ terms). One sees this, of course, by
expanding the expression
$$
l(x_1+\dots +x_n,\dots,x_1+\dots+x_n)
$$
as
a sum of $n^n$ terms, cancelling all except the permutations of
$l(x_1,\dots,x_n)$ and dividing by $n$ factorial.

Here is a corresponding categorical argument: We have
$$\align
(cr_n\Delta_nL)(X_1,\dots,X_n)&=\tfiber((\Delta_nL)\circ\Cal
S(X_1,\dots,X_n))\\
&\cong\tfiber(\underline n - T\mapsto L(\coprod _{s\in T} X_s,\dots,\coprod
_{s\in T} X_s)_{h\Sigma_n})\\
&\sim (\tfiber(\Cal X))_{h\Sigma_n},\\
\endalign$$
where $\Cal X(\underline n - T)=L(\coprod _{s\in T} X_s,\dots,\coprod
_{s\in T} X_s)$ and $\coprod$ denotes whiskered sum in $\Cal T_Y$. The obvious
equivalence
$$
\Cal X(\underline n - T)
@>{\sim}>>\prod_{\underline{n} @>{\pi}>> T}L(X_{\pi(1)},\dots ,X_{\pi (n)})
$$
is natural with respect to $T$. It can be interpreted as an equivalence
of cubes $\Cal X@>>>\prod_{\underline{n} @>{\pi}>>\underline{n}}\Cal Y_{\pi}$
where
$$
\align
\Cal Y_{\pi}(\underline{n}-T)&=L(X_{\pi (1)},\dots ,X_{\pi (n)}), \quad if
\quad\pi(\underline{n})\subset T\\
\Cal Y_{\pi}(\underline{n}-T)&=* , \quad \text{otherwise.}
\endalign
$$
For any $\pi$ that is not a permutation and therefore not surjective, the cube
$\Cal Y_{\pi}$ is cartesian. (Viewed in one way it is an isomorphism of
$(n-1)$-cubes: if
$s\notin\pi(\underline{n})$ then all maps $\Cal Y_\pi(T)@>>>\Cal
Y_\pi(T\cup\lbrace s\rbrace)$ are isomorphisms.) For any permutation $\pi$
we have
$$\tfiber (\Cal Y_\pi)\cong L(X_{\pi(1)},\dots,X_{\pi(n)})$$
Therefore
$$\tfiber(\Cal X)@>{\sim}>>
\prod_{\pi\in\Sigma_n}L(X_{\pi(1)},\dots,X_{\pi(n)})$$ This map respects the
symmetry if the group is made to permute the factors of the 
right-hand side, and so it leads to an equivalence of homotopy orbit 
spectra
$$
(cr_n\Delta_nL)(X_1,\dots,X_n)\sim
(\prod_{\pi\in\Sigma_n}L(X_{\pi(1)},\dots,X_{\pi(n)}))_{h\Sigma_n}\sim
L(X_1,\dots,X_n)
$$
Note for future reference that an explicit inverse equivalence
$$
L(X_1,\dots,X_n)@>{\theta}>>(cr_n\Delta_nL)(X_1,\dots,X_n)
$$
is the map of total homotopy fibers induced by an obvious map of cubes
$$
\Cal Y_1 @>>> \Cal X@>>> \Cal X_{h\Sigma_n},
$$
where $\underline{n}@>1>>\underline{n}$ is the identity map and the map $\Cal
Y_1(\emptyset) @>>> \Cal X(\emptyset)$ is the map
$$
L(X_1,\dots,X_n)@>i>>L(Z,\dots,Z)
$$
induced by the inclusions $X_j@>>>Z=\coprod_{1\leq j\leq n}X_j$. The
following diagram commutes:
$$
\CD
L(X_1,\dots,X_n)@>i>>L(Z,\dots,Z)\\
@V{\theta}VV @VVV\\
(cr_n\Delta_nL)(X_1,\dots,X_n)@>{\epsilon}>>L(Z,\dots,Z)_{h\Sigma_n}.
\endCD
$$
Here, and throughout this proof, $\epsilon$ denotes the projection from the
total homotopy fiber of a cubical diagram to the \lq\lq 
initial\rq\rq\ object in
the diagram.

The remaining task is to exhibit an equivalence
$$
\Delta _ncr_n F@>{\gamma}>>F
$$
for any $n$-homogeneous $\Cal T_Y @>F>> \Cal Sp$. In fact we will define
$\gamma$ for any
$n$-excisive $F$ and then show that it induces an equivalence
$$
cr_n\Delta _n cr_n F@>{cr_n(\gamma)}>>cr_nF.
$$
By 3.4 this suffices.

To define $\gamma$ we use a map $\hat\gamma$, defined for $Z\in \Cal T_Y$
as the composition
$$
(cr_n F)(Z,\dots,Z) @>{\epsilon}>> F(\coprod _{1\leq j\leq n}Z) @>{F(f)}>> F(Z)
$$
where $f$ is the  \lq\lq fold\rq\rq\ map which takes each copy of $Z$
identically to $Z$. The map
$\hat\gamma$ is equivariant with respect to the obvious $\Sigma _n$-actions.
(The action on
$F(Z)$ is trivial.) Define
$\gamma$ as the composition
$$\align
(\Delta_n cr_n F)(Z)& = {((cr_n F)(Z,\dots,Z))}_{h\Sigma_n}\\
&@>{{\hat\gamma}_{h{\Sigma_n}}} >> F(Z)_{h\Sigma_n} \cong F(Z)\wedge
(B\Sigma_n)_+ @>>> F(Z)\endalign
$$
where the last arrow is induced by the nontrivial based map $(B\Sigma_n)_+\a
S^0$.

To see that $cr_n(\gamma)$ is an equivalence we examine the composition
$$
cr_nF@>{\theta}>{\sim}>cr_n\Delta_ncr_nF@>{cr_n(\gamma)}>>cr_nF.
$$
It will be enough if it coincides with the identity, at least on the
level of homotopy groups. In fact it will be enough if the composition
$$
(cr_nF)(X_1,\dots,X_n)@>{cr_n(\gamma)\circ\theta}>>(cr_nF)(X_1,\dots,X_n)
@>{\epsilon}>>F(Z)
$$
is equal to
$\epsilon$, since $\epsilon$ is a split injection. (Here again
$Z=\coprod_{1\leq j\leq n}X_j$.)

That it is is a direct consequence
of the following facts: The diagram
$$
\CD
(cr_nF)(X_1,\dots,X_n)@>i>>(cr_nF)(Z,\dots,Z)\\
@V{\theta}VV @VVV\\
(cr_n\Delta_ncr_nF)(X_1,\dots,X_n)@>{\epsilon}>>(cr_nF)(Z,\dots,Z)_{h\Sigma_n}\\
@V{cr_n(\gamma)}VV
@V{\gamma}VV\\(cr_nF)(X_1,\dots,X_n)@>{\epsilon}>>F(Z);
\endCD
$$
commutes; the composition of the right-hand vertical maps above is
$\hat\gamma$, which is also the composition of the lower horizontal maps below;
the diagram
$$
\CD
cr_nF(X_1,\dots,X_n)@>{\epsilon}>>F(Z)\\
@ViVV @VF(D)VV\\
cr_nF(Z,\dots,Z) @>{\epsilon}>>F(\coprod_{1\leq j\leq n}Z) @>F(f)>> F(Z),
\endCD
$$
commutes, where $Z@>D>>\coprod_{1\leq j\leq n}Z$ sends the
copy of
$X_i$ in
$Z$ to the copy of $X_i$ in the $i^{th}$ copy of $Z$; and the composition
$f\circ D$ is the identity map.
\qed\enddemo

\proclaim{3.6 Corollary}The functor
$\Cal H_n(\Cal T_Y,\Cal T)@>{cr_n}>>\Cal L_n(\Cal T_Y,\Cal T)$
has an inverse up to weak equivalence.
\endproclaim

\demo{Proof of 3.6}Using the commutative
diagram
$$\CD
\Cal H_n(\Cal T_Y,\Cal T) @<{{\Omega^{\infty}}_{*}}<< \Cal H_n(\Cal
T_Y,\Cal Sp) \\
@V{cr_n}VV @V{cr_n}VV\\
\Cal L_n(\Cal T_Y,\Cal T) @<{{\Omega^{\infty}}_{*}}<< \Cal L_n(\Cal
T_Y,\Cal Sp)
\endCD
$$
this will follow from 3.5,  2.1, and the next result.
\qed\enddemo

\proclaim{3.7 Proposition}The functor
$
\Cal L_n(\Cal C,\Cal Sp)@>{{\Omega^{\infty}}_*}>>\Cal L_n(\Cal C,\Cal T)
$
has an inverse up to weak equivalence.
\endproclaim

\demo{Proof}As in proving 2.1, we need a delooping functor from
$\Cal L_n(\Cal C,\Cal T)$ to itself. This is much easier than 2.2.

We refer to Remark 1.22 for notation. If $L$ is symmetric
multilinear then there are natural equivalences
$$
                   L(X_1,\dots ,X_n)@>{\sim}>>(T_{1,\dots,1}L)(X_1,\dots
,X_n)@<{\sim}<<
  {\Omega}^n L(\Sigma_YX_1,\dots ,\Sigma _YX_n).
$$
We will need these equivalences to respect the
$\Sigma_n$-symmetry. In order for this to be true,  permutations of the loop
coordinates must be built into the symmetry of the last expression; a better
expression is ${\Omega}^{V_n} L(\Sigma_YX_1,\dots ,\Sigma _YX_n)$, where
$V_n$ is the standard
$n$-dimensional representation of
$\Sigma_n$.  Since this contains a trivial one-dimensional representation,
we have what we need: if $V_n=\Bbb R \oplus \bar V_n$ then
$L\sim\Omega BL$ where $BL$ is defined by
$$
(BL)(X_1,\dots ,X_n) = {\Omega}^{\bar V_n} L(\Sigma _YX_1,\dots,\Sigma
_YX_n).\eqno{\qed}
$$

We can now justify an assertion made in 1.20. Let $\Cal D$ be either 
$\Cal T$ or
$\Cal Sp$.

\proclaim{3.8 Corollary}If $\Cal T_Y@>F>>\Cal D$ is $n$-homogeneous then
$F$ is determined by
$F\circ\Sigma_Y$.
\endproclaim

\demo{Proof}According to 3.5 and 3.6 $F$ is determined by the cross-effect
$cr_nF$, and this satisfies
$$
(cr_nF)(\Sigma_YX_1,\dots,\Sigma_YX_n)\cong
(cr_n(F\circ\Sigma_Y))(X_1,\dots,X_n).
$$
On the other hand, for any symmetric multilinear functor $L$ we have
$$
L(X_1,\dots,X_n)\sim{\Omega}^{V_n}L(\Sigma_YX_1,\dots,\Sigma_YX_n),
$$
so that
$$
(cr_nF)(X_1,\dots,X_n)
\sim{\Omega}^{V_n}(cr_nF)(\Sigma_YX_1,\dots,\Sigma_YX_n)$$
$$\hbox to 2in{\hss}
\cong{\Omega}^{V_n} cr_n(F\circ\Sigma_Y)(X_1,\dots,X_n)
\eqno{\qed}$$

\head\S4.The role of the base point\endhead

The cross-effect construction applies to
functors $\Cal T_Y@>>>\Cal D$ but not to functors $\Cal U_Y@>>>\Cal D$.
In spite of this we now show that the classification of homogeneous functors
extends without change to the
$\Cal U_Y$ case.

Let
$\Cal T_Y@>{\phi}>>\Cal U_Y$ be the forgetful functor. Because
$\phi$ preserves equivalences and cocartesian square diagrams, composition
with
$\phi$ yields functors
$$
\align
\Cal H_n(\Cal U_Y,\Cal D)&@>{\phi^*}>>\Cal H_n(\Cal T_Y,\Cal D)\\
\Cal L_n(\Cal U_Y,\Cal D)&@>{\phi^*}>>\Cal L_n(\Cal T_Y,\Cal D).
\endalign
$$
\proclaim{4.1 Theorem}If $\Cal D$ is $\Cal T$ or $\Cal Sp$ then both of the
two functors $\phi^*$ above have inverses up to weak equivalence.
\endproclaim

\proclaim{4.2 Corollary}The functor $\Cal L_n(\Cal U_Y,\Cal
Sp)@>{\Delta_n}>>\Cal H_n(\Cal U_Y,\Cal Sp)$ has an inverse up to weak
equivalence.
\endproclaim

\demo{Proof of 4.2}Use 3.5, 4.1, and the
diagram
$$
\CD
\Cal L_n(\Cal U_Y,\Cal Sp)@>{\Delta_n}>>\Cal H_n(\Cal U_Y,\Cal Sp)\\
@V{\phi^*}VV @V{\phi^*}VV\\
\Cal L_n(\Cal T_Y,\Cal Sp)@>{\Delta_n}>>\Cal H_n(\Cal T_Y,\Cal Sp).
\endCD
$$
\vglue -14mm \hbox{}\qed

\demo{Proof of 4.1}We can assume $\Cal D=\Cal Sp$, since by 2.1 and 3.7 the
spectrum-valued case of 4.1 implies the space-valued case.

Let $\Cal U_Y @>{\psi}>> \Cal T_Y$ be the left
adjoint of $\phi$, so that if $X$ is a space over $Y$ then $\psi(X)$ is the
disjoint union of $X$ and $Y$ viewed as a fiberwise based space over $Y$. Like
$\phi$, this functor preserves equivalences and cocartesian square diagrams.
Therefore if
$\Cal T_Y@>F>>\Cal Sp$ is $n$-excisive then $\Cal U_Y@>{\psi^*F}>>\Cal
Sp$ will be $n$-excisive. Unlike
$\phi$, $\psi$ does not preserve the final object $Y$, so
$\psi^*F$ need not be homogeneous if $F$ is. We will show that
$$
\Cal H_n(\Cal T_Y,\Cal Sp)@>{D_n\circ\psi^*}>>\Cal H_n(\Cal
U_Y,\Cal Sp)
$$
is an
inverse, up to weak equivalence, for
$$
\Cal H_n(\Cal U_Y,\Cal Sp)@>{\phi^*}>>\Cal H_n(\Cal T_Y,\Cal Sp).
$$
The unit map $X@>>>\phi\psi X$ of the adjoint pair
induces a map
$F@>>>\psi^*\phi^*F$.
For any morphism $X@>>>X'$ in $\Cal U_Y$ the square diagram
$$\CD
X @>>>\phi\psi X\\
@VVV @VVV\\
X' @>>>\phi\psi X'
\endCD
$$
is cocartesian. Therefore if $\Cal X$ is any strongly
cocartesian
$n$-cube in $\Cal U_Y$ the unit yields a strongly cocartesian $(n+1)$-cube
$$
\Cal X@>>>\phi\psi\Cal X.
$$
It follows that if $\Cal U_Y@>F>>\Cal Sp$ is $n$-excisive then the fiber of
$F@>>>\psi^*\phi^*F$ is
$(n-1)$-excisive. Thus if $F$ is
$n$-homogeneous there are natural equivalences
$$
F@>{\sim}>> P_nF @<{\sim}<< D_nF @>{\sim}>> D_n\psi^*\phi^*F
$$
(The last map is an equivalence because its fiber is contractible. This
implication relies on the fact that these are spectrum-valued functors.)

The counit $\psi\phi X @>>> X$, like the unit, yields a
cocartesian square for every map and a strongly cocartesian $(n+1)$-cube for
every strongly cocartesian $n$-cube. Therefore for any $n$-homogeneous
$\Cal U_Y@>F>>\Cal Sp$ it yields equivalences
$$
F@>{\sim}>> P_nF @<{\sim}<< D_nF @<{\sim}<< D_n\phi^*\psi^*F \cong
\phi^* D_n\psi^*F.
$$
This completes the proof in the homogeneous case.

We sketch the proof in the
multilinear case, which is much the same. The inverse of
$$
\Cal L_n(\Cal U_Y,\Cal Sp)@>{\phi^*}>>\Cal L_n(\Cal T_Y,\Cal Sp)
$$
is
$$
\Cal L_n(\Cal T_Y,\Cal Sp)@>{r\circ\psi^*}>>\Cal L_n(\Cal
U_Y,\Cal Sp),
$$
where $r$ is the operation of \lq\lq
reducing\rq\rq\ a symmetric homotopy functor of $n$ variables in all variables
simultaneously. For example, when $n=2$ then $(rL)(X_1,X_2)$ is the total
homotopy fiber of
$$\CD
L(X_1,X_2) @>>> L(X_1,Y)\\
@VVV @VVV\\
L(Y,X_2) @>>> L(Y,Y).
\endCD$$
The key point is that the maps $L@>>>\phi^*\psi^*L$ and
$\psi^*\phi^*L@>>>L$ induce equivalences
$rL@>>>r\phi^*\psi^*L\cong\phi^*r\psi^*L$ and
$r\psi^*\phi^*L@>>>rL$. For example, in the
composition
$$
L(X_1,X_2)@>>>L(X_1,\psi\phi X_2)@>>>L(\psi\phi X_1,\psi\phi X_2)
$$
the fiber of the first map becomes contractible upon reducing with respect to
$X_2$ while the fiber of the second map becomes contractible upon reducing
with respect to $X_1$, and so they both become contractible upon applying $r$.
\qed\enddemo

  The proof of 4.1 suggests a variant of the notion of multilinear functor,
related to it as unreduced homology is related to reduced homology.

\definition{4.3 Definition}A functor $\Cal U_Y^n@>L>>\Cal Sp$ is
unreduced-multilinear if it is $1$-excisive in each variable and if 
it satisfies
$L(X_1,\dots,X_n)\sim *$ whenever $X_s$ is the initial object
$\emptyset$ for some $s$. The category of symmetric unreduced-multinear
functors is
$\hat\Cal L_n(\Cal U_Y,\Cal Sp)$.
\enddefinition

Since $\psi(\emptyset)\cong Y$, the functor $\psi^*$ maps $\Cal L_n(\Cal
T_Y,\Cal Sp)$ into $\hat\Cal L_n(\Cal
U_Y,\Cal Sp)$.

\proclaim{4.4 Proposition}The functor $\Cal L_n(\Cal
T_Y,\Cal Sp)@>{\psi^*}>>\hat\Cal L_n(\Cal
U_Y,\Cal Sp)$ has an inverse up to weak equivalence.
\endproclaim

\demo{Proof}The proof of 4.1 shows that $\phi^*\circ r$ is an inverse.
\qed\enddemo

\remark{4.5 Warning}Although Definition 4.3 could be extended verbatim to
functors into based spaces, the resulting category
$\hat\Cal L_n(\Cal U_Y,\Cal T)$ would not have the expected
property: the corresponding variant of Proposition 4.4 would be false. An
instructive example is obtained by adding a disjoint basepoint to the excisive
functor $\Cal U@>{P_1J}>>\Cal U$  mentioned near the end of the introduction.
\endremark

\remark{4.6 Remark}3.8 is valid for functors of unbased spaces, in view of 4.1
and the fact that $\phi$ commutes with $\Sigma_Y$ up to isomorphism.
\endremark

\head \S5. The $n^{th}$ differential and the $n^{th}$ derivative\endhead

We can summarize the main results of Sections 2 through 4 by saying
that the following eight categories of functors are equivalent at the
homotopy category level:
$$\CD
\Cal H_n(\Cal U_Y,\Cal Sp) @>{{\Omega^{\infty}}_{*}}>> \Cal H_n(\Cal
U_Y,\Cal T) \\
@V{\phi^*}VV @V{\phi^*}VV\\
\Cal H_n(\Cal T_Y,\Cal Sp) @>{{\Omega^{\infty}}_{*}}>> \Cal H_n(\Cal
T_Y,\Cal T)
\endCD
$$
$$\CD
\Cal L_n(\Cal U_Y,\Cal Sp) @>{{\Omega^{\infty}}_{*}}>> \Cal L_n(\Cal
U_Y,\Cal T) \\
@V{\phi^*}VV @V{\phi^*}VV\\
\Cal L_n(\Cal T_Y,\Cal Sp) @>{{\Omega^{\infty}}_{*}}>> \Cal L_n(\Cal
T_Y,\Cal T)
\endCD
$$
In addition to the arrows displayed, there is also $cr_n$ from the upper
square to the lower in each of the
two $(\Cal T_Y,-)$ cases and $\Delta_n$ from the lower square to the upper in
each of the
$(-,\Cal Sp)$ cases. We have explicitly inverted enough of these arrows to show
that all of them are invertible: We inverted the four called
$\Omega^{\infty}_*$ in 2.1 and 3.7; we inverted the left hand $\phi^*$ of each
square in 4.1; and in 3.1 we connected the two squares by showing that $cr_n$
and $\Delta_n$ are inverses in the $(\Cal T_Y,\Cal Sp)$ case.

Let $\Cal C_Y@>F>>\Cal D$ be a homotopy functor from spaces over $Y$
to either based spaces or spectra, and suppose that we wish to describe the
homogeneous layer $D_nF$ in its Taylor tower for some $n\geq 1$.  By the
results above, knowing
$D_nF$ is the same as knowing a certain symmetric multilinear functor. This
will be called the $n$-fold differential of $F$.

\definition {5.1 Definition}The object in $\Cal
L_n(\Cal C_Y,\Cal D)$ corresponding to the homogeneous functor $D_nF\in
\Cal H_n(\Cal C_Y,\Cal D)$ is called the $n$-fold differential of $F$
and is denoted by $D^{(n)}F$.
\enddefinition

Specifically, $D^{(n)}F$ determines $D_nF$ by the rule:
$$
(D_nF)(X) \sim (D^{(n)}F)(X,\dots,X)_{h\Sigma_n}
$$
in the case when $\Cal D=\Cal Sp$, or
$$
(D_nF)(X) \sim
\Omega^{\infty}((B^{\infty}D^{(n)}F)(X,\dots,X)_{h\Sigma_n})
$$
in the case when $\Cal D=\Cal T$.  Here $B^{\infty}$ is the inverse (up to
natural weak equivalence) to
$$
L_n(\Cal C_Y,\Cal Sp)@>{\Omega^{\infty}_*}>>L_n(\Cal C_Y,\Cal T)
$$
provided by 3.6.

Conversely, to obtain $D^{(n)}F$ from $D_nF$ one simply takes the
$n^{th}$ cross-effect in the case when $\Cal C_Y=\Cal T_Y$, while in the case
when
$\Cal C_Y=\Cal U_Y$ one knows that $\phi^*D^{(n)}F$ determines $D^{(n)}F$ by
4.2 and is given by
$$
\phi^*D^{(n)}F\cong D^{(n)}\phi^*F=cr_nD_n\phi^*F.
$$
In practice we often work with a functor $\Cal U_Y@>{F_Y}>>\Cal D$ that is
the restriction of a functor $\Cal U@>F>>\Cal D$, meaning the composition of
$F$ with the forgetful functor $\Cal U_Y@>>>\Cal U$.

\definition{5.2 Definition}In this case $D^{(n)}F_Y$ is called the $n$-fold
differential of $F$ at $Y$ and denoted by $D^{(n)}_YF$.
\enddefinition

The next theme to be developed is the description
of multilinear functors by their
\lq\lq coefficients\rq\rq\. We discuss this briefly in the important special
case when $Y$ is the one-point space, and then a little more elaborately in the
general case.

If $C$ is a spectrum then the functor
$$
L(X_1,\dots,X_n)=C\wedge(X_1\wedge\dots\wedge X_n)
$$
from ($n$-tuples of) based spaces to spectra is multilinear, and if $C$ has an
action of the symmetric group $\Sigma_n$ then $L$ is symmetric.

Conversely, if $L$ is a
multilinear functor from based spaces to spectra then, taking $C$ to be
$L(S^0,\dots,S^0)$, we have (essentially, see [\Cite{1}; page 5]) a natural
assembly map
$$
C\wedge(X_1\wedge\dots\wedge X_n)@>>>L(X_1,\dots,X_n).
$$
This is an equivalence when $X_j=S^0$ for all $j$, and it follows (see 5.8
below), using the multilinearity of both functors, that it is an equivalence
when the
$X_j$ are finite complexes. If $L$
satisfies a suitable limit axiom (5.10 below), then this even holds 
for all $X_j$.
If $L$ is symmetric then $C$ gets a $\Sigma_n$-action and the assembly map
respects the symmetry.

In short, symmetric multilinear functors from finite based spaces to spectra
correspond precisely, in the sense of an equivalence of homotopy categories, to
spectra with $\Sigma_n$-action. The spectrum (together with its
$\Sigma_n$-action) is called the \it coefficient\rm\ of the symmetric
multilinear functor, or of the corresponding homogeneous functor.

Thus if $F$ is a homotopy functor from based
spaces to based spaces or spectra then for every $n>0$ the $n$-homogeneous
layer
$D_nF$ of its Taylor tower is governed by a certain spectrum with
$\Sigma_n$-action. This will be called the $n^{th}$ derivative of $F$ at the
one-point space and denoted by $\partial^{(n)}F(*)$. In the 
spectrum-valued case
we have
$$
\partial^{(n)}F(*)\simeq(D^{(n)}F)(S^0,\dots,S^0)
$$
and
$$
(D_nF)(X)\simeq(\partial^{(n)}F(*)\wedge X^{\wedge n})_{h\Sigma_n}.
$$
In the space-valued case we have
$$
\Omega^{\infty}\partial^{(n)}F(*)\simeq(D^{(n)}F)(S^0,\dots,S^0)
$$
and
$$
(D_nF)(X)\simeq\Omega^{\infty}((\partial^{(n)}F(*)\wedge X^{\wedge
n})_{h\Sigma_n}).
$$
To fully describe $\partial^{(n)}F(*)$ in the space-valued case,
we examine the proof of 3.7 and see that the
$i^{th}$ space in the spectrum is
$$
\Omega^{i\bar V_n}(D^{(n)}F)(S^i,\dots,S^i)
$$
\remark{5.3 Remark}By the proof of 3.8 we have
$$
\partial^{(n)}F(*)\sim\Omega^{V_n}\partial^{(n)}(F\circ\Sigma)(*)
$$
\endremark

\remark{5.4 Remark}If a subgroup
$G\subset\Sigma_n$ acts on a spectrum $C$ then the functor $(C\wedge
X^{\wedge n})_{hG}$ is homogeneous. When written in standard form it is
$$
(((\Sigma_n)_+\wedge_GC)\wedge X^{\wedge n})_{h\Sigma_n};
$$
the coefficient is the \lq\lq induced spectrum\rq\rq\
$(\Sigma_n)_+\wedge_GC$.
\endremark

We now pass from functors of spaces to functors of
spaces over $Y$ while doing our best to retain the principle that a multilinear
functor is determined by its behavior on a small class of objects.

Here is some notation for naming special objects of
$\Cal T_Y$:  If $Z$ is a based space and $y$ is a point in $Y$ then 
let $Y\vee_y
Z$ be the union of $Y$ and $Z$ with $y$ identified to the basepoint 
of $Z$. This
is to be viewed as a space over $Y$ with all of $Z$ being mapped to $y$. More
generally if several points $y_1,\dots,y_n$ and several based spaces
$Z_1,\dots,Z_n$ are given then
$Y\vee_{y_1}Z_1\vee\dots\vee_{y_n}Z_n$ is the space obtained from $Y$ by
attaching $Z_j$ at $y_j$ for all $j$. Again this is to be viewed as 
an object of
$\Cal T_Y$.

Let $F$ be a homotopy functor from unbased spaces to either based spaces or
spectra.

\definition{5.6 Definition}Let $Y$ be a space and let $y_1,\dots,y_n$ be
points in $Y$. The
$n^{th}$ derivative of $F$ at
$(Y,y_1,\dots,y_n)$, denoted
$$
\partial^{(n)}_{y_1,\dots,y_n}F(Y),
$$
is the coefficient spectrum of the
multilinear functor
$$
(Z_1,\dots,Z_n)\mapsto
(D^{(n)}_YF)(Y\vee_{y_1}Z_1,\dots,Y\vee_{y_n}Z_n).
$$
\enddefinition

Thus in the spectrum-valued case we have
$$
\partial^{(n)}_{y_1,\dots,y_n}F(Y)
\simeq(D^{(n)}_YF)(Y\vee_{y_1}S^0,\dots,Y\vee_{y_n}S^0),
$$
$\partial^{(n)}F$ has a symmetry with respect
to permutations of $(y_1,\dots,y_n)$ and is also functorial in $Y$. To
be more precise, we have here a functor whose domain is the category in which
an object is a \lq\lq space with $n$ basepoints\rq\rq\ $(Y,y_1,\dots,y_n)$ and
a morphism
$(Y,y_1,\dots,y_n)@>>>(W,w_1,\dots,w_n)$ is a pair $(f,\pi)$ with $Y@>f>>W$
being a continuous map and $\pi$ a permutation such that $f(y_{\pi(j)})=w_j$
for all
$j$. The points $y_j$ are not assumed distinct.

It is true in various senses, beginning with 5.9 below, that the
$n^{th}$ derivative of a functor determines the behavior of the $n^{th}$
differential, at least on \it finite\rm\ objects.

\definition{5.7 Definition}An object $Y@>i>>X@>r>>Y$ of $\Cal T_Y$ is finite if
$(X,i(Y))$ is a finite CW pair. An object $X@>r>>Y$ of $\Cal U_Y$ is finite
if the space $X$ is finite CW.
\enddefinition

\proclaim{5.8 Proposition}
Let $L_1@>g>>L_2$ be a map between symmetric
multilinear functors from $\Cal T_Y$ to $\Cal Sp$. In order that
$$
L_1(X_1,\dots,X_n)@>g>>L_2(X_1,\dots,X_n)
$$
should be an equivalence whenever the objects $X_j$ are finite, it is enough if
$g$ is an equivalence in the special case when each $X_j$ is
$Y\vee{y_j}S^0$ for an arbitrary point $y_j\in Y$.
\endproclaim

\demo{Proof}The symmetry is irrelevant here, and it is clear that the case
$n=1$ implies the general case. We therefore give a proof for the case of a map
$L_1@>g>>L_2$ between linear functors from
$\Cal T_Y$ to spectra.

Define $\Cal T_Y@>L>>\Cal Sp$ by letting $L(X)$ be the homotopy fiber of
$L_1(X)@>g>>L_2(X)$. The functor $L$
\lq\lq vanishes\rq\rq\ at
$Y\vee_yS^0$, for every point
$y$ in $Y$, in the
sense that $L(Y\vee_yS^0)\simeq *$. We also know that
$L$ is linear; it preserves weak equivalences, it takes
cocartesian squares to cartesian squares, and it vanishes at
$Y$. We must show that $L$ vanishes at every finite object
$X$.

We name some more objects of $\Cal T_Y$: If $Z$ is a
space and $Z@>f>>Y$ is a map then let $Y+_fZ$ be the disjoint union of Y and Z,
considered as an object of $\Cal T_Y$ in the evident way.

If $X$ is an object of
$\Cal T_Y$ obtained by attaching an $m$-cell to another object $X'$, then there
is a cocartesian diagram
$$
\CD
Y+_{f|_\partial} S^{m-1} @>>> X'\\
@VVV @VVV\\
Y+_f D^m @>>> X
\endCD
$$
for some $D^m@>f>>Y$, and thus a cartesian diagram
$$
\CD
L(Y+_{f|_\partial} S^{m-1}) @>>> L(X')\\
@VVV @VVV\\
L(Y+_f D^m) @>>> L(X).
\endCD
$$
$L(X)$ will be contractible if the other three spectra are. Thus an 
induction on
the number of cells in $X-i(Y)$ will be possible as soon as we have dealt with
the cases $X=Y+_fD^m$ and $X=Y+_fS^m$. The sphere case follows from the disk
case by induction on $m$ (beginning with the case $Y+_fS^{-1}=Y$), and the
disk case is taken care of by a
weak equivalence
$$
Y\vee_yS^0 @>>>Y+_f D^n
$$
obtained by choosing a point $y\in f(D^n)$.
\qed\enddemo

Note that if Y is path-connected then it is only necessary to verify
the hypothesis of 5.8 for one choice of $(y_1,\dots,y_n)$, since a path
$I@>f>>Y$ from
$y$ to $y'$ in
$Y$ yields a diagram of equivalences in $\Cal T_Y$:
$$Y\vee_yS^0@>{\sim}>>Y+_fI@<{\sim}<<Y\vee_{y'}S^0.
$$
There is a variant of 5.8, with the same proof, for unreduced-multilinear
functors from $\Cal U_Y$ to spectra, with a point over $y_j$ replacing
$Y\vee_{y_j}S^0$. The proof of 5.8 also generalizes rather easily to 
prove that a
map between
$d$-excisive homotopy functors from $\Cal T_Y$ [resp.\ $\Cal U_Y$] to
$\Cal Sp$ must be an equivalence for all objects if it is an 
equivalence for all
objects of the form
$Y+_fS$ [resp.\ $S$] where $S$ is a discrete set of at most $d$ points with a
map to $X$.

For us the main consequence of 5.8 is:

\proclaim{Corollary 5.9}If a map $F@>>>G$ of homotopy functors $\Cal
U@>>>\Cal T$ induces an equivalence of $n^{th}$ derivatives
$$
\partial_{y_1,\dots,y_n}F(Y)@>>>\partial_{y_1,\dots,y_n}G(Y)
$$
for every point $(y_1,\dots,y_n)\in Y^n$ then it induces an equivalence of
$n^{th}$ differentials
$$
D_nF_Y(X)@>>>D_nG_Y(X)
$$
for every finite object $X$ of $\Cal U_Y$ or $\Cal T_Y$.
\endproclaim

\demo{Proof}Apply 5.8 to $D^{(n)}F_Y@>>>D^{(n)}G_Y$, noting that the 
behavior of
$D_nF_Y$ on finite objects is determined by the behavior of $D^{(n)}F_Y$ on
finite objects.
\qed\enddemo

These statements definitely require spectra rather than
spaces as the output of the functors. For example, if $C$ is a spectrum and $L$
is the linear functor from based spaces to based spaces given by
$$
L(X)=\Omega^{\infty}(C\wedge X),
$$
then the statement that $L(S^0)$ is contractible means only that the
homotopy groups $\pi_j(C)$ are trivial for $j\geq 0$.

The restriction to finite objects in the results above can be removed if $F$
satisfies a suitable limit axiom.

\definition{5.10 Definition}A homotopy functor $\Cal C@>F>>\Cal D$ is
\it finitary\rm\ if it preserves filtered homotopy colimits up to weak
equivalence, that is, if the natural map
$$
\hocolim_{\alpha}F(X_\alpha)@>>>F(\hocolim_\alpha X_\alpha)
$$
is a weak equivalence whenever $\lbrace X_\alpha\rbrace$ is a diagram in
$\Cal C$ indexed by a filtering category.
\enddefinition

In the case of linear functors from spaces to spaces, this condition
means that the corresponding homology theory satisfies Milnor's wedge axiom.
It is clear that the $n^{th}$ Taylor approximation of a finitary 
functor is itself
finitary.

\remark{5.11 Remark}This was called the \it limit
axiom\rm\ in \cite{2}, but it seems useful to have an adjective available. The
term \lq\lq continuous functor\rq\rq\ has been used, but we would rather
reserve that for something else (a functor that behaves continuously on
morphisms). \endremark

Since every space is equivalent to a filtered $\hocolim$ of finite CW
complexes, a finitary functor of spaces is determined (up to natural weak
equivalence) by its behavior on finite complexes. More generally a finitary
functor of objects in $\Cal C_Y$ is determined by its behavior on objects
that are finite in the sense of 5.7.

According to 5.9  the $n^{th}$ differential (which knows all about the
$n^{th}$ layer of the Taylor tower) is in some sense controlled by the
$n^{th}$ derivative. This is only a weak sense, however, since 5.9 cannot
produce an equivalence between the differentials of $F$ and $G$ unless a map
between $F$ and $G$ is already given. It
would be better to have a way of building the $n^{th}$ differential from
the $n^{th}$ derivative.

The $n^{th}$ derivative of $F$ at $Y$ gives a spectrum for each ordered
$n$-tuple of points in $Y$. To identify the
$n^{th}$ differential, one needs just a little more information. We will sketch
one of several possible answers to the following vague question:

\remark{5.13 Vague question}How can $D^{(n)}F$ be assembled from the spectra\nl
$\partial^{(n)}_{y_1,\dots,y_n}F(Y)$?
\endremark

Answers in the case $n=1$ tend to generalize easily to the general case, so we
will concentrate on that case.

A useful point of view is that a (finitary) linear
functor of spaces over $Y$ corresponds to a \lq\lq coefficient system\rq\rq\
on $Y$ which assigns a spectrum $E_y$ to each point $y$. The spectrum $E_y$
must depend on the point $y$ continuously in some sense, so the linear functor
cannot really be specified by merely giving $E_y$ for each $y$. (That would be
like trying to specify a vector bundle by giving all of its fibers, a mere
collection of vector spaces.)

One can give rigorous sense to
this idea by making the following definition: A \it system of spectra\rm\ on
$Y$ consists of objects
$\lbrace E_n\rbrace$ in $\Cal T_Y$ such that the structural maps $E_n@>>>Y$
are fibrations, related by maps
$$
\Sigma_YE_n @>>>E_{n+1}.
$$
We may write $E_{y,n}$ for the fiber of $E_n@>>>Y$ over $y\in Y$. The
spaces $\lbrace E_{y,n}\rbrace$ for fixed $y$ constitute a spectrum $E_y$, the
\it fiber\rm\ of $E$ over $y$.

Such a system $E$ determines a prespectrum whose
$n^{th}$ space is the homotopy cofiber of the structural map $Y@>>>E_n$. The
associated spectrum, which may be thought of as the homology of $Y$ with
coefficients in
$E$, will be denoted by
$$
\int_{y\in Y} E_y \,dy.
$$
Of course this notation is deceptive, since it appears not to matter how the
various spectra $E_y$ are related.

The integral signs are not meant to suggest antidifferentiation.

We draw attention to the familiar special case when $E$ is a \lq\lq
trivial bundle of spectra\rq\rq\ over
$Y$: Let $C$ be a spectrum and take $E_n$ to be $Y\times C_n$. Then
$$
\int_{y\in Y} E_y \,dy = \int_{y\in Y} C \,dy \cong C\wedge Y.
$$
In the general case there is a spectral sequence of Atiyah-Hirzebruch type,
with
$$
E^2_{p,q}H_p(Y;\pi_q(E_y))\Rightarrow
\pi_{p+q}\int_{y\in Y} E_y \,dy.
$$
One constructs it by taking the direct limit over $n$
of a spectral sequence with
$$
E^2_{p,q}=H_p(Y;\pi^S_{q+n}(E_{n,y}))\Rightarrow
\pi^S_{p+q+n}\hocofiber(Y@>>>E_n).
$$
Systems of spectra pull back: a map
$X@>f>>Y$ and a system $E$ on $Y$ determine a system
$f^*E$ on $X$ given by $(f^*E)_n=E_n\times_Y X$, whose fiber at $x\in 
X$ is isomorphic to
$E_{f(x)}$. Thus a system $E$ on $Y$ gives a functor from
$\Cal U_Y$ to
$\Cal Sp$:
$$
X=(X,X@>f>>Y)\mapsto \int_{x\in X} E_{f(x)} \,dx.
$$
This is a homotopy functor, and it is finitary and $1$-excisive; all of these
assertions can proved by spectral sequence comparison arguments. This 
functor vanishes
at the empty set, so it is not linear  but rather
unreduced-linear (4.5). Reducing it,
one gets a linear functor, which sends
$X$ to the homotopy fiber of
$$
\int_{x\in X} E_{f(x)} \,dx@>>>\int_{y\in Y} E_y \,dy.
$$
It is fairly clear that one could classify the finitary linear functors on
$\Cal T_Y$ along these lines, but we will not pursue that here. The
generalization to $n>1$ involves systems of spectra over
$Y^n$ with a $\Sigma_n$-symmetry, or alternatively systems of spectra over
$Y^n\times_{\Sigma_n}E\Sigma_n$.

If $Y$ is equivalent to the
classifying space of a group $G$ then systems of spectra on $Y$ are the same
(at the homotopy category level) as spectra with $G$-action. This idea can be
extended to simplicial groups $G$, so that it applies to all based connected
$Y$.

The idea of constructing linear functors by systems of spectra over a space was
implicitly present in sections 2 and 3 of \cite{1}. The systems that 
arose there
were mainly \lq\lq fiberwise suspension spectra\rq\rq\ in the following
sense: If
$W$ is an object of
$\Cal T_Y$ whose structural map $W@>>>Y$ is a fibration, then the repeated
fiberwise suspensions
$\Sigma_Y^n W$ form a system
$E$ whose fiber $E_y$ is the suspension spectrum of the fiber of $W$ over $y$.

There is also a dual construction, in which $E$ is used as coefficients for
twisted cohomology rather than twisted homology. If
$E$ is a system of spectra on a CW complex $K$ then the spectrum
$$
\int^{k\in K} E_k \,dk.
$$
is defined by letting the $n^{th}$ space be the space of sections of 
the fibration
$E_n@>>>Y$. If $K$ is locally compact then
there is also the compactly supported version
$$
\int_c^{k\in K} E_k \,dk,
$$
made out of spaces of compactly supported sections.

If $X\mapsto E(X)$ is a homotopy functor from spaces to systems of spectra on
the finite complex
$K$ then \lq\lq differentiation under the integral\rq\rq\ is valid; there is
an equivalence
$$
\partial_y\int^{k\in K}E_k(Y)\, dk
\sim \int^{k\in K}\partial_yE_k(Y)\, dk.
$$
The proof is by induction on the number of cells in $K$. This extends to the
compactly supported case if for example the one-point compactification of $K$
is a finite complex.

We will also need the following simple principle: If $P@>{\pi}>>B$ is a
principal
$G$-bundle ($G$ being a topological group) and
$E$ is a system of spectra on the locally compact space $B$, then when $G$ is
made to act on the spectrum
$$
\int^{p\in P}E_{\pi (p)}\, dp
$$
in the obvious way the homotopy fixed point spectrum is equivalent to
$$
\int^{b\in B}E_b\, db.
$$
\remark{Confession}We should really distinguish between (fiberwise)
unreduced suspension and (fiberwise) reduced suspension. The former,
(fiberwise) join with a two-point set, is what we ordinarily denote by
$\Sigma_Y$ here, and it has the pleasant feature that it takes fibrations to
fibrations. The latter, on the other hand, is much better for making
spectra. When they are inequivalent, we have to choose the former, but this
must be paid for by growing whiskers. We have not systematically imposed
a solution of this small technical difficulty on the reader,  because we do not
have a neat solution.
\endremark

\head \S6. Multilinearized cross-effects\endhead

Before working on some examples, we need one more tool.

\proclaim{6.1 Theorem}Let $\Cal D$ be $\Cal T$ or $\Cal Sp$. The $n^{th}$
differential $D^{(n)}F=cr_nD_nF$ of a homotopy functor
$\Cal T_Y@>F>>\Cal D$ is (naturally, weakly) equivalent to the
multilinearization of the $n^{th}$ cross-effect $cr_nF$ of $F$
itself.\endproclaim

This means that
$(D^{(n)}F)(X_1,\dots,X_n)$ is essentially the homotopy colimit of
$$
\Omega^{i_1+\dots
+i_n}(cr_nF)(\Sigma^{i_1}_YX_1,\dots,\Sigma^{i_n}_YX_n)
$$
over
$(i_1,\dots,i_n)$. The latter may
also be described as
$$
\hocolim_i \Omega^{iV_n}(cr_nF)(\Sigma^i_YX_1,\dots,\Sigma^i_YX_n).
$$
In particular the $n^{th}$ derivative
$\partial^{(n)}_{y_1,\dots,y_n}F(Y)$ is equivalent to the associated 
spectrum of a prespectrum
whose
$i^{th}$ space is
$$
\Omega^{i\bar V_n}(cr_nF)(Y\vee_{y_1}S^i,\dots,Y\vee_{y_n}S^i).
$$
6.1 can be a key tool for identifying $D_nF$ in examples. While the definitions
of
$P_nF$ and $D_nF$ are difficult to use for explicit calculation when $n$ is
greater than one, we do know how to recover
$D_nF$ from $cr_nD_nF$, and 6.1 says that this in turn can be obtained
rather directly from $F$ itself.

\remark{6.2 Remark}Like 1.7  and 2.1 , 6.1 has an easier proof in the case when
the functor
$F$ is analytic. In fact, in that case using 1.5(2) one sees easily that if the
maps
$X_j@>>>Y$ are all $k$-connected then the map
$$
(cr_nF)(X_1,\dots,X_n)@>>>(cr_nP_nF)(X_1,\dots,X_n)
$$
induced by $p_nF$ is
$((n+1)k-c')$-connected for some constant $c'$, from which it follows that the
canonical map from
$$
\Omega^{iV_n}(cr_nF)(\Sigma^iX_1,\dots,\Sigma^iX_n)
$$
to
$$
\Omega^{iV_n}(cr_nP_nF)(\Sigma^iX_1,\dots,\Sigma^iX_n)\sim
(cr_nP_nF)(X_1,\dots,X_n)
$$
has a connectivity tending to $\infty$ with $i$.
For that matter, in applying 6.1 to a particular $F$ one often uses
this same kind of reasoning again: One identifies the multilinearization
of
$cr_nF$ with a given functor $L$ by exhibiting a (natural and
symmetry-preserving) map
$$
(cr_nF)(X_1,\dots,X_n)@>>>L(X_1,\dots,X_n)
$$
and checking that it is
$((n+1)k-c')$-connected for some constant $c'$ when the
$X_j@>>>Y$ are $k$-connected. Therefore some readers may prefer to skip to \S7
after the proof of 6.3.
\endremark

A consequence of 6.1 is that what we are calling a second
derivative (that is, the coefficient of the bilinear functor corresponding to a
$2$-homogeneous layer) can actually be seen as the derivative of a derivative.
More generally, we have the following useful interpretation of 6.1:

\proclaim{6.3 Corollary}
$$
\partial^{(p+q)}_{y_1,\dots, y_{p+q}}F(Y)\sim\partial^{(p)}_{y_1,\dots,
y_p}\partial^{(q)}_{y_{p+1},\dots, y_{p+q}}F(Y).
$$
\endproclaim
\demo{Explanation and sketch of proof}
If $F$ is a homotopy functor from spaces to either based spaces or spectra,
then, as we have already observed, the spectrum
$\partial^{(q)}_{y_{p+1,\dots, y_p+q}}F(Y)$ depends functorially on
$Y=(Y,y_{p+1}\dots,y_{p+q})$ with appropriate definitions. In stating 6.3 we
are extending the \lq\lq partial derivative\rq\rq\ notation
to cover functors of \lq\lq spaces with several base points\rq\rq\. Thus,
for example, if $F$ is a functor of based spaces then $\partial_yF(Y,y_0)$
depends functorially on $(Y,y,y_0)$ and is defined as the coefficient spectrum
of the linear part of the functor
$$
Z\mapsto F(Y\vee_yZ).
$$
Here $y_0\in Y$ is serving as basepoint in $Y\vee_yZ$ for
purposes of applying the functor $F$ while $y$ is serving
as basepoint in $Y$ for wedging with $Z$.

In view of 6.1, to prove 6.3 one
has only to see that for a functor of $p+q$ based spaces, such as
$$
(Z_1,\dots,Z_{p+q})\mapsto F(Y\vee_{y_1} Z_1,\dots,Y\vee_{y_n}
Z_{p+q}),
$$
the following two processes are equivalent: (1) reducing in all variables
followed by multilinearizing in all variables, (2) reducing and 
multilinearizing
in the last $q$ variables, followed by reducing and multilinearizing
in the first $p$ variables. This is easy.
\qed\enddemo

Another consequence of 6.1 is:

\proclaim{6.4 Corollary}$F$ is $m$-reduced if and only if $F$ is reduced and
for every
$0<n<m$  the multilinearization of the $n^{th}$ cross-effect of $F$ is
contractible.
\endproclaim
\demo{Proof}
This is clear from 6.1 and the fact that $cr_nD_nF$ determines $D_nF$.
\qed\enddemo

  The proof of 6.1 is connected with the idea of \lq\lq multivariable
Taylor series\rq\rq\ (Remark 1.22).
A key observation is that the product category $\Cal
C_{Y_1}\times\dots\times
\Cal C_{Y_n}$ is itself the category of spaces over a space: it is 
equivalent to
$\Cal C_{Y_1
{\coprod}\dots{\coprod} Y_n}$ where the \lq\lq $\coprod$\rq\rq\ denotes
disjoint union. In this equivalence of categories,
$n$-tuples of weak equivalences correspond to weak equivalences and
$n$-tuples of cocartesian cubes correspond to
cocartesian cubes. When it is necessary to distinguish between a functor
$$
\Cal C_{Y_1}\times\dots\times
\Cal C_{Y_n}@>F>>\Cal D
$$
of $n$ variables and the associated functor
$$
\Cal C_{Y_1 {\coprod}\dots{\coprod} Y_n}@>>>\Cal D
$$
of one variable, we will denote the latter by $\lambda F$.

6.1 is proved using the following statement, which will be
proved at the end of \S 6:

\proclaim{6.5 Lemma}If
$\Cal C_{Y_1}\times\dots\times \Cal C_{Y_n}@>G>>\Cal D$
  is $(1,\dots,1)$-reduced then
$\lambda P_{1,\dots,1}G\sim P_n\lambda G$.\endproclaim

\demo{Proof of 6.1}Applying 6.5 in the special case
$Y=Y_1=\dots =Y_n$, with $G=cr_nF$, we find that what we need is a natural
equivalence
$$
P_n\lambda cr_nF \sim \lambda cr_nP_nF.
$$
This will follow from a natural equivalence $
T_n\lambda cr_nF \sim \lambda cr_nT_nF
$, which in turn will follow from a natural equivalence $
J_U\lambda cr_nF \sim \lambda cr_nJ_UF
$, where (as in the proof of 1.8) $J_U$ is composition with the fiberwise join
with a finite set
$U$. $\lambda\circ cr_n$ is a three-step process: Compose with the
(whiskered) sum
$$
\Cal T_Y^n@>{+}>>\Cal T_Y,
$$
then reduce in all variables (this was called $r$ in the proof of 4.2), then
compose with the equivalence of categories
$$
\Cal T_{Y\coprod\dots\coprod Y}@>>>\Cal T_Y^n.
$$
Each of these steps commutes with $J_U$ up to natural equivalence.
\qed\enddemo

To get to 6.5 we must revisit and strengthen some results that were discussed
in \S 3.  The following is a strengthening of [\Cite{2}; 3.4].

\proclaim{6.6 Lemma}If a homotopy functor $\Cal
C_{Y_1}\times\dots\times
\Cal C_{Y_n}@>F>>\Cal D$ is $(d_1,\dots,d_n)$-excisive then $\lambda F$
is $(d_1+\dots +d_n)$-excisive.\endproclaim

\demo{Proof}In fact the proof of [\Cite{2}; 3.4] becomes a proof of
6.6 if one uses $(X_1,\dots,$ $X_n)$ in place of $(X,\dots,X)$
throughout.\qed\enddemo

In particular $\lambda F$ is $n$-excisive if $F$ is $(1,\dots,1)$-excisive. We
need to know also that
$\lambda F$ is
$n$-homogeneous if $F$ is multilinear. The following is a strengthening of 3.1.

\proclaim{6.7 Lemma}If $\Cal
C_{Y_1}\times\dots\times
\Cal C_{Y_n}@>F>>\Cal D$ is $(1,\dots,1)$-reduced, then
$\lambda F$ is $n$-reduced.\endproclaim

\demo{Proof}This follows from the next statement as 3.1 followed from
3.2.\qed\enddemo

\proclaim{6.8 Lemma}If $\Cal
C_{Y_1}\times\dots\times
\Cal C_{Y_n}@>F>>\Cal D$ is $(1,\dots,1)$-reduced, then the map
$$
\lambda F @>{t_{n-1}\lambda F}>> T_{n-1}\lambda F
$$
factors through a weakly contractible functor.\endproclaim

\demo{Proof}Again, the proof of 3.2 applies with no change except
$(X_1,\dots,X_n)$ for $(X,\dots,X)$.\qed\enddemo

We need this partial converse to 6.6:

\proclaim{6.9 Lemma}If $\Cal
C_{Y_1}\times\dots\times
\Cal C_{Y_n}@>L>>\Cal D$ is $(1,\dots,1)$-reduced and $\lambda
L$ is $n$-excisive, then $L$ is multilinear.\endproclaim

\demo{Proof}We have to show that $L$ is $1$-excisive in each variable. It is
sufficient to consider the last variable. Fix an object $X_j$ of $\Cal C_{Y_j}$
for each
$1\leq j \leq {n-1}$ and let
$$
\CD
A @>>> B\\
@VVV @VVV\\
C@>>> D
\endCD
$$
be any cocartesian square in $\Cal C_{Y_n}$. We have to show that
the square
$$
\CD
L(X_1,\dots,X_{n-1},A) @>>> L(X_1,\dots,X_{n-1},B)\\
@VVV @VVV\\
L(X_1,\dots,X_{n-1},C)@>>> L(X_1,\dots,X_{n-1},D)
\endCD
$$
is cartesian.

Define an $(n+1)$-cube in $\Cal
C_{Y_1}\times\dots\times
\Cal C_{Y_n}$ as follows: For each subset $S$ of
$\lbrace 1,\dots,{n-1}\rbrace$, let $X_j(S)$ be $Y_j$ if $j\in S$ and $X_j$ if
$j\notin S$. 

Then $\lbrace (X_1(S),\dots,X_{n-1}(S))\rbrace$ constitutes an
$(n-1)$-cube in $\Cal
C_{Y_1}\times\dots\times
\Cal C_{Y_{n-1}}$ and our $(n+1)$-cube will consist of the
squares
$$
\CD
(X_1(S),\dots,X_{n-1}(S),A) @>>> (X_1(S),\dots,X_{n-1}(S),B)\\
@VVV @VVV\\
(X_1(S),\dots,X_{n-1}(S),C)@>>> (X_1(S),\dots,X_{n-1}(S),D).
\endCD
$$
The $(n+1)$-cube is strongly cocartesian, so $L$ yields a cartesian cube. On
the other hand, for each nonempty
$S$ the square
$$
\CD
L(X_1(S),\dots,X_{n-1}(S),A) @>>> L(X_1(S),\dots,X_{n-1}(S),B)\\
@VVV @VVV\\
L(X_1(S),\dots,X_{n-1}(S),C)@>>> L(X_1(S),\dots,X_{n-1}(S),D)
\endCD
$$
  is cartesian; in fact it is made up of contractible objects because $L$ is
reduced in each variable. It follows by 1.6 of  \cite{2} that
the
square corresponding to $S=\emptyset$ is also cartesian.
\qed\enddemo

\demo{Proof of 6.5}Now it is convenient to drop the
distinction between $F$ and $\lambda F$. $P_{1,\dots,1}F$ is the universal
example of a
$(1,\dots,1)$-excisive functor under $F$, and is
also $n$-excisive (by 6.6).
$P_nF$ is the universal example of an $n$-excisive functor under $F$,
and is also $(1,\dots,1)$-excisive (by 6.9). It follows that they
are the same. \qed\enddemo

The reader, looking at 6.6 and 6.7, might have wondered about:

\proclaim{6.10 Lemma}If a homotopy functor $\Cal
C_{Y_1}\times\dots\times
\Cal C_{Y_n}@>F>>\Cal D$ is $(d_1,\dots,d_n)$-reduced then $\lambda F$
is $(d_1+\dots +d_n)$-reduced.\endproclaim

In fact this is true, and it can be deduced from 6.4.

\head\S7. Example: Suspension spectra of mapping spaces\endhead

For an unbased space $X$ let $\Sigma^{\infty}_+X$ be the suspension
spectrum of the based space $X_+$ obtained by adding an extra point to $X$.
We will call this the unreduced suspension spectrum (and hope that 
this does not
lead anyone to confuse the unreduced suspension of
$X$ with
$S^1\wedge X_+$). Similarly, if $X$ is a space
fibered over $Y$ we can speak of its unreduced fiberwise suspension spectrum,
meaning the fiberwise suspension spectrum of the
fiberwise based space $\psi X$, where as in \S4  this means
the disjoint union of $X$ and $Y$ considered as an object of $\Cal T_Y$.

For a
finite CW complex
$K$, the functor
$F(X)=\Sigma^{\infty}_+X^K$ is analytic by [\Cite{2}; 4.4], and its 
first derivative
was found in [\Cite{1}; 2.4]. We now find its
$n^{th}$ derivative. We begin by recalling what the formula for the first
derivative is and where that formula came from.

In the notation of
\S5 the formula for the first derivative is
$$
\partial_y\Sigma^{\infty}_+Y^K\sim \int^{k\in
K}\Sigma^{\infty}_+(Y,y)^{(K,k)}\,dk.
\tag 7.1
$$
Recall that this means that the linearization of the functor
$$
\align
Z&\mapsto \hofiber(\Sigma^{\infty}_+(Y\vee_y
Z)^K@>>> \Sigma^{\infty}_+Y^K)\\
&\sim \Sigma^{\infty}\hocofiber(Y^K@>>>(Y\vee_y Z)^K)
\tag 7.2\endalign
$$
is the functor
$$
Z\mapsto Z\wedge \int^{k\in
K}\Sigma^{\infty}_+(Y,y)^{(K,k)}\,dk.
$$
Because $K$ is finite, this last can also be written
$$
\int^{k\in
K}\Sigma^{\infty}(Z\wedge {{(Y,y)^{(K,k)}}_+})\,dk.
\tag 7.3
$$
Implicitly in  7.1 we are using a certain
system of spectra on $K$, namely the unreduced fiberwise suspension spectrum
of a certain space over
$K$, let us call it $W$, whose fiber over $k\in K$ is $W_k=(Y,y)^{(K,k)}$; $W$
is the subspace of $Y^K\times K$ consisting of
pairs
$(f,k)$ such that $f(k)=y$. Likewise in 7.3 we are using the fiberwise
suspension spectrum of a space over
$K$, call it
$Z\wedge_K \psi W$, whose fiber is
$Z\wedge {W_k}_+$, namely
$$
\colim(Z\times W @<<< W @>>>K).
$$

\remark{7.4 Remark}A formula for the differential
and not just the derivative was given in \cite{1}. In the present 
notation it says
that the unreduced-linear functor corresponding to $D_YF$ takes the object
$X@>f>>Y$ to
$$
\int_{x\in X}\int^{k\in K}
\Sigma^{\infty}_+(Y,f(x))^{(K,k)}\,dk\,dx.
$$
In other words, the linear functor is given by a coefficient system
on $Y$ which may be obtained by \lq\lq integration over the fiber\rq\rq\ from
a system on $K\times Y$, the fiberwise unreduced suspension spectrum of
the fibration $K\times Y^K@>>>K\times Y$ whose fiber over $(k,y)$ may be
identified with
$(Y,y)^{(K,k)}$.
\endremark

The method of
proof for 7.1 was this: First give a natural map from 7.2 to 7.3, 
then show that
its connectivity is roughly twice that of the space
$Z$.

To produce the map, it was enough to give a natural map of based spaces
$$
\hocofiber(Y^K@>>>(Y\vee_y Z)^K)@>>>
\Omega^{\infty}\int^{k\in
K}\Sigma^{\infty}(Z\wedge {(Y,y)^{(K,k)}}_+)\,dk.
$$
This was done by means of a tautological map from $(Y\vee_y Z)^K$  to the
space of sections of $Z\wedge_K
\psi W@>>>K$. To specify this map we say where it sends the map
$K@>{\tilde f}>>Y\vee_{y}Z$. Let $f$ be the composed map
$$
K@>{\tilde f}>>Y\vee_{y}Z@>>>Y
$$
Then $\tilde f$ is sent to the section whose value at
$k$ is $\tilde f(k)\wedge f\in Z\wedge {(Y,y)^{(K,k)}}_+$ if $\tilde 
f(k)\in Z$ and
otherwise is the (fiberwise) basepoint.

It is a good precaution to add a \lq\lq whisker\rq\rq\ to $Z$, replacing
$Y\vee_y Z$ by
$Y\cup_y I\cup_z Z$, before making the construction just described, to insure
that the map
$Z\wedge_Y \psi W@>>>K$ is a fibration.

The proof that the resulting map from 7.2 to 7.3 is highly connected 
will not be
repeated here.

There is a variant of 7.1 for spaces of based maps: If $K$ has
a basepoint
$k_0$, then we obtain
$$
\partial_y\Sigma^{\infty}_+{(Y,y_0)^{(K,k_0)}}\sim \int_c^{k\in
K-{\lbrace k_0\rbrace}}\Sigma^{\infty}_+{(Y,y,y_0)^{(K,k,k_0)}}\,dk.
$$
(This time the whisker is even more important, because it insures that
even if $y=y_0$ the section being constructed will have compact
support.)

Now we are in a position to compute a second derivative, using 6.2. We have:
$$
\align
\partial_{y_1}\partial_{y_2}\Sigma^{\infty}_+{Y^K}
&\sim
\partial_{y_1}\int^{{k_2}\in K}\Sigma^{\infty}_+{(Y,y_2)^{(K,k_2)}}\,dk_2\\
&\sim
\int^{{k_2}\in K}\partial_{y_1}\Sigma^{\infty}_+{(Y,y_2)^{(K,k_2)}}\,dk_2\\
&\sim
\int^{{k_2}\in K}\int_c^{k_1\in
K-{\lbrace k_2\rbrace}}\Sigma^{\infty}_+{(Y,y_1,y_2)^{(K,k_1,k_2)}}\
dk_1\,dk_2\\
&\sim \int_c^{(k_1,k_2)\in
K^{(2)}}\Sigma^{\infty}_+{(Y,y_1,y_2)^{(K,k_1,k_2)}}\
d\,(k_1,k_2).
\tag 7.5 \endalign
$$
The differentiation under the integral sign depends on the hypothesis that
$K$ is finite.
$K^{(2)}$ is the complement of the diagonal in $K\times K$. The
last step, replacing a double integral by a single integral, is a
tautology.

The
expression in the last line implicitly refers to a system of spectra on
$K^{(2)}$ whose fiber at $(k_1,k_2)$ is
$\Sigma^{\infty}_+{(Y,y_1,y_2)^{(K,k_1,k_2)}}$, namely the unreduced fiberwise
suspension spectrum of a certain fibration
$$
W^{[2]}@>>>K^{(2)}.
$$
The space $W^{[2]}$ is the subspace of $Y^K\times K^{(2)}$
consisting of all $(f,k_1,k_2)$ such that $f(k_1)=y_1$ and
$f(k_2)=y_2$, so that the fiber $W^{[2]}_{(k_1,k_2)}$ over $(k_1,k_2)$ is
$(Y,y_1,y_2)^{(K,k_1,k_2)}$.

This is not enough to determine the quadratic functor $D_2F_Y$,
because we do not yet know the $\Sigma_2$-symmetry in the second
derivative. We can make a good guess about that, because the last
expression in 7.5 has an
obvious symmetry. To verify the guess, we can proceed as follows:

We have to study the second-order cross-effect of $F_Y$ as applied to objects
$Y\vee_y Z$, in other words the total homotopy fiber of
$$
\CD
\Sigma^{\infty}_+({Y\vee_{y_1}Z_1\vee_{y_2}Z_2})^K
@>>>\Sigma^{\infty}_+({Y\vee_{y_1}Z_1})^K\\
@VVV @VVV\\
\Sigma^{\infty}_+({Y\vee_{y_2}Z_2})^K@>>>\Sigma^{\infty}_+({Y})^K.
\endCD
\tag 7.6
$$
We know that the bilinearization of this is equivalent to
$$
Z_1\wedge Z_2\wedge\int_c^{(k_1,k_2)\in
K^{(2)}}\Sigma^{\infty}_+{(Y,y_1,y_2)^{(K,k_1,k_2)}}\
d(k_1,k_2),
$$
or equivalently
$$
\int_c^{(k_1,k_2)\in
K^2-\Delta}\Sigma^{\infty}(Z_1\wedge
Z_2\wedge{(Y,y_1,y_2)^{(K,k_1,k_2)}}_+)\ d(k_1,k_2).
\tag 7.7
$$
This last expression refers to the fiberwise suspension spectrum of a
certain space over
$K^{(2)}$ whose fibers are
$Z_1\wedge Z_2\wedge {W^{[2]}_{(k_1,k_2)}}_+.$ Call it $(Z_1\wedge
Z_2)\wedge_{K^{(2)}} \psi W^{[2]}$.

We are seeking to show that that equivalence can be chosen to preserve the
$\Sigma_2$-symmetry, so we should look for a
symmetry-preserving map from the total homotopy fiber of 7.6 to 7.7.

The
total homotopy fiber of 7.6 can be rewritten as the suspension spectrum of the
total cofiber of
$$
\CD
(Y)^K @>>>(Y\vee_{y_1}Z_1)^K\\
@VVV @VVV\\
(Y\vee_{y_2}Z_2)^K @>>>(Y\vee_{y_1}Z_1\vee_{y_2}Z_2)^K.
\endCD
\tag 7.8
$$
There is a tautological map from this total cofiber to the zeroth space of
7.7.
It arises from a tautological map
from $(Y\vee_{y_1}Z_1\vee_{y_2}Z_2)^K$ to the space of compactly supported
sections of $(Z_1\wedge Z_2)\wedge_Y
W^{[2]}@>>>K^{(2)}$. To specify this we say where it sends the map
$K@>{\tilde f}>>Y\vee_{y_1}Z_1\vee_{y_2}Z_2$. Let $f$ be the composed map
$$
K@>{\tilde f}>>Y\vee_{y_1}Z_1\vee_{y_2}Z_2@>>>Y
$$
Then $\tilde f$ is sent to the section whose value at
$(k_1,k_2)$ is $\tilde f(k_1)\wedge\tilde f(k_2)\wedge (f,k_1,k_2)$ if
$\tilde f(k_1)\in Z_1$ and
$\tilde f(k_2)\in Z_2$ and otherwise the (fiberwise) basepoint.

We claim, leaving the remaining details to the reader, that this results in a
map
$$
\partial^{(2)}_{y_1,y_2}\Sigma^{\infty}({Y^K}_+)
@>>>
\int_c^{(k_1,k_2)\in
K^{(2)}}\Sigma^{\infty}_+{(Y,y_1,y_2)^{(K,k_1,k_2)}}\
d(k_1,k_2)
$$
that corresponds to 7.5 under 6.2 and is therefore an equivalence.

The same method gives the $n^{th}$ derivative. The
conclusion is:

\proclaim{7.10 Theorem}For a finite complex $K$ we have a symmetry-preserving
equivalence
$$\partial^{(n)}_{y_1,\dots y_n}\Sigma^{\infty}_+{Y^K}
\sim \int_c^{k\in
K^{(n)}}\Sigma^{\infty}_+{(Y,y_1,\dots,y_n)^{(K,k_1,\dots,k_n)}}\,
dk.
$$
\endproclaim
Here $K^{(n)}$ is the space of all ordered $n$-tuples $k=(k_1,\dots,k_n)$ of
distinct points in $K$.

It is interesting to work out what this says in the case when $K$ is 
a finite set
of cardinality $m$, and to compare it with the formula
$$
m(m-1)\dots (m-n-1)y^{m-n}.
$$
for the $n^{th}$ derivative of $y^m$ in ordinary calculus.

In the case when $Y$ is a
one-point space, the right-hand side of 7.10 becomes
$$
\int_c^{k\in
K^{(n)}}\Sigma^{\infty}S^0\,dk\sim  Map_*(K^{(n)c},\Sigma^{\infty}S^0)
=(K^{(n)c})^*
$$
In other words, the $n^{th}$ derivative of $F(X)=\Sigma^{\infty}_+{X^K}$ at a
point is
the $S$-dual $(K^{(n)c})^*$ of a certain based $\Sigma_n$-space $K^{(n)c}$, the
one-point compactification of
$K^{(n)}$ (or the quotient of $K^n$ by the fat diagonal $K^n-K^{(n)}$). The
$n^{th}$ homogeneous functor is
$$
(D_nF)(X)\sim ((K^{(n)c})^*\wedge
X^{\wedge n})_{h\Sigma_n},
$$
and this can be identified with
$$
  Map_*(K^{(n)c},\Sigma^{\infty}X^{\wedge n})_{h\Sigma_n}
$$
because $K^{(n)c}$ is finite. It can also be identified with
$$
  Map_*(K^{(n)c},\Sigma^{\infty}X^{\wedge n})^{\Sigma_n}
$$
because the group action on $K^{(n)c}$ is free.

The analogous conclusion for based $K$ says
$$\partial^{(n)}_{y_1,\dots y_n}\Sigma^{\infty}_+{(Y,y_0)^{(K,k_0)}}
\sim \int_c^{k\in
{(K-{\lbrace
k_0\rbrace})}^{(n)}}\Sigma^{\infty}_+{(Y,y_0,y_1,\dots,y_n)^{(K,k_0,k_1,\dots,k_n)}}\,
dk.
$$
When $Y$ is a point we find that the $n^{th}$ coefficient of
$F(X)=\Sigma^{\infty}_+{(X,x_0)}^{(K,k_0)}$ is the
$S$-dual of the one-point compactification of ${(K-{\lbrace
k_0\rbrace})}^{(n)}$.

The case when $K$ is a based circle is particularly simple. Since
${(S^1-{\lbrace k_0\rbrace})}^{(n)}$ is the disjoint union of open $n$-cells
freely and transitively permuted by the group, the formula for the $n^{th}$
homogeneous layer of
$\Sigma^{\infty}_+\Omega X$ becomes
$\Omega^n\Sigma^{\infty}(X^{\wedge n})$, as already pointed out in 1.20.

\head\S8. Example: The identity functor\endhead

The identity functor $\Cal T@>I>>\Cal T$ from based
spaces to based spaces is a central example, and its $n^{th}$ derivative
$\partial^{(n)}I(*)$ is a basic object in homotopy theory. This spectrum with
$\Sigma_n$-action turns out to be
$S$-dual to a certain finite complex with $\Sigma_n$-action. We summarize
and discuss some known results.

Note that the problem of determining the $n^{th}$ derivative of the 
identity is equivalent to that of
determining the $n^{th}$ derivative of the functor
$\Omega\Sigma$. In fact, by 5.3 we have
$$
\partial^{(n)}\Omega\Sigma(*)\simeq S^{\bar V_n}\wedge\partial^{(n)}I(*)
$$
and therefore
$$
\partial^{(n)}I(*)\simeq\Omega^{\bar V_n} \partial^{(n)}\Omega\Sigma(*)
$$
To begin with the obvious, the first derivative of $I$ is the sphere spectrum.

The second derivative can be identified by using
the second cross-effect. The total homotopy
fiber of
$$
\CD
X_1\vee X_2 @>>> X_1\\
@VVV @VVV \\
X_2 @>>> *
\endCD
$$
is the homotopy fiber of $X_1\vee X_2 @>>>X_1\times X_2$, and the
bilinearization of this is $\Omega Q(X_1 \wedge X_2)$, with the obvious
$\Sigma_2$-symmetry (trivial action on the loop coordinate), simply because
there is a map (natural and symmetry-preserving)
$$
\hofiber (X_1\vee X_2 @>>>X_1\times X_2)
@>>>\Omega Q(X_1 \wedge X_2)
$$
that is approximately $3k$-connected when $X_1$ and $X_2$ are $k$-connected.
It follows that $\partial^{(2)}I(*)$ is the $-1$-sphere with trivial action.

For the $n^{th}$ derivative, partial information
can be obtained from the Hilton-Milnor theorem \cite{9}.
Recall that this describes the space
$\Omega\Sigma(X\vee Y)$ as a weak product (direct limit of finite
products) of factors each of which has the form $\Omega\Sigma
(X^{\wedge a}\wedge Y^{\wedge b})$. The factors for a given pair
$(a,b)$ correspond to certain nested commutator expressions, an integral basis
for the bidegree $(a,b)$ summand of a free Lie ring on two generators
whose bidegrees are $(1,0)$ and $(0,1)$. Iteration yields a description of
$\Omega\Sigma(X_1\vee\dots\vee X_n)$ as a weak product of factors of the
form $\Omega\Sigma(X_1^{\wedge a_1}\wedge\dots\wedge X_n^{\wedge
a_n})$. The cross-effect $(cr_n I)(X_1,\dots,X_n)$ is the product of those
factors for which $a_j\geq 1$ for all $j$. The multilinearized
cross-effect sees only those factors for which $a_j=1$  for all $j$.
Thus the
$n^{th}$ differential of the functor $\Omega\Sigma$ is a product of copies of
$Q(X_1\wedge\dots\wedge X_n)$, and the $n^{th}$ derivative is the product of a
corresponding number of copies of the sphere spectrum. The number of copies is
$(n-1)!$, and they correspond to a basis for the group $Lie(n)$ generated
by all Lie monomials in
$n$ variables such that each variable occurs just once (and of course
considered modulo the Jacobi identity and antisymmetry).
A standard choice of basis for $Lie(n)$ consists of
the monomials
$$
\lbrack x_{\pi (1)},\lbrack x_{\pi (2)},\lbrack\dots \lbrack x_{\pi 
(n-1)} , x_n\rbrack\dots\rbrack\rbrack\rbrack
$$
one for each permutation $\pi$ belonging to the subgroup
$\Sigma_{n-1}\subset\Sigma_n$. It follows
that the $n^{th}$ derivative of
$\Omega\Sigma$, regarded as a spectrum with
$\Sigma_{n-1}$-action, is the suspension spectrum of the finite set
$(\Sigma_{n-1})_+$.

Of course for determining $D_nI$ one needs the full $\Sigma_n$-action. The
method outlined above even identifies the action of $\Sigma_n$ on the 
homology of
the spectrum $\partial ^{(n)}I(*)$ (a free abelian group of rank
$(n-1)!$ concentrated in degree $1-n$); it is the obvious action of 
$\Sigma_n$ on
$Lie(n)$, twisted by signs of permutations. But this is  still insufficient for
determining the homogeneous functor.

Johnson \cite{5} defined a based finite complex $K_n$ with $\Sigma_n$-action
whose
$S$-dual is $\partial^{(n)}I(*)$. Her $K_n$ was designed to admit an 
interesting
map
$$
(cr_nI)(X_1,\dots,X_n)@>>>Map_*(K_n,X_1\wedge \dots\wedge X_n),
$$
both natural and symmetry-preserving, and she showed that after
multilinearization this map leads to an equivalence
$$
(D^{(n)}I)(X_1,\dots,X_n)@>>>Map_*(K_n,Q(X_1\wedge \dots\wedge X_n)).
$$
Arone and Mahowald \cite{6} later came up with another answer to the same
question. It is defined in terms of the poset of all partitions
(equivalence relations) on the set $\underline n=\lbrace 
1,\dots\,n\rbrace$. The
poset has a maximal element (the trivial partition) and a minimal element (the
improper partition), so that the nerve of the poset is (for two reasons)
contractible. The Arone-Mahowald version of $K_n$ can be described as the
double suspension of the nerve of the poset of proper nontrivial partitions, or
alternatively as the  quotient of the nerve of all partitions by the 
union of the
nerve of the nontrivial partitions and the nerve of the proper partitions.

In \cite{6} this model is justified by proving directly that it is 
(equivariantly)
homotopy equivalent to Johnson's space, but it was actually 
discovered from a very
different point of view, which is worked out in detail in \cite{14}. There is a
cosimplicial functor from spaces to spaces which has in degree $d$ the functor
$Q^{d+1}$, iterated composition of $Q$ with itself. When applied to
$1$-connected spaces it serves as a resolution of the identity functor. (In
general it gives the Bousfield-Kan  integral completion functor.) 
For each $d$ the
functor
$Q^{d+1}$ has a split Taylor tower, which can be read off from the Snaith
splitting formula. It is rather easy to see that the
$n^{th}$ coefficient of the $d^{th}$ functor is (functorially in $d$, up to
homotopy) the
$S$-dual of the (discrete) space of $d$-simplices in $K_n$. (The details worked
out in \cite{14} dispose of that unfortunate \lq\lq up to homotopy\rq\rq\.)

Arone and Mahowald use this description of $K_n$ to investigate the mod $p$
homology of the spectrum $(\partial^{(n)}I(*)\wedge X^{\wedge n})_{h\Sigma_n}$
whose zeroth space is $(D_nI)(X)$. Their main results concern the 
case when $X$ is
a sphere. For simplicity take it to be an odd sphere. When
$n$ is not a power of $p$ they find that the layer
$(D_nI)(S^{2m-1})$ is
$p$-locally trivial. (This is equivalent to the
statement that the homology of $\Sigma_n$ with coefficients in $Lie(n)$
localized at $p$ is trivial.) When
$n=p^k$ they use Dyer-Lashof operations to calculate the homology as a module
over the Steenrod algebra, finding in particular that the layer
$(D_{p^k}I)(S^{2m-1})$ has trivial $v_j$-periodic homotopy when $j<k$.

For further insight into these matters, see \cite{13}.

We close with some remarks about the derivatives of the
identity at a general space $Y$.  Recall from the discussion 
following 6.3 that the
$n^{th}$ derivative of
$\Cal T@>I>>\Cal T$ will be a functor of spaces with $n+1$ basepoints,
symmetric with respect to permutations of the last $n$ points. It is clear that
$$
\partial_y (Y,y_0)\sim \Sigma^{\infty}_+P_{y_0}^yY.
$$
where $P_{y_0}^yY$ is the path space $(Y,y_0,y)^{(I,0,1)}$. In general, the
$n^{th}$ derivative of the identity at an arbitrary space may be described in
terms of
the
$n^{th}$ derivative of the identity at a point. For example, given a
space $Y$ and points
$y_0,y_1,y_2$, the total homotopy fiber (with respect to $y_0$) of
$$
\CD
Y\vee_{y_1}Z_1\vee_{y_2}Z_2@>>>Y\vee_{y_1}Z_1\\
@VVV @VVV\\
Y\vee_{y_2}Z_2@>>>Y
\endCD
$$
is equivalent to the total homotopy fiber of
$$
\CD
(P_{y_0}^{y_1}Y_+\wedge Z_1)\vee
(P_{y_0}^{y_2}Y_+\wedge Z_2)
@>>>P_{y_0}^{y_1}Y_+\wedge Z_1\\ @VVV @VVV\\
P_{y_0}^{y_2}Y_+\wedge Z_2@>>>*.
\endCD
$$
(The second square consists essentially of the homotopy fibers over $y_0\in Y$
of the spaces in the first square.) Bilinearizing with respect to
$Z_1$ and
$Z_2$, we find that
$$
(D^{(2)}_{(Y,y_0)}I)(Y\vee_{y_1} Z_1,Y\vee_{y_2} Z_2)\sim
(D^{(2)}_*I)(P_{y_0}^{y_1}Y_+\wedge Z_1,P_{y_0}^{y_2}Y_+\wedge Z_2).
$$
The
same argument applies for any $n$ and yields a natural and symmetrical
equivalence
$$
\partial^{(n)}_{y_1,\dots,y_n}(Y,y_0)\sim
P_{y_0}^{y_1}Y_+\wedge\dots\wedge P_{y_0}^{y_n}Y_+\wedge K_n^*
$$
These observations can be used to give an alternative to the Hilton-Milnor
argument above; if one is willing to settle for 
$\Sigma_{n-1}$-symmetry rather than
$\Sigma_n$-symmetry, then by 6.3 one can write
$$
\partial^{(n)}_{y_1,\dots,y_n}(Y,y_0)\sim\partial^{(n-1)}_{y_1,\dots,y_{n-1}}
\Sigma^{\infty}_+P_{y_0}^{y_n}Y,
$$
which by another very slight generalization of 7.10 is equivalent to
$$
\int_c^{k\in
{(I-{\lbrace
0,1\rbrace})}^{(n-1)}}\Sigma^{\infty}_+{(Y,y_0,y_1,\dots,y_{n-1},y_n)^{(I,0,k_1,\dots,k_{n-1},1)}}\,
dk.
$$
In the case $Y=*$ this becomes the $S$-dual of the one-point compactification
of ${(I-{\lbrace
0,1\rbrace})}^{(n-1)}$, in other words the $S$-dual of a wedge of
$(n-1)$-spheres freely and transitively permuted by $\Sigma_{n-1}$.

In fact, by taking a different point of view, one can
see that the
$n^{th}$ differential or derivative of the identity has a
$\Sigma_{n+1}$-symmetry and not just a $\Sigma_n$-symmetry, just as if the
identity were itself a derivative.  We hope to return to this point 
in a future paper.

\head\S9. Example: Waldhausen K-theory\endhead

Let $A$ be Waldhausen's algebraic $K$-theory functor from spaces to
spectra \cite{9}. In \cite{1} it was shown that
$$
\partial_yA(Y)
\sim \Sigma^{\infty}_+\Omega Y^.
\tag 9.1
$$
We will give similar formulas for the higher derivatives of $A$.

Really what was shown in \cite{1} was
$$
\partial_yP^{Diff}(Y)
\sim \Omega^2 \Sigma^{\infty}\Omega Y^,
\tag 9.2
$$
where $P^{Diff}$ is stable smooth pseudoisotopy theory. Then 9.1 was 
a corollary
in view of Waldhausen's relation
$$
  A(X)\sim Wh^{Diff}X\times \Sigma^{\infty}_+X.
\tag 9.3
$$
where the Whitehead functor $Wh^{Diff}$ satisfies
$\Omega^2\Omega^{\infty}Wh^{Diff}X\sim P^{Diff}X$.

A natural map
$$
P^{Diff}X@>>>\Omega^2 Q(X^{S^1}/X)
$$
played a key role in obtaining 9.2, where $X^{S^1}/X$ is the quotient 
of the free
loopspace
$X^{S^1}$ by the constant loops. It was then clear that a more direct 
account of 9.1
ought to involve some analogous map
$$
A(X)@>{\tau}>>L(X),
$$
where $L(X)=\Sigma^{\infty}_+X^{S^1}$, a map that ought to have a $K$-theoretic
rather than a manifold-theoretic description. It is not
hard to make such a map (the \it trace\rm\ ), using the methods of \cite{10}.

\remark{9.4 Remark}This was generalized by Bokstedt. The trace map $\tau$ is
reminiscent of the Dennis trace map from the algebraic
$K$-theory of a ring to its Hochschild homology, and this observation 
pointed the
way to the generalization. Since
$A(BG)$ can be defined as the algebraic
$K$-theory of a generalized ring which may be thought of as the group ring of
$G$ over the sphere spectrum, and since the Hochschild homology of a group
ring $k[G]$ is the homology of $(BG)^{S^1}$ with coefficients in $k$, it was
natural to imagine that $\tau$ should be a special case of a construction
$$
K(R)@>>>THH(R)
$$
defined for reasonable ring spectra $R$, where $THH(R)$ is a kind of
\lq\lq Hochschild homology over the sphere spectrum\rq\rq\. Bokstedt
invented (and named) the object $THH$ and defined the trace 
(originally for a class
of generalized rings called \it functors with smash product \rm\ ). See
\cite{12}.
\endremark

The trace induces a
map
$$
\partial_yA(Y)@>{\partial_y \tau}>>\partial_yL(Y)
$$
By 7.1 we have
$$
\partial_yL(Y)
\sim \int^{k\in S^1}\Sigma^{\infty}_+{(Y,y)^{(S^1,k)}}\, dk.
$$
Using rotations of the circle to continuously identify $(S^1,k)$ with 
$(S^1,1)$,
where $1$ is one point in $S^1$, this last spectrum may be identified with
$$
\int^{k\in S^1}\Sigma^{\infty}_+{(Y,y)^{(S^1,1)}}\, dk
\cong Map_*(S^1_+,\Sigma^{\infty}_+\Omega_yY),
$$
and thus split into two factors
$\Sigma^{\infty}_+\Omega_yY\times\Omega\Sigma^{\infty}_+\Omega_yY$.
Projecting on the first factor and composing with $\partial_y \tau$ we get a
map
$$
\partial_yA(Y)@>>>\Sigma^{\infty}_+\Omega_yY,
$$
which is in fact an equivalence.

\remark{9.5 Remark} It would have been very tedious to prove this last fact
directly from 9.2. That would have meant combining the proofs of 9.2 and 9.3
and the construction of the trace and undoubtedly dealing with 
several different
models for
$A(X)$. Instead in \cite{3} we took a shortcut, observing that the functors
$\partial_yA(Y)$ and
$\Sigma^{\infty}_+\Omega_yY$ are abstractly equivalent by 9.1, and then
arguing by universal examples that a natural map between them must be an
equivalence for all $Y$ if this is so in some special cases where things can be
checked.
\endremark

It is possible, and very convenient, to use a homotopy fixed point
construction to single out the factor of $\partial_yL(Y)$ that is to 
correspond to
$\partial_yA(Y)$. Consider the obvious action of the circle group $T$ 
on $L$. The
trace can easily be made to factor through the homotopy fixed point spectrum
$$
A(Y)@>{\tilde\tau}>>L(Y)^{hT}@>>>L(Y),
$$
and it was shown in \cite{3} that the map
$$
\partial_yA(Y)@>>>({\partial_yL(Y)})^{hT}
\tag 9.6
$$
resulting  from this refined trace $\tilde\tau$ is an equivalence. This form of
9.1 will be very useful for getting the higher derivatives.

We emphasize that the right hand side of 9.6 is not
$\partial_y({(L(Y)}^{hT})$, and that $L^{hT}$ is not an analytic functor.
There is a canonical map
$$
\partial_y({(L(Y)}^{hT})@>>>({\partial_yL(Y)})^{hT},
$$
because $T$ acts continuously on the homotopy functor $L$; but this map
does not happen to be an equivalence.

\proclaim{9.7 Theorem}For all
$n\geq 1$ the (natural and symmetry-preserving) map
$$
\partial^{(n)}_{y_1,\dots,y_n}A(Y)@>>>({\partial^{(n)}_{y_1,\dots,y_n}L(Y)})^{hT}
$$
induced by $\tilde\tau$ is an equivalence.
\endproclaim
\demo{Proof}In proving this we are
allowed to ignore the symmetry; thus for this purpose we can get away with
using 6.3 to view
$n^{th}$ derivatives as first derivatives of $(n-1)^{st}$ derivatives.

The key point now is that the operations
$\partial_y$ and
$()^{hT}$, which did not commute when applied to $L(Y)$, do
commute when applied to ${\partial^{(n-1)}_{y_2,\dots,y_n}L(Y)}$. Once this is
established we can argue by induction on $n$:
$$
\partial_{y_1}\partial^{(n-1)}_{y_2,\dots,y_n}A(Y)
\sim \partial_{y_1}(({\partial^{(n-1)}_{y_2,\dots,y_n}L(Y)})^{hT})
\sim (\partial_{y_1}{\partial^{(n-1)}_{y_2,\dots,y_n}L(Y)})^{hT}.
$$
What is needed to establish it is the equation
$$
{\partial^{(n-1)}_{y_2,\dots,y_n}L(Y)}
\sim \int_c^{k\in
(S^1)^{(n-1)}}\Sigma^{\infty}_+{(Y,y_2,\dots,y_n)^{(S^1,k_2,\dots,k_n)}}\,
dk
\tag 9.8
$$
(an instance of 7.10), together with the observation that
$T$ is acting freely on $(S^1)^{(n-1)}$.

Now apply the principle given at the end of \S 5, taking $G$ to be $T$, $P$ to
be
$(S^1)^{(n-1)}$ and
$B$ to be the orbit space. It is clear from 9.6 that for a suitable system $E$,
functorial in
$(Y,y_2,\dots,y_n)$, we have
$$
{\partial^{(n-1)}_{y_2,\dots,y_n}L(Y)}\sim \int^{p\in P}E_{\pi (p)}\, dp
$$
equivariantly and
$$
\align
\partial_y(({\int^{p\in P}E_{\pi (p)}\, dp})^{hT})
&\sim \partial_y \int^{b\in B}E_b\, db\\
&\sim \int^{b\in B}\partial_yE_b\, db\\
&\sim (\int^{p\in P}{\partial_yE_{\pi(p)}\, dp})^{hT}\\
&\sim (\partial_y {\int^{p\in P}E_{\pi (p)}\, dp})^{hT}.
\endalign
$$
The differentiations under the integral are justified by the finiteness
of $B$ and of $P$. (This step fails when $n=1$, basically because the
appropriate $B$ would then be the space $BT$, which is not so finite.)
\qed\enddemo

In the special case $Y=*$ the conclusion of 9.7 is that the $n^{th}$ 
derivative of
$A$ is the homotopy fixed point spectrum for $T$ acting on the $S$-dual of the
one-point compactification of $(S^1)^{(n)}$, where $(S^1)^{(n)}$ has the
obvious commuting actions of
$T$ and $\Sigma_n$. Since the $T$-action is free, the answer may be rewritten
as the $S$-dual of the one-point compactification of $(S^1)^{(n)}/T$. As a
$\Sigma_n$-space, $(S^1)^{(n)}/T$ is isomorphic to
$\Sigma_n\times_{C_n}(S^1\times\bar V_{C_n})$, where the cyclic group
$C_n$ acts on $S^1$ as a subgroup of the rotations and acts linearly on
the vector space $\bar V_{C_n}$ by the reduced regular representation. It
follows that the $n^{th}$ derivative of $A$ at a one-point space is 
induced from
$C_n$ and the $n^{th}$ layer of the Taylor tower is given by
$$
(D_nA)(X)\sim
Map_*( {S^1}_+\wedge S^{\bar V_{C_n}} , \Sigma^{\infty}(X^{\wedge n}) )^{C_n}.
$$

\remark{9.9 Remark}In the special case when $X$ is a suspension $\Sigma Y$
this becomes $\Sigma Map_*({S^1}_+,\Sigma^{\infty}(X^{\wedge
n}))^{C_n}$. In fact, in that case the Taylor tower splits:
$$
A(\Sigma Y)\sim A(*)\times \prod_{n\geq 1}
\Sigma
Map_*({S^1}_+,\Sigma^{\infty}(X^{\wedge n}))^{C_n}
$$
if $Y$ is connected. A proof of this was sketched in \cite{11} and corrected in
\cite{3}. This special case of the conclusion was used in the shortcut
mentioned in 9.5.
\endremark

\remark{9.10 Remark}There is also a direct $K$-theoretic approach to
all of this. Dundas and McCarthy \cite{15} generalized 9.1 to 
(generalized) rings.
Their statement is that the trace map induces an equivalence from \lq\lq stable
$K$-theory\rq\rq\ to $THH$. This led in \cite{16} to an extension of the main
result of \cite{3} to such generalized rings.

\end

%% file: gtmacros.tex
%
%
%
%
%
%
\magnification=\magstephalf      
%
%
\vsize=7.5truein                 
\hsize=5.2truein                 
\newskip\stdskip                 
\stdskip=6pt plus3pt minus3pt    
\medskipamount=\stdskip          
\parindent=0pt                   
\parskip=\stdskip                
\abovedisplayskip=\stdskip       
\belowdisplayskip=\stdskip       
\mathsurround=0.75pt             
\overfullrule=0pt                
%
%
\def\ppar{\par\goodbreak\vskip 8pt plus 4pt minus 4pt}     
%
%
\def\stdspace{\hskip 0.75em plus 0.15em\ignorespaces}
\let\qua\stdspace 
%
%
%
%
%
%
%
\def\hexnumber#1{\ifcase#1 0\or 1\or 2\or 3\or 4\or 5\or 6\or 7\or 8\or
 9\or A\or B\or C\or D\or E\or F\fi}
%
%
\font\thirtnmsa=msam10 scaled 1315    
\font\tenmsa=msam10          \font\ninemsa=msam9
\font\sevenmsa=msam7         \font\sixmsa=msam6
\font\fivemsa=msam5
%
%
\newfam\msafam                  \textfont\msafam=\tenmsa
\scriptfont\msafam=\sevenmsa    \scriptscriptfont\msafam=\fivemsa
\edef\hexa{\hexnumber\msafam}        
\def\msa{\fam\msafam\tenmsa}         
%
%
\font\thirtnmsb=msbm10 scaled 1315   
\font\tenmsb=msbm10      \font\ninemsb=msbm9
\font\sevenmsb=msbm7     \font\sixmsb=msbm6
\font\fivemsb=msbm5
%
\newfam\msbfam                   \textfont\msbfam=\tenmsb       
\scriptfont\msbfam=\sevenmsb     \scriptscriptfont\msbfam=\fivemsb
\edef\hexb{\hexnumber\msbfam}    
\def\msb{\fam\msbfam\tenmsb}     
%
%
\font\thirtneufm=eufm10 scaled 1315   
\font\teneufm=eufm10                 \font\nineeufm=eufm9
\font\seveneufm=eufm7                \font\sixeufm=eufm6
\font\fiveeufm=eufm5
%
\newfam\eufmfam                    \textfont\eufmfam=\teneufm
\scriptfont\eufmfam=\seveneufm     \scriptscriptfont\eufmfam=\fiveeufm
\edef\hexf{\hexnumber\eufmfam}      
\def\frak{\fam\eufmfam\teneufm}     
%
%
%
\font\thirtnrm=cmr10 scaled 1315    
\font\ninerm=cmr9                   \font\sixrm=cmr6   
%
\font\thirtni=cmmi10 scaled 1315    
\font\ninei=cmmi9                   \font\sixi=cmmi6  
%
\font\thirtnsy=cmsy10 scaled 1315   
\font\ninesy=cmsy9                  \font\sixsy=cmsy6  
%
\font\thirtnbf=cmbx10 scaled 1315   
\font\ninebf=cmbx9                  \font\sixbf=cmbx6  
%
%
\font\thirtnex=cmex10 scaled 1315   
\font\nineex=cmex9                  
%
%
\font\thirtnit=cmti10 scaled 1315  
\font\nineit=cmti9                  
%
\font\thirtnsl=cmsl10 scaled 1315  
\font\ninesl=cmsl9                  
%
\font\thirtntt=cmtt10 scaled 1315  
\font\ninett=cmtt9                  
%
%
%
%
\def\small{%
%
%
\textfont0=\ninerm \scriptfont0=\sixrm \scriptscriptfont0=\fiverm
\def\rm{\fam0\ninerm}
%
%
\textfont1=\ninei \scriptfont1=\sixi \scriptscriptfont1=\fivei
%
%
\textfont2=\ninesy \scriptfont2=\sixsy \scriptscriptfont2=\fivesy
%
%
\textfont3=\nineex \scriptfont3=\nineex \scriptscriptfont3=\nineex
%
%
\textfont\bffam=\ninebf \scriptfont\bffam=\sixbf
\scriptscriptfont\bffam=\fivebf \def\bf{\fam\bffam\ninebf}%
%
%
\textfont\itfam=\nineit \def\it{\fam\itfam\nineit}%
\textfont\slfam=\ninesl \def\sl{\fam\slfam\ninesl}%
\textfont\ttfam=\ninett \def\tt{\fam\ttfam\ninett}%
%
%
%
\textfont\msafam=\ninemsa \scriptfont\msafam=\sixmsa
\scriptscriptfont\msafam=\fivemsa \def\msa{\fam\msafam\ninemsa}%
%
%
\textfont\msbfam=\ninemsb \scriptfont\msbfam=\sixmsb
\scriptscriptfont\msbfam=\fivemsb \def\msb{\fam\msbfam\ninemsb}%
%
%
\textfont\eufmfam=\nineeufm  \scriptfont\eufmfam=\sixeufm
\scriptscriptfont\eufmfam=\fiveeufm \def\frak{\fam\eufmfam\nineeufm}%
%
%
%
\normalbaselineskip=11pt%
\setbox\strutbox=\hbox{\vrule height8pt depth3pt width0pt}%
%
%
\normalbaselines\rm
%
%
\stdskip=4pt plus2pt minus2pt    
\medskipamount=\stdskip          
\parskip=\stdskip                
\abovedisplayskip=\stdskip       
\belowdisplayskip=\stdskip       
\def\ppar{\par\goodbreak\vskip 6pt plus 3pt minus 3pt}%
%
%
\def\section##1{\global\advance\sectionnumber by 1
\vskip-\lastskip\penalty-800\vskip 20pt plus10pt minus5pt 
\egroup{\bf\number\sectionnumber\quad##1}\bgroup\small         
\vskip 6pt plus3pt minus3pt
\nobreak\resultnumber=1}
}    
%
\def\beginsmall{\bgroup\small}
\let\endsmall\egroup
%
%
%
%
\def\large{%
\textfont0=\thirtnrm \scriptfont0=\ninerm \scriptscriptfont0=\sevenrm
\def\rm{\fam0\thirtnrm}%
\textfont1=\thirtni \scriptfont1=\ninei \scriptscriptfont1=\seveni
\textfont2=\thirtnsy \scriptfont2=\ninesy \scriptscriptfont2=\sevensy
\textfont3=\thirtnex \scriptfont3=\thirtnex \scriptscriptfont3=\thirtnex
\textfont\bffam=\thirtnbf \scriptfont\bffam=\ninebf
\scriptscriptfont\bffam=\sevenbf \def\bf{\fam\bffam\thirtnbf}%
\textfont\itfam=\thirtnit \def\it{\fam\itfam\thirtnit}%
\textfont\slfam=\thirtnsl \def\sl{\fam\slfam\thirtnsl}%
\textfont\ttfam=\thirtntt \def\tt{\fam\ttfam\thirtntt}%
\textfont\msafam=\thirtnmsa \scriptfont\msafam=\ninemsa
\scriptscriptfont\msafam=\sevenmsa \def\msa{\fam\msafam\thirtnmsa}%
\textfont\msbfam=\thirtnmsb \scriptfont\msbfam=\ninemsb
\scriptscriptfont\msbfam=\sevenmsb \def\msb{\fam\msbfam\thirtnmsb}%
\textfont\eufmfam=\thirtneufm  \scriptfont\eufmfam=\nineeufm
\scriptscriptfont\eufmfam=\seveneufm \def\frak{\fam\eufmfam\teneufm}%
\normalbaselineskip=16pt%
\setbox\strutbox=\hbox{\vrule height11.5pt depth4.5pt width0pt}%
\normalbaselines\rm}%
\let\Large\large   
%
\def\Bbb#1{{\msb#1}}

%

%
\mathchardef\plussquare="0\hexa01
\mathchardef\nge="3\hexb0B
\mathchardef\maltesecross="0\hexa7A
\mathchardef\del="0\hexf01
%
%
%
%
\font\sc=cmcsc10
%
%
%
%
\def\sqr#1#2{{\vcenter{\vbox{\hrule  height.#2truept
	\hbox{\vrule width.#2truept height#1truept 
	\kern#1truept \vrule width.#2truept}
	\hrule height.#2truept}}}}
\def\sq{\sqr55}    
%
%
%
%
\newcount\sectionnumber            
\newcount\resultnumber             
\sectionnumber=0\resultnumber=1    
%
%
%
\def\section#1{\global\advance\sectionnumber by 1
\xdef\nextkey{\number\sectionnumber}
\vskip-\lastskip\penalty-800\vskip 20pt plus10pt minus5pt 
{\large\bf\number\sectionnumber\quad#1}         
\vskip 8pt plus4pt minus4pt
\nobreak\resultnumber=1}      
%
%
%
%
%
\def\sh#1{\vskip-\lastskip\ppar{\bf #1}\par\nobreak\medskip}         
%
%
%
%

%
\def\proc#1{\xdef\nextkey{\number\sectionnumber.\number\resultnumber}%
\vskip-\lastskip\ppar\bf%
\noindent#1\ \number\sectionnumber.\number\resultnumber
\stdspace\sl\global\advance\resultnumber by 1\ignorespaces}
 
%
%
\def\qed{\hfill$\sq$\par\goodbreak\rm}   
%
%
%
%
%
%
%
%
\def\proclaim#1{\vskip-\lastskip\ppar\bf%
\noindent#1\stdspace\sl\ignorespaces} 
\let\endproclaim\endproc
%
%
%
%
\def\rk#1{\vskip-\lastskip\ppar{\bf #1}\stdspace\ignorespaces}                

%
%
%
%
%
%
\def\label{\xdef\nextkey{\number\sectionnumber.\number\resultnumber}%
\number\sectionnumber.\number\resultnumber
\global\advance\resultnumber by 1}
%
%
%
%
%
%
%
%
%
%
%
%
%
%
%
%
\newcount\refnumber              
\refnumber=1                     
\long\def\reflist#1\endreflist{%
\long\def\thereflist{#1}{\def\refkey##1##2\par{\xdef##1{\number\refnumber}%
\global\advance\refnumber by 1}%
\def\key##1##2\par{\expandafter\xdef%
\csname##1\endcsname{\number\refnumber}%
\global\advance\refnumber by 1}#1\par}}
\long\def\references{%
\penalty-800\vskip-\lastskip\vskip 15pt plus10pt minus5pt 
{\large\bf References}\ppar 
{\leftskip=25pt\frenchspacing    
\small\parskip=3pt plus2pt       
\def\refkey##1##2\par{\noindent  
\llap{[##1]\stdspace}\ignorespaces##2\par}         
\def\key##1##2\par{\noindent  
\llap{[\ref{##1}]\stdspace}\ignorespaces##2\par}  
\def\,{\thinspace}\thereflist\par}}
%
%
%
\newcount\footnotenumber         
\footnotenumber=1                
\def\fnote#1{\xdef\nextkey{\number\footnotenumber}%
{\small\ifnum\footnotenumber>9\parindent=14pt%
\else\parindent=10pt\fi\footnote{$^{\number\footnotenumber}$}%
{\hglue-5pt#1}\global\advance\footnotenumber by 1}}
%
%
%
%
%
%
%
\newcount\figurenumber          
\figurenumber=1                 
\def\caption#1{\xdef\nextkey{\number\figurenumber}%
\cl{\small Figure \number\figurenumber: #1}%
\global\advance\figurenumber by 1}
\def\figurelabel{\xdef\nextkey{\number\figurenumber}%
\cl{\small Figure \number\figurenumber}%
\global\advance\figurenumber by 1}
\long\def\figure#1\endfigure{{\xdef\nextkey{\number\figurenumber}%
\let\captiontext\relax\def\caption##1{\xdef\captiontext{##1}}%
\midinsert\cl{\ignorespaces#1\unskip\unskip\unskip\unskip}\vglue6pt\cl{\small 
Figure \number\figurenumber\ifx\captiontext\relax\else: \captiontext
\fi}\endinsert\global\advance\figurenumber by 1}}
%
%
%
%
%
%
%
\def\nextkey{??}   
%
\def\key#1{\expandafter\xdef\csname #1\endcsname{\nextkey}}
\def\ref#1{\expandafter\ifx\csname #1\endcsname\relax
\immediate\write16{Reference {#1} undefined}??\else
\csname #1\endcsname\fi}
%
%
%
%
%
%
%
\newread\gtinfile
\newwrite\gtreffile
\def\useforwardrefs{
\openin\gtinfile\jobname.ref
\ifeof\gtinfile
\closein\gtinfile
\immediate\write16{No file \jobname.ref}
\else
\closein\gtinfile
\input \jobname.ref
\fi
\immediate\openout\gtreffile \jobname.ref
%
%
\def\key##1{{\def\\{\noexpand}%
\expandafter\xdef\csname ##1\endcsname{\nextkey}%
\immediate\write\gtreffile{\\\expandafter\\\def\\\csname ##1\\\endcsname%
{\nextkey}}}}
%
%
\long\def\reflist##1\endreflist{%
\long\def\thereflist{##1}{\def\refkey####1####2\par{\xdef####1{%
\number\refnumber}{\def\\{\noexpand}\immediate\write\gtreffile
{\\\def\\####1{\number\refnumber}}}\global\advance\refnumber by 1}%
\def\key####1####2\par{\expandafter\xdef%
\csname####1\endcsname{\number\refnumber}%
{\def\\{\noexpand}\immediate\write\gtreffile
{\\\expandafter\\\def\\\csname ####1\\\endcsname{\number\refnumber}}}
\global\advance\refnumber by 1}##1\par}}
\long\def\biblio##1\endbiblio{\reflist##1\endreflist\references}%
%
%
\def\numkey##1{{\def\\{\noexpand}%
\xdef##1{\number\sectionnumber.\number\resultnumber}
\immediate\write\gtreffile{\\\def\\##1%
{\number\sectionnumber.\number\resultnumber}}}}
\def\seckey##1{{\def\\{\noexpand}\xdef##1{\number\sectionnumber}
\immediate\write\gtreffile{\\\def\\##1{\number\sectionnumber}}}}
\def\figkey##1{\xdef##1{\number\figurenumber}%
{\def\\{\noexpand}\immediate\write\gtreffile%
{\\\def\\##1{\number\figurenumber}}}
\number\figurenumber\global\advance\figurenumber by 1}
}   
%
%
%
%
\def\figkey#1{\xdef#1{\number\figurenumber}%
\number\figurenumber\global\advance\figurenumber by 1}
\def\fig#1#2\endfig{%
\midinsert\cl{#2}\vglue6pt\cl{\small Figure #1}\endinsert}
\def\newfig{\number\figurenumber\global\advance\figurenumber by 1}
\def\numkey#1{\xdef#1{\number\sectionnumber.\number\resultnumber}}
\def\seckey#1{\xdef#1{\number\sectionnumber}}
%
%
%
%
%
%
%
%
%
\def\verb{\catcode`\"=\active}       
\def\brev{\catcode`\"=12}            
\brev                                
\verb                                
{\obeyspaces\gdef {\ }}              
{\catcode`\`=\active\gdef`{\relax\lq}}
\def"{%
\begingroup\baselineskip=12pt\def\par{\leavevmode\endgraf}%
\tt\obeylines\obeyspaces\parskip=0pt\parindent=0pt%
\catcode`\$=12\catcode`\&=12\catcode`\^=12\catcode`\#=12%
\catcode`\_=12\catcode`\~=12%
\catcode`\{=12\catcode`\}=12\catcode`\%=12\catcode`\\=12%
\catcode`\`=\active\let"\endgroup}
\brev      
%
%
%
%
%
%
\def\items{\par\leftskip = 25pt}           
\def\enditems{\par\leftskip = 0pt}         
\def\item#1{\par\leavevmode\llap{#1\stdspace}%
\ignorespaces}                             
%
%

%
%
\def\np{\vfil\eject}         
\def\nl{\hfil\break}         
\def\cl{\centerline}         
\def\gt{{\mathsurround=0pt\it $\cal G\mskip-2mu$eometry \&\ 
$\cal T\!\!$opology}}        
\def\agt{{\mathsurround=0pt\it$\cal A\mskip-.7mu$lgebraic \&\ 
$\cal G\mskip-2mu$eometric $\cal T\!\!$opology}}  
%
%
%

%
%
%
%
%
\def\title#1{\def\thetitle{#1}}

\def\author#1{\edef\previousauthors{\theauthors}
 \ifx\theauthors\relax\def\theauthors{#1}\else
 \def\theauthors{\previousauthors\par#1}\fi}

%
\def\address#1{\edef\previousaddresses{\theaddress}
 \ifx\theaddress\relax\def\theaddress{#1}\else
 \def\theaddress{\previousaddresses\par\vskip 2pt\par#1}\fi}
\def\secondaddress#1{\edef\previousaddresses{\theaddress}
 \ifx\theaddress\relax\def\theaddress{#1}\else
 \def\theaddress{\previousaddresses\par{\rm and}\par#1}\fi}   

\def\email#1{\edef\previousemails{\theemail}
 \ifx\theemail\relax\def\theemail{#1}\else
 \def\theemail{\previousemails\hskip 0.75em\relax#1}\fi}
\def\secondemail#1{\edef\previousemails{\theemail}
 \ifx\theemail\relax\def\theemail{#1}\else
 \def\theemail{\previousemails\hskip 0.75em{\rm and}\hskip 0.75em
 \relax#1}\fi}
\def\url#1{\edef\previousurls{\theurl}
 \ifx\theurl\relax\def\theurl{#1}\else
 \def\theurl{\previousurls\hskip 0.75em\relax#1}\fi}
\def\secondurl#1{\edef\previousurls{\theurl}
 \ifx\theurl\relax\def\theurl{#1}\else
 \def\theurl{\previousurls\hskip 0.75em{\rm and}\hskip 0.75em
 \relax#1}\fi}
\long\def\abstract#1\endabstract{\long\def\theabstract{#1}}
\def\primaryclass#1{\def\theprimaryclass{#1}}
\let\subjclass\primaryclass                        
\def\secondaryclass#1{\def\thesecondaryclass{#1}}
\def\keywords#1{\def\thekeywords{#1}}
%
%
\let\\\par\let\thetitle\relax\let\theshorttitle\relax
\let\theauthors\relax\let\theshortauthors\relax
\let\theaddress\relax\let\theshortaddress\relax
\let\theemail\relax\let\theurl\relax
\let\theabstract\relax\let\theprimaryclass\relax
\let\thesecondaryclass\relax\let\thekeywords\relax
%
%
%
%
\long\def\maketitlepage{    

\vglue 0.2truein   

%
{\parskip=0pt\leftskip 0pt plus 1fil\def\\{\par\smallskip}{\large
\bf\thetitle}\par\medskip}   

\vglue 0.15truein 

%
{\parskip=0pt\leftskip 0pt plus 1fil\def\\{\par}{\sc\theauthors}
\par\medskip}%
 
\vglue 0.1truein 

%
{\small\parskip=0pt
{\leftskip 0pt plus 1fil\def\\{\par}{\sl\theaddress}\par}
\ifx\theemail\relax\else  
\vglue 5pt \def\\{\stdspace{\rm and}\stdspace} 
\cl{Email:\stdspace\tt\theemail}\fi
\ifx\theurl\relax\else    
\vglue 5pt \def\\{\stdspace{\rm and}\stdspace} 
\cl{URL:\stdspace\tt\theurl}\fi\par}

\vglue 7pt 

{\bf Abstract}

\vglue 5pt

\theabstract

\vglue 7pt 

{\bf AMS Classification numbers}\quad Primary:\quad \theprimaryclass\par

Secondary:\quad \thesecondaryclass

\vglue 5pt 

{\bf Keywords:}\quad \thekeywords

\np  

}    
%
%
\long\def\makeshorttitle{    


%
{\parskip=0pt\leftskip 0pt plus 1fil\def\\{\par\smallskip}{\large
\bf\thetitle}\par\medskip}   

\vglue 0.05truein 

%
{\parskip=0pt\leftskip 0pt plus 1fil\def\\{\par}{\sc\theauthors}
\par\medskip}%
 
\vglue 0.03truein 

%
{\small\parskip=0pt
{\leftskip 0pt plus 1fil\def\\{\par}{\sl\ifx\theshortaddress\relax
\theaddress\else\theshortaddress\fi}\par}
\ifx\theemail\relax\else  
\vglue 5pt \def\\{\stdspace{\rm and}\stdspace} 
\cl{Email:\stdspace\tt\theemail}\fi
\ifx\theurl\relax\else    
\vglue 5pt \def\\{\stdspace{\rm and}\stdspace} 
\cl{URL:\stdspace\tt\theurl}\fi\par}

\vglue 10pt 


{\small\leftskip 25pt\rightskip 25pt{\bf Abstract}\stdspace\theabstract

{\bf AMS Classification}\stdspace\theprimaryclass
\ifx\thesecondaryclass\relax\else; \thesecondaryclass\fi\par
{\bf Keywords}\stdspace \thekeywords\par}
\vglue 7pt
}    
\let\maketitle\makeshorttitle        
%
%

\def\volumenumber#1{\def\thevolumenumber{#1}}
\def\volumeyear#1{\def\thevolumeyear{#1}}
\def\pagenumbers#1#2{\def\startpage{#1}\def\finishpage{#2}}
\def\published#1{\def\publishdate{#1}}
\def\received#1{\def\receiveddate{#1}}

\let\reviseddate\relax
\volumenumber{X}
\volumeyear{20XX}
\pagenumbers{1}{XXX}
\published{XX Xxxember 20XX}

\long\def\makeagttitle{   
\agt\hfill      
\hbox to 60truept{\vbox to 0pt{\vglue -14truept{\bf [Logo here]}\vss}\hss}
\break
{\small Volume \thevolumenumber\ (\thevolumeyear)
\startpage--\finishpage\nl
Published: \publishdate}

\vglue .2truein

{\parskip=0pt\leftskip 0pt plus 1fil\def\\{\par\smallskip}{\large
\bf\thetitle}\par\medskip}   
\vglue 0.05truein 

%
{\parskip=0pt\leftskip 0pt plus 1fil\def\\{\par}{\sc\theauthors}
\par\medskip}%
 
\vglue 0.03truein 


{\small\leftskip 25truept\rightskip 25truept{\bf Abstract}\stdspace\theabstract

{\bf AMS Classification}\stdspace\theprimaryclass
\ifx\thesecondaryclass\relax\else; \thesecondaryclass\fi\par
{\bf Keywords}\stdspace \thekeywords\par}\vglue 7truept

}   


\def\Addresses{\bigskip
{\small \parskip 0pt \leftskip 0pt \rightskip 0pt plus 1fil \def\\{\par}
\sl\theaddress\par\medskip \rm Email:\stdspace\tt\theemail\par
\ifx\theurl\relax\else\smallskip \rm URL:\stdspace\tt\theurl\par\fi}}

\def\agtart{
\hoffset 14truemm
\voffset 31truemm
\font\phead=cmsl9 scaled 950
\font\pnum=cmbx10 scaled 913
\font\pfoot=cmsl9 scaled 950
\headline{\vbox to 0pt{\vskip -4.5mm\line{\small\phead\ifnum
\count0=\startpage ISSN numbers are printed here
\hfill {\pnum\folio}\else\ifodd\count0\def\\{ }%
\ifx\theshorttitle\relax\thetitle\else\theshorttitle\fi\hfill{\pnum\folio}
\else\def\\{ and }{\pnum\folio}\hfill\ifx\theshortauthors\relax\theauthors
\else\theshortauthors\fi\fi\fi}\vss}}
\footline{\vbox to 0pt{\vglue 0mm\line{\small\pfoot\ifnum\count0=\startpage
Copyright declaration is printed here\hfill\else
\agt, Volume \thevolumenumber\ (\thevolumeyear)\hfill\fi}\vss}}
\let\maketitle\makeagttitle\let\makeshorttitle\makeagttitle}

%% file: amsnames.tex
%
%
%
%
%
%
%
\def\newsymbol#1#2#3#4#5{%
\ifodd#2\let\hexno\hexa\else\let\hexno\hexb\fi
\mathchardef#1="#3\hexno#4#5}
%
%
\def\mathsurroundoff{\mathsurround=0pt}
\def\mathhexbox#1#2#3{\relax
\ifmmode\mathpalette{}{\mathsurroundoff\mathchar"#1#2#3}%
\else\leavevmode\hbox{$\mathsurroundoff\mathchar"#1#2#3$}\fi}
%
%
\mathchardef\dashbar"0\hexa39

%
%
\newsymbol\boxdot 1200
\newsymbol\boxplus 1201
\newsymbol\boxtimes 1202
\newsymbol\square 1003
\newsymbol\blacksquare 1004
\newsymbol\centerdot 1205
\newsymbol\lozenge 1006
\newsymbol\blacklozenge 1007
\newsymbol\circlearrowright 1308
\newsymbol\circlearrowleft 1309
\let\rightleftharpoons\undefined
\newsymbol\rightleftharpoons 130A
\newsymbol\leftrightharpoons 130B
\newsymbol\boxminus 120C
\newsymbol\Vdash 130D
\newsymbol\Vvdash 130E
\newsymbol\vDash 130F
\newsymbol\twoheadrightarrow 1310
\newsymbol\twoheadleftarrow 1311
\newsymbol\leftleftarrows 1312
\newsymbol\rightrightarrows 1313
\newsymbol\upuparrows 1314
\newsymbol\downdownarrows 1315
\newsymbol\upharpoonright 1316
 
\newsymbol\downharpoonright 1317
\newsymbol\upharpoonleft 1318
\newsymbol\downharpoonleft 1319
\newsymbol\rightarrowtail 131A
\newsymbol\leftarrowtail 131B
\newsymbol\leftrightarrows 131C
\newsymbol\rightleftarrows 131D
\newsymbol\Lsh 131E
\newsymbol\Rsh 131F
\newsymbol\rightsquigarrow 1320
\newsymbol\leftrightsquigarrow 1321
\newsymbol\looparrowleft 1322
\newsymbol\looparrowright 1323
\newsymbol\circeq 1324
\newsymbol\succsim 1325
\newsymbol\gtrsim 1326
\newsymbol\gtrapprox 1327
\newsymbol\multimap 1328
\newsymbol\therefore 1329
\newsymbol\because 132A
\newsymbol\doteqdot 132B
 
\newsymbol\triangleq 132C
\newsymbol\precsim 132D
\newsymbol\lesssim 132E
\newsymbol\lessapprox 132F
\newsymbol\eqslantless 1330
\newsymbol\eqslantgtr 1331
\newsymbol\curlyeqprec 1332
\newsymbol\curlyeqsucc 1333
\newsymbol\preccurlyeq 1334
\newsymbol\leqq 1335
\newsymbol\leqslant 1336
\newsymbol\lessgtr 1337
\newsymbol\backprime 1038
\newsymbol\risingdotseq 133A
\newsymbol\fallingdotseq 133B
\newsymbol\succcurlyeq 133C
\newsymbol\geqq 133D
\newsymbol\geqslant 133E
\newsymbol\gtrless 133F
\newsymbol\sqsubset 1340
\newsymbol\sqsupset 1341
\newsymbol\vartriangleright 1342
\newsymbol\vartriangleleft 1343
\newsymbol\trianglerighteq 1344
\newsymbol\trianglelefteq 1345
\newsymbol\bigstar 1046
\newsymbol\between 1347
\newsymbol\blacktriangledown 1048
\newsymbol\blacktriangleright 1349
\newsymbol\blacktriangleleft 134A
\newsymbol\vartriangle 134D
\newsymbol\blacktriangle 104E
\newsymbol\triangledown 104F
\newsymbol\eqcirc 1350
\newsymbol\lesseqgtr 1351
\newsymbol\gtreqless 1352
\newsymbol\lesseqqgtr 1353
\newsymbol\gtreqqless 1354
\newsymbol\Rrightarrow 1356
\newsymbol\Lleftarrow 1357
\newsymbol\veebar 1259
\newsymbol\barwedge 125A
\newsymbol\doublebarwedge 125B
\let\angle\undefined
\newsymbol\angle 105C
\newsymbol\measuredangle 105D
\newsymbol\sphericalangle 105E
\newsymbol\varpropto 135F
\newsymbol\smallsmile 1360
\newsymbol\smallfrown 1361
\newsymbol\Subset 1362
\newsymbol\Supset 1363
\newsymbol\Cup 1264
 
\newsymbol\Cap 1265
 
\newsymbol\curlywedge 1266
\newsymbol\curlyvee 1267
\newsymbol\leftthreetimes 1268
\newsymbol\rightthreetimes 1269
\newsymbol\subseteqq 136A
\newsymbol\supseteqq 136B
\newsymbol\bumpeq 136C
\newsymbol\Bumpeq 136D
\newsymbol\lll 136E
 
\newsymbol\ggg 136F
 
\newsymbol\circledS 1073
\newsymbol\pitchfork 1374
\newsymbol\dotplus 1275
\newsymbol\backsim 1376
\newsymbol\backsimeq 1377
\newsymbol\complement 107B
\newsymbol\intercal 127C
\newsymbol\circledcirc 127D
\newsymbol\circledast 127E
\newsymbol\circleddash 127F
\newsymbol\lvertneqq 2300
\newsymbol\gvertneqq 2301
\newsymbol\nleq 2302
\newsymbol\ngeq 2303
\newsymbol\nless 2304
\newsymbol\ngtr 2305
\newsymbol\nprec 2306
\newsymbol\nsucc 2307
\newsymbol\lneqq 2308
\newsymbol\gneqq 2309
\newsymbol\nleqslant 230A
\newsymbol\ngeqslant 230B
\newsymbol\lneq 230C
\newsymbol\gneq 230D
\newsymbol\npreceq 230E
\newsymbol\nsucceq 230F
\newsymbol\precnsim 2310
\newsymbol\succnsim 2311
\newsymbol\lnsim 2312
\newsymbol\gnsim 2313
\newsymbol\nleqq 2314
\newsymbol\ngeqq 2315
\newsymbol\precneqq 2316
\newsymbol\succneqq 2317
\newsymbol\precnapprox 2318
\newsymbol\succnapprox 2319
\newsymbol\lnapprox 231A
\newsymbol\gnapprox 231B
\newsymbol\nsim 231C
\newsymbol\ncong 231D
\newsymbol\diagup 201E
\newsymbol\diagdown 201F
\newsymbol\varsubsetneq 2320
\newsymbol\varsupsetneq 2321
\newsymbol\nsubseteqq 2322
\newsymbol\nsupseteqq 2323
\newsymbol\subsetneqq 2324
\newsymbol\supsetneqq 2325
\newsymbol\varsubsetneqq 2326
\newsymbol\varsupsetneqq 2327
\newsymbol\subsetneq 2328
\newsymbol\supsetneq 2329
\newsymbol\nsubseteq 232A
\newsymbol\nsupseteq 232B
\newsymbol\nparallel 232C
\newsymbol\nmid 232D
\newsymbol\nshortmid 232E
\newsymbol\nshortparallel 232F
\newsymbol\nvdash 2330
\newsymbol\nVdash 2331
\newsymbol\nvDash 2332
\newsymbol\nVDash 2333
\newsymbol\ntrianglerighteq 2334
\newsymbol\ntrianglelefteq 2335
\newsymbol\ntriangleleft 2336
\newsymbol\ntriangleright 2337
\newsymbol\nleftarrow 2338
\newsymbol\nrightarrow 2339
\newsymbol\nLeftarrow 233A
\newsymbol\nRightarrow 233B
\newsymbol\nLeftrightarrow 233C
\newsymbol\nleftrightarrow 233D
\newsymbol\divideontimes 223E
\newsymbol\varnothing 203F
\newsymbol\nexists 2040
\newsymbol\Finv 2060
\newsymbol\Game 2061
\newsymbol\mho 2066
\newsymbol\eth 2067
\newsymbol\eqsim 2368
\newsymbol\beth 2069
\newsymbol\gimel 206A
\newsymbol\daleth 206B
\newsymbol\lessdot 236C
\newsymbol\gtrdot 236D
\newsymbol\ltimes 226E
\newsymbol\rtimes 226F
\newsymbol\shortmid 2370
\newsymbol\shortparallel 2371
\newsymbol\smallsetminus 2272
\newsymbol\thicksim 2373
\newsymbol\thickapprox 2374
\newsymbol\approxeq 2375
\newsymbol\succapprox 2376
\newsymbol\precapprox 2377
\newsymbol\curvearrowleft 2378
\newsymbol\curvearrowright 2379
\newsymbol\digamma 207A
\newsymbol\varkappa 207B
\newsymbol\Bbbk 207C
\newsymbol\hslash 207D
\let\hbar\undefined
\newsymbol\hbar 207E
\newsymbol\backepsilon 237F

%% file: gtoutput.tex

\def\ifplaintex{\expandafter\ifx\csname documentclass\endcsname\relax}


\ifplaintex 
\hoffset 14truemm
\voffset 31truemm
\else
\headsep 23pt
\footskip 35pt
\hoffset -4truemm
\voffset 12.5truemm
\fi

\expandafter\ifx\csname beginpicture\endcsname\relax
\expandafter\ifx\csname documentclass\endcsname\relax
\input pictex \else\font\fiverm=cmr5
\input prepictex \input pictex \input postpictex \fi\fi

\def\gt{{\mathsurround=0pt\it $\cal G\mskip-2mu$eometry \&\ 
$\cal T\!\!$opology}}        

\def\gtp{{\mathsurround=0pt\it $\cal G\mskip-2mu$eometry \&\ 
$\cal T\!\!$opology $\cal P\!$ublications}}  


\def\lognumber#1{\def\thelognumber{#1}}
\def\volumenumber#1{\def\thevolumenumber{#1}}
\def\papernumber#1{\def\thepapernumber{#1}}
\def\volumeyear#1{\def\thevolumeyear{#1}}

\def\pagenumbers#1#2{\def\startpage{#1}\def\finishpage{#2}}
\def\published#1{\def\publishdate{#1}}
\def\proposed#1{\def\theproposer{#1}}
\def\seconded#1{\def\theseconders{#1}}
\def\received#1{\def\receiveddate{#1}}

\def\accepted#1{\def\accepteddate{#1}}

\long\def\asciiabstract#1{\long\def\theasciiabstract{#1}}


\let\\\par\let\thelognumber\relax
\let\thevolumenumber\relax\let\thepapernumber\relax
\let\thevolumeyear\relax\let\thesamplenumber\relax\let\startpage\relax
\let\finishpage\relax\let\publishdate\relax\let\receiveddate\relax
\let\reviseddate\relax\let\accepteddate\relax\let\theasciititle\relax
\let\theasciiauthors\relax
\let\theasciiabstract\relax
\let\theasciiemail\relax\let\theshortauthors\relax\let\theshorttitle\relax

\long\def\maketitlep{   

\count0=\startpage

\gt\hfill      
\beginpicture
\setcoordinatesystem units <0.33truein, 0.33truein> point at 2.2 0.9
\setplotsymbol ({$\cal G$})
\plotsymbolspacing=9truept
\circulararc 315 degrees from 0 1 center at 0 0
\setplotsymbol ({$\cal T$})
\circulararc 315 degrees from 1 -1 center at 1 0
\endpicture
%
\break
{\small\ifx\thesamplenumber\relax 
Volume \else Sample
\fi\thevolumenumber\ (\thevolumeyear)
\startpage--\finishpage\nl
Published: \publishdate}
\vglue 0.5truein plus 0.4fil minus 0.1truein

{\parskip=0pt\leftskip 0pt plus 1fil\def\\{\par\smallskip}{\ifplaintex\large
\else\Large\fi\bf\thetitle}\par\medskip}   

\vglue 0pt plus 0.1fil 

{\parskip=0pt\leftskip 0pt plus 1fil\def\\{\par}{\sc\theauthors}
\par\medskip}

\vglue 0pt plus 0.1fil 

{\small\parskip=0pt\let\newline\\
{\leftskip 0pt plus 1fil\def\\{\par}{\sl\theaddress}\par}
\expandafter\ifx\theemail\relax    
\relax\else\vglue 5pt plus 0.02fil minus 2pt\def\\{\stdspace{\rm 
and}\stdspace} 
\cl{Email:\stdspace\tt\theemail}\fi
\ifx\theurl\relax                  
\relax\else\vglue 5pt plus 0.02fil minus 2pt\def\\{\stdspace{\rm 
and}\stdspace}
\cl{URL:\stdspace\tt\theurl}\fi\par}

\vglue 7pt plus 0.3fil minus 3pt

{\bf Abstract}
\vglue 5pt plus 0.1fil minus 2pt

\theabstract

\vglue 7pt plus 0.3fil minus 3pt

{\bf AMS Classification numbers}\quad Primary:\quad \theprimaryclass

Secondary:\quad \thesecondaryclass

\vglue 5pt plus 0.3fil minus 2pt

{\bf Keywords}\quad \thekeywords

\vglue 10pt plus 0.5fil minus 5pt

{\small  Proposed: \theproposer\hfill Received: \receiveddate\nl
Seconded: \theseconders\hfill 
\ifx\reviseddate\relax                         
Accepted: \accepteddate                        
\else
Revised: \reviseddate                          
\fi}
\eject
}       

\let\maketitlepage\maketitlep
\let\maketitle\maketitlepage


\font\phead=cmsl9 scaled 950
\font\lhead=cmsl9 scaled 1050
\font\pnum=cmbx10 scaled 913
\font\lnum=cmbx10 
\font\pfoot=cmsl9 scaled 950
\font\lfoot=cmsl9 scaled 1050
\ifplaintex
\headline{\vbox to 0pt{\vskip -4.5mm\line{\small\phead\ifnum
\count0=\startpage ISSN 1364-0380 (on line)
1465-3060 (printed) \hfill {\pnum\folio}\else\ifodd\count0\def\\{ }%
\ifx\theshorttitle\relax\thetitle\else\theshorttitle\fi\hfill{\pnum\folio}
\else\def\\{ and }{\pnum\folio}\hfill\ifx\theshortauthors\relax\theauthors
\else\theshortauthors\fi\fi\fi}\vss}}
\footline{\vbox to 0pt{\vglue 0mm\line{\small\pfoot\ifnum\count0=\startpage
\copyright\ \gtp\hfill\else
\gt, Volume \thevolumenumber\ (\thevolumeyear)\hfill\fi}\vss
}}
\else
\makeatletter
\def\@oddhead{{\small\lhead\ifnum\count0=\startpage ISSN 1364-0380 (on line)
1465-3060 (printed) \hfill {\lnum\number\count0}\else\ifodd\count0
\def\\{ }\ifx\theshorttitle\relax \thetitle \else\theshorttitle\fi\hfill
{\lnum\number\count0}\else\def\\{ and }{\lnum\number\count0}
\hfill\ifx\theshortauthors\relax 
\theauthors\else\theshortauthors\fi\fi\fi}}\def\@evenhead{\@oddhead}
\def\@oddfoot{\small\lfoot\ifnum\count0=\startpage\copyright\ \gtp\hfill\else
\gt, Volume \thevolumenumber\ (\thevolumeyear)\hfill\fi}
\def\@evenfoot{\@oddfoot}
\makeatother
\fi


\newwrite\gtoutfile
\long\gdef\makeheadfile{  
{\def\\{, }\def\s{ }
\immediate\openout\gtoutfile head.xxx
\immediate\write\gtoutfile{To: math@arxiv.org}
\immediate\write\gtoutfile{Subject: put or rep NNNNN:pppp}
\immediate\write\gtoutfile{--text follows this line--}
\immediate\write\gtoutfile{Proxy-for: \ifx\theasciiauthors\relax
\theauthors\else\theasciiauthors\fi\s<\ifx\theasciiemail\relax\theemail\else\theasciiemail\fi>}
\immediate\write\gtoutfile{\noexpand\\}
\immediate\write\gtoutfile{Authors: \ifx\theasciiauthors\relax
\theauthors\else\theasciiauthors\fi}
{\def\\{ }\immediate\write\gtoutfile{Title: \ifx\theasciititle\relax
\thetitle\else\theasciititle\fi}}
\immediate\write\gtoutfile{Subj-class: GT or SG or MG etc}
\immediate\write\gtoutfile{MSC-class: \theprimaryclass\ifx\thesecondaryclass\relax\else, \thesecondaryclass\fi}
\immediate\write\gtoutfile{Journal-ref: Geom. Topol. \thevolumenumber
(\thevolumeyear) \startpage-\finishpage}
\immediate\write\gtoutfile{Comments: Published by Geometry and Topology at}
\immediate\write\gtoutfile{\s\s http://www.maths.warwick.ac.uk/gt/GTVol\thevolumenumber/paper\thepapernumber.abs.html}
\immediate\write\gtoutfile{\noexpand\\}
\immediate\write\gtoutfile{}
\ifx\theasciiabstract\relax
\immediate\write\gtoutfile{\theabstract}\else
\immediate\write\gtoutfile{\theasciiabstract}\fi
\immediate\write\gtoutfile{}
\immediate\write\gtoutfile{\noexpand\\}
\immediate\write\gtoutfile{}
\immediate\closeout\gtoutfile}}  

\def\maketitlepage{\maketitlep\makeheadfile}
\let\maketitle\maketitlepage